\documentclass[11pt]{amsart} 

\usepackage{amssymb,amsmath,amscd,amsfonts,graphicx,latexsym,epsfig}
\usepackage[all]{xy}

\usepackage[all]{xy}
\newtheorem{teo}{Theorem}[section]
\newtheorem{lem}[teo]{Lemma}

\newtheorem{cor}[teo]{Corollary}

\newtheorem{prop}[teo]{Proposition}
\newtheorem{defi}[teo]{Definition}

\newtheorem{remark}[teo]{Remark}

\newcommand{\Dd}{{\mathcal D}}

\newcommand{\Kk}{{\mathcal K}}
\newcommand{\Ll}{{\mathcal L}}

\newcommand{\Pp}{{\mathcal P}}

\newcommand{\Ss}{{\mathcal S}}

\newcommand{\I}{{\mathcal I}}
\newcommand{\Z}{{\mathbb Z}}

\newcommand{\X}{{\mathbb X}}

\newcommand{\N}{{\mathbb N}}
\DeclareMathOperator{\cut}{cut}
\DeclareMathOperator{\rk}{rk}
\DeclareMathOperator{\crk}{crk}

\newcommand{\compl}{\mathsf{C}}

\DeclareMathOperator{\lk}{lk}

\theoremstyle{definition}
\newtheorem{defn}[teo]{Definition}

\newcommand{\matN}{\ensuremath {\mathbb{N}}}

\newcommand{\calS} {\ensuremath {\mathcal{S}}}

\begin{document}

\title{Levels of knotting of spatial handlebodies}

\author{R. Benedetti, R. Frigerio}

\email{benedett@dm.unipi.it, frigerio@dm.unipi.it}

\address{Dipartimento di Matematica \\
Universit\`a di Pisa \\
Largo B.~Pontecorvo 5 \\
56127 Pisa, Italy}



\begin{abstract} Given a (genus 2) cube--with--holes $M=S^3\setminus H$, where $H$ is a handlebody,
  we relate intrinsic properties of $M$ (like its cut number) with extrinsic features depending on the way the handlebody
  $H$ is knotted in $S^3$. Starting from a first level of knotting that requires
  the non--existence of a planar spine for $H$, we define several instances of
  knotting of $H$ in terms of the non--existence of 
  spines with special properties. Some of these instances are implied by an intrinsic
  counterpart in terms of the non--existence of special cut--systems
  for $M$. We study a natural partial order on these instances of
  knotting, as well as its intrinsic counterpart, and the relations
  between them. To this aim, we recognize a few effective
  ``obstructions'' based on recent quandle--coloring invariants
  for spatial handlebodies, on the extension to the context of spatial
handlebodies of tools coming from the theory of homology boundary links, on the
analysis of appropriate coverings of $M$, and on the very classical use of Alexander
  elementary ideals of the fundamental group of $M$.
Our treatment of the matter
also allows us to revisit a few
  old--fashioned beautiful themes of 3--dimensional geometric
  topology.
\end{abstract}

\maketitle

\section{Introduction}\label{Intro}

Let $M$ be a compact connected 3--dimensional proper submanifold of
$S^3$ with smooth or PL boundary $\partial M$. We also assume that no boundary component of $M$ is spherical. Sometimes the
interior of $M$ is called a {\it domain} of $S^3$.

\smallskip

The {\it cut number} of $M$, henceforth denoted by 
$\cut (M)$, is the largest number 
of disjoint connected
two--sided properly embedded surfaces whose union does not disconnect $M$
(see~\cite{Stallings}).
If we denote by $\chi$ the Euler characteristic,  it is not difficult to show
that
$$1\leq \cut(M) \leq \rk H_1 (M)= n-\chi(M)\, , $$
 where $\rk H_1 (M)$ is the rank of the first homology group
of $M$ with integer coefficients, and
$n$ is the number of boundary components of $M$ (note that 
$\chi(M)=\chi(\partial M)/2$). From now on, homology and cohomology are always 
understood to have integer coefficients.

The {\it corank} of a group $G$, henceforth denoted by $\crk G$, 
is the largest rank of a free group
isomorphic to a quotient of $G$. It is known~\cite{Stallings} 
that the cut number of a manifold
coincides with the corank of its fundamental group 
(see \emph{e.g.}~the expository paper \cite{RRR}, where we have
also discussed in detail some basic properties of the 
cut number of domains, including the inequalities mentioned above -- we summarize
in Subsections~\ref{cut} and~\ref{cut-corank} the results that are relevant to our purposes).
We renewed our interest in these classical invariants of manifolds and 
groups after having
pointed out in~\cite{RRR} that the old {\it Helmholtz's cuts method}, used to
implement the Hodge decomposition of vector fields on spatial domains
(see \cite{Cant}), actually can be applied when the domains
have maximal cut number (\emph{i.e.}~they have a fundamental group
of maximal corank).

 \smallskip

 When $M=\compl (L)$, that is the closure of the complement of a
 tubular neighbourhood of an $n$--component link $L$, then by
 definition $L$ is a {\it homology boundary link}~\cite{Smythe} if $\cut(\compl (L))$ 
(or, equivalently, $\crk \pi_1 (\compl(L))$) is maximal, \emph{i.e.}~equal to $n$. 
In that case, the $n$ disjoint surfaces whose complement is not disconnected
are call \emph{generalized Seifert surfaces of $L$}.  
Homology boundary links form a distinguished
 class of links arisen at least since the discovering of Milnor's link
 invariants, and which have been widely studied in classical knot
 theory.

\smallskip

In this paper we are mostly interested to the somehow
opposite situation when the boundary is connected and consists of a
surface $S$ of genus $g\geq 2$. In this case we have  $1\leq \cut(M) \leq g $.

By Fox's {\it reimbedding theorem} \cite{Fox}, up to reimbedding $M$ 
can be realized as a ``cube with holes''  
$$M = \compl (H)\, , $$ 
that is the closure of the complement of a genus $g$ handlebody $H$
embedded in $S^3$ with $\partial H=\partial M$.  In such a situation it
is natural to look for relations between {\it intrinsic} properties of
$M$ (like its cut number) and {\it extrinsic} properties depending on the
way $H$ is knotted in $S^3$.  
\smallskip

Although several considerations can be generalized to an arbitrary
genus (see Section~\ref{where} at the end of the paper), in this paper we 
concentrate our attention on the case $g=2$ which already displays a rich phenomenology, as
well as some peculiar features. So, from now on, {\it both $H$ and $S=\partial H=\partial M$
  will be assumed of genus 2}.  \smallskip

Fox's reimbedding Theorem implies that one can study
the topology of spatial domains using techinques coming from knot theory.  More
precisely, a given (genus 2) spatial handlebody $H$ (considered up to
isotopy) is the regular neighborhood of infinitely many {\it handcuff
  spines} $\Gamma$, and every such spine carries a 2--component
{\it constituent link} $L_\Gamma$ (see Figure~\ref{planargraphs}). 
The handlebody $H$ is \emph{knotted} if it does not admit any planar spine.
Starting from this
first level of knotting, 
we define, in quite a natural way, several
{\it instances of knotting} of $H$ in terms of the {\it
  non}--existence of any spine $\Gamma$ enjoying less and less restrictive
properties (which sometimes can be expressed in terms of 
the constituent link $L_\Gamma$). The sets of spatial
handlebodies that realize the different instances of knotting are
partially ordered by inclusion (in the sense that smaller sets
correspond to higher levels of knotting).

Some of the
above instances of knotting are implied by intrinsic counterparts, 
which may be expressed in
terms of the non--existence of systems of disjoint properly embedded surfaces
satisfying 
less and less restrictive conditions (see Section~\ref{cut-levels}). This also leads to a certain partial
ordering relation.  Our aim is to clarify as much as possible these
extrinsic as well as intrinsic partial orders, and the
relations between them. 
As a key step in our investigation, we recognize
a few spatial handlebody invariants
that provide effective obstructions to be unknotted 
with respect to some of our instances of knotting.
  
\section{Basic results on spatial handlebodies}\label{ingredients}  
\begin{figure}[htbp]
\begin{center}
 \includegraphics[height=3cm]{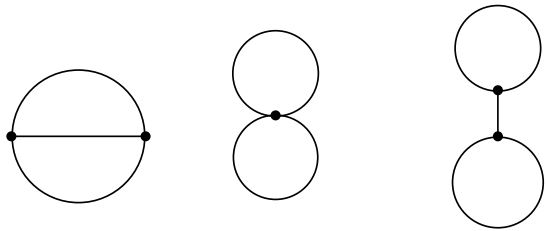}
\caption{\label{planargraphs} 
From the left
to the right: the $\theta$--, (f8)-- and (hc)--graphs.}
\end{center}
\end{figure}

Consider the planar realizations of the basic $\theta$--,
``figure-eight'' (f8)-- and ``handcuff'' (hc)--graphs
shown in Figure \ref{planargraphs}.

By a {\it (genus $2$) spatial graph} we mean a tame embedding
$\Gamma$ of any such graph in $S^3$. A {\it (genus $2$) spatial
handlebody} $H=H(\Gamma)$ is by definition a closed regular
neighbourhood of a spatial graph $\Gamma$, which is called a {\it
spine} of $H$. The symbol
$\compl (H)$ denotes the closure of the complementary domain.
Sometimes we write also $\compl(\Gamma)$, if $H=H(\Gamma)$.

\begin{figure}[htbp]
\begin{center}
 \includegraphics[height=2cm]{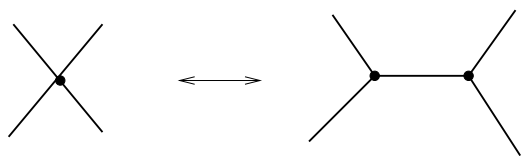}
\caption{\label{Wmove} Whitehead move.}
\end{center}
\end{figure}

Spatial handlebodies are considered up to 
isotopy. Two spatial graphs are the spines of isotopic spatial handlebodies
if and only if they are equivalent according to
the equivalence relation
generated by ambient isotopy together with
the {\it Whitehead move} shown in Figure \ref{Wmove}.
 
\begin{figure}[htbp]
\begin{center}
 \includegraphics[height=8cm]{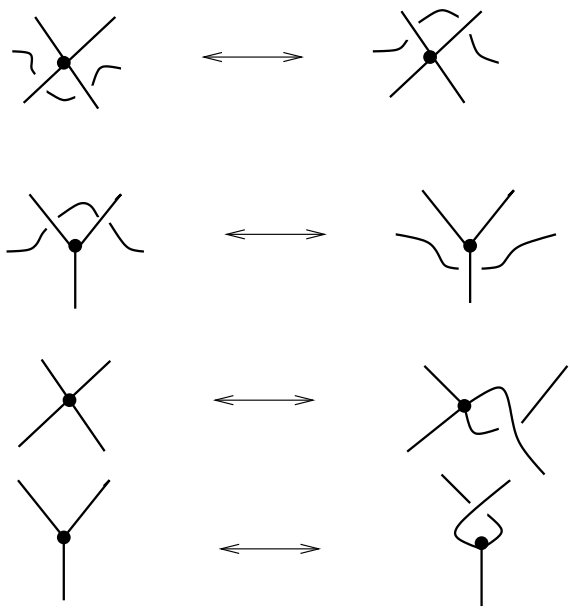}
\caption{\label{isomove} Diagram isotopy moves.}
\end{center}
\end{figure}

A spine is {\it unknotted} if it is isotopic to one of the above
planar graph. A handlebody $H$ is {\it unknotted} if it admits an
unknotted spine.  Sometimes it is useful to deal with planar diagrams
associated to generic projections, that encode the spines, and diagram
moves that recover the isotopy of either spatial graphs or
handlebodies. The graph isotopy moves are the usual local {\it
  Reidemeister moves}, like for link diagrams, and a few additional
moves at vertices (see \cite{kauffman}).  Some of these {\it vertex
  moves} are shown in Figure \ref{isomove}; to get the full set it is
enough to change simultaneously the (over/under) crossings in each
move of the picture. In order to recover the handlebody isotopy in
terms of spine diagrams, it is enough to add the above Whitehead move,
interpreted now as acting on diagrams as well (see \cite{ishii1}). We
call this whole set of moves {\it spine moves}.

\subsection{Constituent knots and links}\label{costituentK} 
A {\it constituent knot} of any spatial graph $\Gamma$ as above is a
(spatial) subgraph of $\Gamma$ homeomorphic to $S^1$.  

If $\Gamma$ is a (f8)--graph, then the two constituent knots $K_1$ and
$K_2$ intersect exactly at the unique singular point (the vertex) $p$
of $\Gamma$; in this case we write $\Gamma = K_1 \bigvee_p K_2$.

If $\Gamma$ is a (hc)--graph, then the two constituent knots are disjoint,
hence they form a constituent link $L_\Gamma $. This is obtained by
removing from $\Gamma$ the interior of its {\it isthmus}, \emph{i.e.}~the edge
that connects the two knots. 

Let us denote by $\Kk(\Gamma)$ the knot isotopy classes realized by
the constituent knots of $\Gamma$. If $\Gamma$ is either a (f8)-- or a (hc)--graph,
then $|\Kk(\Gamma)|\leq 2$; if it is a $\theta$--graph, then
$|\Kk(\Gamma)|\leq 3$.

If $\tilde \Gamma \to \Gamma$ is a Whitehead move producing a
(f8)-graph $\Gamma$, then $\Kk(\Gamma) \subset \Kk(\tilde \Gamma)$.

The following Lemma is immediate but important.

\begin{lem}\label{triv-inv} 
  The following sets of isotopy classes of spatial graphs and links
  respectively are isotopy invariants of the handlebody $H$:
\begin{itemize}
\item the set $\Ss(H)$ of isotopy classes of (hc)--spines of $H$;

\item the set $\Ll(H)$ of isotopy classes of links $L_\Gamma$,  where
$\Gamma$ varies in $\Ss(H)$.
\end{itemize}
\end{lem}
\smallskip

\begin{figure}[htbp]
\begin{center}
 \includegraphics[height=5cm]{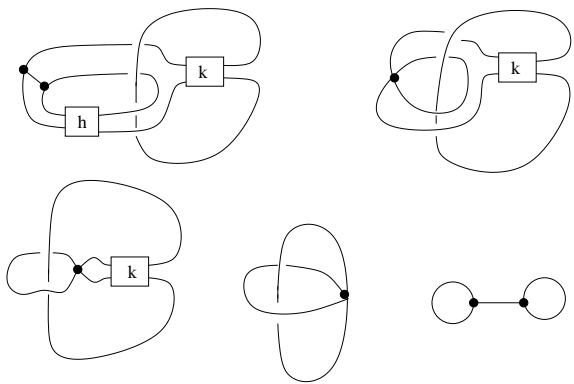}
\caption{\label{handcuff2} Spines of the unknotted handlebody.}
\end{center}
\end{figure}

\begin{figure}[htbp]
\begin{center}
 \includegraphics[height=4cm]{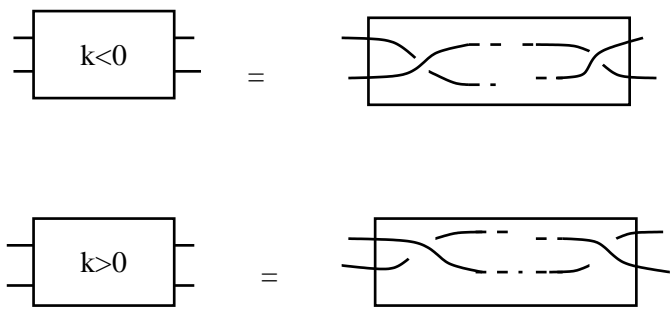}
\caption{\label{Boxes} The boxes used in Figure~\ref{handcuff2}.}
\end{center}
\end{figure}

\begin{remark}\label{manyspine}{\rm
\noindent (1) 
We will often confuse links and graphs with the respective isotopy
classes. The observations before Lemma~\ref{triv-inv} imply that every
constituent knot of any (not necessarily handcuff) spine of $H$
arises as a component of some $L_\Gamma \in \Ll(H)$. 

(2) 
For a single graph $\Gamma$, its finite set of constituent knots or
links is a rather informative, largely used invariant, to which one
can apply all the formidable machinery of classical knot theory. A
complication in the case of handlebodies arises from the fact that 
spines and constituent knots or links are infinitely
many. Moreover, there is not an immediate relationship between the
knotting of a given spine and the knotting of the
associated handlebody.  Sometimes this sounds a bit
anti--intuitive. For example, essentially by definition every {\it
  tunnel number $1$} knot (resp. link) arises as a component of some $L_\Gamma$ 
  (resp.~as some $L_\Gamma$) of
an unknotted (genus 2) handlebody. These knots form a richly structured
family of knots (see for instance~\cite{CMc}). It is not difficult to construct 
tunnel number $1$ links with knotted components (see \emph{e.g.}~Figure~\ref{tunnel},
where we describe an example taken from~\cite{lee}).

\begin{figure}[htbp]
\begin{center}
 \includegraphics[height=6cm]{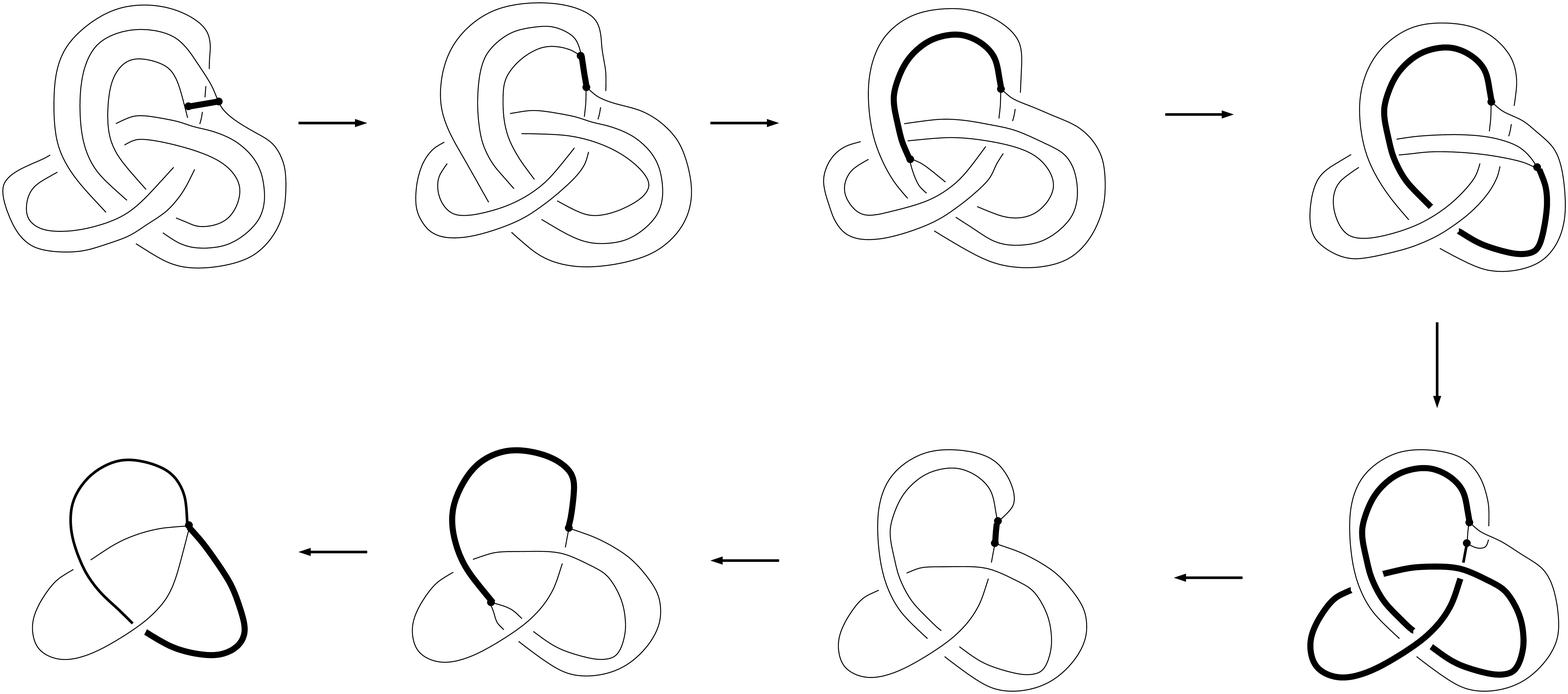}
\caption{\label{tunnel} The diagram on the left of the top row shows an (hc)--graph $\Gamma$ whose constituent
knots are both knotted. The sequence of moves shows that $\Gamma$ is a spine
of the
unknotted handlebody.}
\end{center}
\end{figure}

In the same spirit, 
by
varying the integers $k$ and $h$ in Figure~\ref{handcuff2} (the meaning
of the boxes is defined in Figure \ref{Boxes}), we get
infinitely many non--isotopic (hc)--spines $\Gamma$ of an unknotted
handlebody, with arbitrary linking number of the components of
$L_\Gamma$. See also Remark \ref{knot-vs-unknot} below.}
\end{remark}

\section{Instances of knotting, intrinsic and extrinsic knotting:\\ statement of the problems and main results}\label{levels}

We are going to use the invariants of $H$ introduced in Lemma~\ref{triv-inv}
above, in order to define certain {\it instances of knotting} of
$H$. We will distinguish the instances defined in terms of the set $\mathcal{S}(H)$ of 
the handcuff spines of $H$, from
the ones defined only in terms of the set $\mathcal{L}(H)$ of the constituent links.
From now on, unless otherwise stated every spine is understood to be of handcuff type.

\subsection{Spine--defined instances of knotting} 
 Let us recall that a $2$-component link $L=K_1\cup K_2$ is a \emph{boundary link}
if $K_1,K_2$ admit \emph{disjoint} Seifert surfaces $S_1,S_2$. It is readily seen
that such surfaces, if they exist, do not disconnect $\compl (L)$.
Therefore, a boundary link is in particular a homology boundary link.
 Let $H$ be a genus 2 handlebody in $S^3$.  
 
  \begin{itemize} 
  \item $H$ is {\it
      $(1)_S$-knotted} if it does not admit any planar spine, \emph{i.e.}~if it is knotted
      in the usual sense.
  \item A spine $\Gamma$ of $H$ is a \emph{split spine} if there exists
    an embedded 2--sphere $\Sigma \ $ in $S^3$ that intersects
    $\Gamma$ transversely at just one regular point of its isthmus (in
    particular, $L_\Gamma$ is a split link). A
    handlebody $H$ is {\it $(2)_S$-knotted} if it does not admit any
      split spine.
\item A  spine $\Gamma$ of $H$ is a {\it boundary spine} if its constituent link $L_\Gamma$ is 
a boundary link that admits a pair of disjoint Seifert 
surfaces whose interiors are contained in
$S^3\setminus \Gamma$. A handlebody $H$ is {\it $(3)_S$-knotted} if it does
not admit any boundary spine.
\item
A  spine $\Gamma$ of $H$ is a {\it homology boundary spine} if its constituent link $L_\Gamma$ is 
a homology boundary link that admits a pair of disjoint generalized Seifert 
surfaces whose interiors are contained in
$S^3\setminus \Gamma$. A handlebody $H$ is {\it $(4)_S$-knotted} if it does
not admit any homology boundary spine.
\end{itemize}

\subsection{Link--defined instances of knotting.}
Let $H$ be a genus 2 handlebody in $S^3$. 

\begin{itemize}
\item $H$ is {\it $(1)_L$-knotted} if there is not any {\it
trivial} link $L_\Gamma \in \Ll(H)$.
\item $H$ is {\it $(2)_L$-knotted} if there is not any {\it
split link} $L_\Gamma \in \Ll(H)$.
\item $H$ is {\it $(3)_L$-knotted} if there is not any {\it
boundary link} $L_\Gamma \in \Ll(H)$.
\item $H$ is {\it $(4)_L$-knotted} if there is not any {\it
homology boundary link} $L_\Gamma \in \Ll(H)$.
\end{itemize}

\subsection{The general structure of knotting levels}
The following tautological result describes some obvious relations 
between the instances of knotting we have introduced.

\begin{prop}\label{tautological}
For every $k$ such that the following statements make sense, we have:
$$
\begin{array}{ccc}
H\ {\rm is}\ (k+1)_S-{\rm knotted}\ & \Longrightarrow &\ H\ {\rm is}\ (k)_S-{\rm knotted}\vspace{.1cm}\\
H\ {\rm is}\ (k+1)_L-{\rm knotted}\ & \Longrightarrow &\ H\ {\rm is}\ (k)_L-{\rm knotted}\vspace{.1cm}\\
H\ {\rm is}\ (k)_L-{\rm knotted}\ & \Longrightarrow &\ H\ {\rm is}\ (k)_S-{\rm knotted}\ .
\end{array}
$$
\end{prop}

Hence, we have a few five--steps {\it staircases} of increasing
levels of knotting going up from the $(1)_S$-level to the
$(4)_L$-level, each one starting with a substring of spine--defined
and ending with a substring of link--defined instances. To study the
partial ordering  mentioned in the Introduction, we should in particular: 
\smallskip

\begin{enumerate}
\item[(a)]
Show that every instance of knotting is non--empty. 
%\smallskip
\item[(b)]
Clarify to which extent such staircases are {\it strictly}
increasing.
%\smallskip
\item[(c)]
Study whether or not one can establish
any implication between the instances
$(k)_L$ and $(k)_S$.
\end{enumerate}

\begin{remark}\label{link-vs-spine}
{\rm
  The above $(*)_L$-knotting conditions describe increasing
  levels of knotting of 2--component links. The fact that these conditions determine
  {\it strictly} decreasing sets of links is a classical result, which
  cannot be immediately translated into the context of handlebodies (for example, because of the
  considerations in Remark \ref{manyspine}).
}
\end{remark}

\subsection{Levels of knotting: main results}
The following Theorem summarizes our results about the relations between the 
instances of knottings we have introduced.

\begin{teo}\label{summarize1}
We have the following facts:
$$\xymatrix{H\ {\rm is}\ (1)_L-{\rm knotted} \ar@{=>}[r]|| & 
 H\ {\rm is}\ (2)_S-{\rm knotted} }
 $$\vspace{-1cm}\par
 $$\xymatrix
{H\ {\rm is}\ (2)_L-{\rm knotted} \ar@{=>}[r]|| & 
 H\ {\rm is}\ (3)_S-{\rm knotted} }
 $$\vspace{-.9cm}\par
$$\xymatrix{H\ {\rm is}\ (3)_L-{\rm knotted} \ar@{=>}[r]|| & 
 H\ {\rm is}\ (4)_L-{\rm knotted} }
 $$\vspace{-.9cm}\par
$$\xymatrix{H\ {\rm is}\ (4)_L-{\rm knotted} \ar@{<=>}[r] & 
 H\ {\rm is}\ (4)_S-{\rm knotted} }
 $$\vspace{-.9cm}\par
$$\xymatrix{H\ {\rm is}\ (3)_S-{\rm knotted} \ar@{=>}[r]|| & 
 H\ {\rm is}\ (1)_L-{\rm knotted} }
 $$
\end{teo}

The following diagram summarizes the results described in Proposition~\ref{tautological}
and Theorem~\ref{summarize1}, that completely characterize
the relations among the levels of knotting we have introduced (see Corollary~\ref{distinct:cor}).
We have labelled every non-tautological arrow
by a reference to the Proposition where the corresponding implication is proved.

\smallskip

$$
\xymatrix{
*+[F]{(4)_L-{\rm knotting}} \ar[dd]  \ar[rrr]& {} & {} &  *+[F-:<3pt>]{(4)_S-{\rm knotting}}\ar@/_.8pc/[lll]_{\rm Prop.~\ref{L4S4}}\ar[dd]\\ & & & \\
*+[F]{(3)_L-{\rm knotting}} \ar[rrr]\ar[dd] \ar@/^1pc/[uu]|=^{\rm Prop.~\ref{nuovaprop}\ } 
&  &  & *+[F-:<3pt>]{(3)_S-{\rm knotting}}\ar[dd] 
\ar `r/20pt[d] `[ddddd]|=^{\ \rm Prop.~\ref{alternative}} [llldddd] \\ & & & \\ 
*+[F]{(2)_L-{\rm knotting}} \ar@/^1pc/[uu]|=^{\rm Prop.~\ref{2Lvs3L}\ } \ar[rrr]\ar@/_.5pc/[uurrr]|\|^{\rm Prop.~\ref{H4}}\ar[dd] 
 &  & & *+[F-:<3pt>]{(2)_S-{\rm knotting}}\ar[dd]
\ar@/_.7pc/[ddlll]|\|_{\rm Prop.~\ref{1 knotted}\ }
\\ & & &\\
*+[F]{(1)_L-{\rm knotting}}  \ar[rrr]\ar@/_.7pc/[uurrr]|\|_{\rm\quad\ Prop.~\ref{1Lvs1S:lem}}
 &  & & *+[F-:<3pt>]{(1)_S-{\rm knotting}}\ar@/_1pc/[uu]|=_{\rm \, Prop.~\ref{K1}}\\
& & &
}
$$

\smallskip

In Theorem~\ref{summarize1} we did not mention the results proved in
Propositions~\ref{K1}, \ref{2Lvs3L} and \ref{1 knotted} because,
once the tautological statements of Proposition~\ref{tautological} are established, 
Proposition~\ref{K1} is a consequence of Proposition~\ref{1Lvs1S:lem},
Proposition~\ref{2Lvs3L} is a consequence of Proposition~\ref{H4}, and
Proposition~\ref{1 knotted} is a consequence of Proposition~\ref{alternative}.
\smallskip

However, it is maybe worth mentioning that the arguments proving 
Propositions~\ref{K1}, \ref{2Lvs3L} and \ref{1 knotted} are independent, and sometimes
quite different in nature, from the ones proving Propositions~\ref{1Lvs1S:lem},
\ref{H4}, and
\ref{alternative}. For example, the proof of Proposition~\ref{alternative} 
exploits the notion of \emph{handlebody pattern} introduced in Subsection~\ref{handlepat}, while 
the proof of Proposition~\ref{1 knotted} relies on the use
of quandle invariants. 
We also give a different proof of Proposition~\ref{alternative}
in Proposition~\ref{via-A-obs}, where we exploit  
suitable obstructions coming from Alexander invariants. 
\smallskip

In our opinion, a good reason for providing more proofs of some of our results is that
the interplay of different techniques has been of help
to us for
figuring out 
the global picture of the space of knotted handlebodies. 

\medskip

Let us take two distinct instances of knottings $(k)_J$ and $(k')_{J'}$, where
$J,J'\in \{L,S\}$ and $k,k'\in\{1,2,3,4\}$. Then, 
putting together Proposition~\ref{tautological} and Theorem~\ref{summarize1}, 
one easily gets the following:

\begin{cor}\label{distinct:cor}
With the exception of the case $k_J=(4)_L$ and $(k')_{J'}=(4)_S$
(or viceversa),
we have that
\begin{center}
$(k)_J$-knotting $\Longrightarrow$ $(k')_{J'}$-knotting
\end{center}
\noindent
if and only if $(k')_{J'}$-knotting descends from $(k)_{J}$-knotting
via a ``staircase'' chain of tautological implications (see Proposition~\ref{tautological}). 
\end{cor}

By Theorem~\ref{summarize1}, the class of $(4)_L$-knotted handlebodies coincides with the class of $(4)_S$-knotted handlebodies.
Corollary~\ref{distinct:cor} tells us that, with this exception,
the different instances of knotting we have introduced indeed describe
distinct classes of spatial handlebodies.

\smallskip

In Theorem~\ref{nonempty} below we also show that every instance of knotting is indeed non-empty,
\emph{i.e.}~that there exist examples of $(4)_L$-knotted (hence, $(k)_J$-knotted for every $k=1,\ldots,4$, $J\in \{L,S\}$) spatial handlebodies.

\subsection{First remarks on unknotted and maximally knotted handlebodies}
A $(1)_S$-knotted $H$ is knotted in the ordinary sense; recall that by
Waldhausen's Theorem \cite{WALD, scharlemann2} a
spatial handlebody $H$ is unknotted if and
only if also $\compl (H)$ is a handlebody, \emph{i.e.}~if and only if 
$H$ and $\compl(H)$ form a (genus
2) {\it Heegaard splitting of $S^3$} (which is
unique up to isotopy).  
\smallskip

On the opposite side, we have the following:
\begin{lem}\label{corank1} If $\cut(\compl (H))$ is not
  maximal (\emph{i.e.}~not equal to 2), then $H$ is $(4)_L$-knotted.
\end{lem}
\begin{proof}
If $H$ is not $(4)_L$-knotted, then by definition there exists a link $L_\Gamma
\in \Ll(H)$ which is homology boundary, so $\crk\pi_1(\compl(L))=2$. 
By Lemma~\ref{onto-i*} below, this implies that $\crk \pi_1(\compl(H))=2$, so
$\cut(\compl (H)) = 2$. 
\end{proof}

\begin{lem}\label{onto-i*} 
Let $L_\Gamma \in \Ll(H)$. Up to isotopy, we can assume that $\compl
(H) \subset \compl (L_\Gamma)$. Then $i_*: \pi_1(\compl (H)) \to
\pi_1(\compl (L_\Gamma))$ is an epimorphism. In particular, 
$\crk \pi_1(\compl(H))\geq \crk\pi_1 (\compl(L_\Gamma))$.
\end{lem}
\begin{proof}  
 By a general position argument it is easy to see that every loop
in $S^3\setminus L_\Gamma$ is homotopic to a loop that does not
intersect the isthmus of $\Gamma$; this easily implies that $i_*$
is onto.
\end{proof}

\subsection{Cut systems}\label{cut}
In this Subsection we provide a brief  discussion of the notion of cut system, which will also allow us to 
sketch a proof of the well-known equality between the cut number of a manifold and the corank of its fundamental group.
We focus on
the special situation we are interested in, addressing the reader
\emph{e.g.}~to \cite{RRR} for a general and detailed discussion and the proofs. 

\smallskip

Let $M$ be equal either to $\compl(H)$ (where $H$ is a genus--2 spatial handlebody) or to $\compl(L)$ ($L$
being a 2--component link).

\begin{defi}\label{cut-s} {\rm A {\it cut system} $\Ss=\{S_1, S_2\}$ of
    $M$ is a pair of {\it disjoint} connected oriented surfaces 
    properly embedded in $M$
    such
    that $M \setminus (S_1\cup  S_2)$ is connected.}
\end{defi}

We now list some classical well--known properties and characterizations
of cut systems, also giving a hint about the proof of the equality $\cut(M)=\crk\pi_1 (M)$
already mentioned in the Introduction. 
%For a detailed account  see \emph{e.g.}~the expository paper~\cite{RRR}.
Let $\Ss=(S_1,S_2,\ldots,S_k)$ be a set of disjoint connected oriented surfaces 
properly embedded in $M$. Then $M\setminus(S_1\cup\ldots\cup S_k)$ is connected
if and only if $S_1,\ldots,S_k$ define linearly indepenent elements of $H_2(M,\partial M)\cong \mathbb{Z}^2$. In particular,
if this is the case, then necessarily $k\leq 2$ (whence $\cut (M)\leq 2$), and 
one can see
that 
$\Ss$ is a cut system (\emph{i.e.}~$k=2$)
if and only if 
$S_1$ and $S_2$ actually define 
a geometric basis of $H_2(M,\partial M)$.  In this case, 
by Alexander duality, 
$\partial
\Ss = \{\partial S_1, \partial S_2\}$ provides a geometric basis of
$\ker i_*$, where $i$ is the inclusion of the boundary $\partial M$ into
$M$ (in particular, $\partial S_i\neq\emptyset$ for $i=1,2$). With a slight abuse, ``$\partial$'' denotes 
here both the boundary
homomorphism in the homology long exact sequence of the pair
$(M,\partial M)$, and the geometric boundary of $\Ss$.  

By definition, a link is homology boundary if and only if its complement admits a cut system.
Something more can be achieved when considering classical boundary links.
If $M$ is the complement of a boundary link $L$, then it is easily seen
that every pair of disjoint Seifert surfaces for the components of $L$ define
a basis of $H_2(M,\partial M)$, whence a cut system for $M$. Such a cut system
enjoys the nice property that the boundaries of its surfaces are connected.

\subsection{Cut number and corank}\label{cut-corank}
Since
 $\pi_1(M)/[\pi_1 (M),\pi_1 (M)]\cong H_1 (M)\cong \mathbb{Z}^2$, it is easily
 seen that $1\leq \crk \pi_1 (M)\leq 2$. In fact,
one can prove that $\crk \pi_1 (M)=2$ if and only if $M$ maps onto the bouquet
$S^1\vee S^1$ of two circles via a map which induces a surjective homomorphism between the fundamental
groups. An elementary instance of the Pontryagin--Thom construction
may now be used to show that the required $\pi_1$--surjective map between $M$ and $S^1\vee S^1$ 
exists if and only if 
$M$ admits a cut system.
Theferore, $M$ admits a cut system if and only if $\crk\pi_1 (M)=2$. 

The argument just sketched may be easily generalized for showing
that the equality $\cut(N)=\crk\pi_1 (N)$ holds in fact for every compact
$3$--manifold $N$ (and actually for every $n$--manifold, $n\in\mathbb{N}$, once
one has replaced ``surfaces"  with ``hypersurfaces" in the definition of cut number).

%We have:

%\begin{prop}\label{corank-char} The following are
%equivalent facts: 

%(1) $\crk \pi_1(M)= 2 \ .$

%(2) $M$ admits a cut system.
%\end{prop} 

%The equivalence between (1) and (2) (as well the remarks after the
%above Definition) are very classical well-known facts.  Roughly, first
%one shows that (1) holds if and only if there is a map $r$ from $M$
%{\it onto} a bouquet of 2 circles. Then one achieves the equivalence
%between (1) and (2) by applying to this map $r$ an elementary instance
%of the Pontryagin-Thom construction.  
%\medskip

\subsection{Cut systems and levels of knotting}\label{cut-levels} 
Note that the two surfaces of a cut system of $M$ (if any), have not
necessarily connected boundary, and when $M=\compl(H)$, there could be
also components of $\partial \Ss$ that separate $S=\partial M$.  This
suggests the existence of different intrinsic ``levels of
complication'' of $M=\compl(H)$, which are defined in terms of the
non--existence of cut systems having boundaries that satisfy some
special conditions.

\begin{defi}\label{reduced-bound} 
{\rm 

  Let $\Ss=(S_1, S_2)$ be a cut system of $M=\compl(H)$. The {\it
    reduced boundary} $\partial_R S_j$ of $S_j$ is made by the
  boundary components of $S_j$ that do not separate $S=\partial M$. A cut
  system $\Ss$ is said to be {\it $\partial$-connected} (resp. {\it
    $\partial_R$-connected}) if both surfaces $S_j$ have connected
  boundary (resp. connected reduced boundary).  }

\end{defi}
\smallskip

We have proved in~Lemma \ref{corank1} that if $M$ has no cut systems,
then $H$ is $(4)_L$-knotted. The following Lemma provides more relations
between the (non)existence of special cut systems for $M$ and the knotting level
of $H$.

\begin{lem}\label{ex=>in} Let $M=\compl (H)$. Then:
%\smallskip
\begin{enumerate}
\item[(a)] If $M$ has no $\partial$-connected cut systems, then $H$ is
$(3)_S$-knotted. 
\item[(b)] If $M$ has no $\partial_R$-connected cut systems, then
$H$ is $(3)_L$-knotted. 
\end{enumerate}
\end{lem}
\begin{proof} 
(a): Suppose $H$ is $(3)_S$--unknotted, and let $\Gamma$ be a boundary spine
for $H$. Then, the boundary link $L_\Gamma$ admits 
disjoint Seifert surfaces that do not
intersect the interior of the isthmus, thus defining also a cut
system of $M$ with connected boundaries. 
\smallskip

(b): Suppose $H$ is $(3)_L$--unknotted, and let $\Gamma$ be a spine
for $H$ with constituent boundary link $L_\Gamma$.
Up to isotopy,
we
can assume that $L_\Gamma$ admits a pair of disjoint Seifert surfaces
that are transverse to the interior of the
isthmus. Then the intersection of these surfaces with $M$ form a cut
system with connected reduced boundaries. 
\end{proof}

\subsection{Boundary preserving maps onto handlebodies}
Let $W$ be a genus 2 handlebody.
Following Lambert~\cite{lambert}, we 
say that a continuous map $\varphi\colon M\to W$ is a
\emph{$(M \to W)$--boundary--preserving--map} 
if  
$\phi_{|S}$ is a homeomorphism of $S=\partial H$ onto $\partial W$
(such a $\phi$ is necessarily surjective).
\smallskip

The following result is proved in \cite[Theorem 2]{lambert}.

\begin{prop}\label{lambert}  Let $M=\compl(H)$. Then the
following facts are equivalent:
\begin{enumerate}
\item
%(1) 
$M$ admits a $\partial$-connected cut system. 
\item
%(2) 
There exists a {\rm $(M \to W)$--boundary--preserving--map}. 
\end{enumerate}
\end{prop}

In Sections~\ref{more-bpm} and~\ref{Alex-obs} we describe some obstructions that prevent
$M$ to admit a $\partial$-connected cut system. Such obstructions are obtained
from the study of handlebody patterns as defined in Subsection~\ref{handlepat},
from the analysis
of the maximal free covering of $M$, and 
from the study of the Alexander ideals of $M$.

\subsection{Intrinsic and extrinsic knotting levels: main results}
Let us underline that every property of $H$ which can be expressed only in terms
of  the existence of cut systems of $\compl (H)$ (possibly with specific properties)
does not depend on the embedding of $H$ into $S^3$, but
is an intrinsic property
of the complement of $H$. 
The question whether one can characterize the instances of knotting 
of $H$ via intrinsic properties of $\compl (H)$ seems to be very difficult. 
An easy result in this context is described in Proposition~\ref{1-K-char},
where it is shown that $H$ is $(2)_S$--knotted if and only if
$M=\compl (H)$ is boundary--irreducible, \emph{i.e.}~it has incompressible boundary.

More difficult questions are suggested by
Lemmas \ref{corank1} and \ref{ex=>in}:  
is it possible to reverse the implications proved there?
In Section~\ref{Ext-vs-intr}
we obtain in particular
the following interesting result:

\begin{teo}\label{ex<=>in:teo} $M=\compl(H)$ admits a 
  $\partial$-connected cut system if and only if $H$ is not
  $(3)_S$-knotted.
\end{teo}

The following result shows that 
the non--existence of cut systems satisfying \emph{strictly} decreasingly demanding conditions actually corresponds to
{\it strictly} increasing levels of (intrinsic) knotting.  

\begin{teo}\label{summarize3}
We have the following facts:
\begin{enumerate}
 \item There exist handlebody complements having incompressible boundary
and admitting a $\partial$-connected cut system.
\item
There exist handlebody complements which do not admit any $\partial$-connected cut system
but admit a $\partial_R$-connected cut system.
\item
There exist handlebody complements which do not admit any $\partial_R$-connected
cut system
but admit a generic cut system (\emph{i.e.}~they have cut number equal to 2).
\end{enumerate}
\end{teo}
\begin{proof}
By Proposition~\ref{1-K-char} and Theorem~\ref{ex<=>in:teo},
point~(1) is equivalent to the fact that $(2)_S$-knotting does not
imply $(3)_S$-knotting, which is a consequence of 
Corollary~\ref{distinct:cor} (more precisely, 
Proposition~\ref{H4} implies that there exist examples of
$(3)_S$-unknotted handlebodies which are $(2)_L$-knotted,
whence $(2)_S$-knotted).

By Proposition~\ref{alternative}, there exist $(1)_L$-unknotted handlebodies
whose complements do not admit any $\partial$-connected cut system. By Lemma~\ref{ex=>in},
every such complement admits a $\partial_R$-connected cut system, whence point~(2). 

Point~(3) is a consequence of Proposition~\ref{nuovaprop}. 
\end{proof}

Putting together Lemmas~\ref{corank1}, \ref{ex=>in} and Theorems~\ref{ex<=>in:teo}, \ref{summarize3} we can summarize 
the relations between the intrinsic and extrinsic levels of knotting as follows:

\smallskip 

$$
\xymatrix{
*+[F]{{\rm no\ cut\ systems}} \ar[d]  \ar@/_.4pc/[rr]& {} &   *+[F-:<3pt>]{(4)_L-{\rm knotting}}\ar@/_.4pc/@{-->}[ll]_{\rm ?}\ar[d]\\ 
*+[F]{{\rm no}\ \partial_R-{\rm connected\ cut\ system}} \ar@/_.3pc/[rr] 
\ar@/^1pc/[u]|=
\ar[d] &  &  *+[F-:<3pt>]{(3)_L-{\rm knotting}}\ar[d] \ar@/_.3pc/@{-->}[ll]_(.4){\rm ?}\\ 
*+[F]{{\rm no}\ \partial-{\rm connected\ cut\ system}}\ar[d] \ar@{<->}[rr] 
\ar@/^1pc/[u]|=
& & *+[F-:<3pt>]{(3)_S-{\rm knotting}}\ar[d]\\
*+[F]{{\rm incompressible\ boundary}} \ar@{<->}[rr] 
\ar@/^1pc/[u]|=
& & *+[F-:<3pt>]{(2)_S-{\rm knotting}}
}
$$

\vspace{.8cm}\par
 
Let us now turn to the question whether every instance of knotting is non-empty.
Of course, in order to give an affirmative answer to the question it is sufficient
to show that there exist $(4)_L$-knotted spatial handlebodies. 
As already mentioned in Lemma~\ref{corank1}, if the complement of a spatial
handlebody $H$ has cut number equal to 1, then $H$ is $(4)_L$-knotted.
Moreover, Jaco exhibited in~\cite{jaco} a spatial handlebody $H$ such that
$\cut( \compl(H))=1$, so every instance of knotting is non-empty. Following~\cite{suzuki},
in Section~\ref{Alex-obs} we show how Alexander invariants can be used to provide obstructions
for a handlebody complement to have maximal cut number. Moreover,
in Subsection~\ref{infinitelym} we refine Jaco's result and 
prove the following:

\begin{teo}\label{nonempty}
There exist infinitely many non-isotopic spatial handlebodies
$\{H_i\}_{i\in I}$ such that 
$\cut(\compl(H))=1$ for every $i\in I$. In particular,
there exist infinitely many non-isotopic $(4)_L$-knotted spatial handelbodies.
\end{teo}

\section{About  the lower instances of knotting}\label{low-instance}
Let $H$ be a spatial handlebody, and set
as usual $M=\compl(H)$ and $S=\partial H=\partial M$. 
Recall that we denote by $\mathcal{S}(H)$ and $\mathcal{L}(H)$ the sets of isotopy
classes respectively of (hc)--spines and constituent links of $H$.

\begin{prop}\label{K1} Let $L_\Gamma \in \Ll(H)$. Then:
\begin{enumerate}

\item If $L_\Gamma$ is a non--trivial homology boundary link, then
$H$ is $(1)_S$-knotted.

\item Assume that $H$ is unknotted. Then $L_\Gamma$ is a homology
boundary link if and only if the spine $\Gamma$ is unknotted (\emph{i.e.}~planar).

\item $(1)_S$-knotting does not imply  $(2)_S$-knotting.
\end{enumerate}
\end{prop}

\begin{proof}
 (1): Assume by contradiction that 
$H$ is unknotted. Then also $\compl (H)$ is
a genus 2 handlebody, hence $\pi_1(\compl (H))= \Z^{*2}$.
On the other hand, Papakyriakopolous' unknotting theorem for knots
\cite{papa} (which is based on his ``Loop Theorem'') generalize to links
 (see for instance Theorem 1.1 in
\cite{hillman}), and implies that a $n$-component link is trivial if and only
if the fundamental group of its complement is free on $n$ generators.
As a consequence,
if we
assume that $L_\Gamma$ is a non--trivial homology boundary link, 
by Lemma~\ref{onto-i*} 
we could realize $\Z^{*2}$ as a proper quotient of itself and this is
not possible because $\Z^{*2}$ is a Hopfian group \cite{MKS}.  

(2):
If $H$ is unknotted, then by (1) a homology boundary $L_\Gamma\in\Ll(H)$ is
necessarily trivial. Then we can apply the main theorem of
\cite{scharlemann1} and conclude that the spine $\Gamma$ itself is
unknotted. 

(3):
Take a  split (hc)--graph $\Gamma$ with associated non--trivial split link $L_\Gamma$, and 
let $H$ be a regular neighbourhood of $\Gamma$. 
By construction, $H$ is not $(2)_S$--knotted, while
it is $(1)_S$--knotted by (1).
\end{proof}

\begin{remark}\label{knot-vs-unknot}{\rm

    The above Lemma can be rephrased by saying that the existence of a
    very special (\emph{i.e.}~unknotted) spine forces all the non--planar spines of
    an unknotted handlebody to be at the highest level of knotting (\emph{i.e.}~to have corresponding links
    which are not homology boundary).} 
    \end{remark}

The following Proposition provides a characterization of
$(2)_S$-knotting. It turns out that in spite of its definition, $(2)_S$--knotting 
reflects an {\it intrinsic} property of $M$.
 
\begin{prop}\label{1-K-char}  The following facts are equivalent:
\begin{enumerate}
\item $H$ is $(2)_S$-knotted.
\item $M=\compl (H)$ has incompressible boundary (\emph{i.e.}~it is
$\partial$-irreducible).
\item $\pi_1(M)$ is not decomposable with respect to free
products.
\end{enumerate}
\end{prop}
\begin{proof}
The implication
(2) $\Rightarrow$ (3) is a consequence of the results proved in~\cite{jaco2}, and the implication
(3) $\Rightarrow$ (1) is trivial,
so it is sufficient to show that 
(1) $\Rightarrow$ (2). 
\smallskip

So,
we assume that $M$ is boundary--compressible, and show that
$H$ admits a split spine.
We follow the argument of
\cite[Theorem 4]{scharlemann4}, getting here a stronger conclusion
due to the fact that we are dealing with genus 2 handlebodies.

Let $D$ be 
a properly embedded compressing 2-disk in $M$, such that
$\partial D$ is essential in $S=\partial H$. 
Of course, we may suppose
that $H$ does not determine a Heegaard
splitting of $S^3$ (otherwise $H$ is unknotted, and we are done).
Then we can assume that $D$
is disjoint from a suitable compressing disk $E$ in $H$ having
boundary that does not non separate $S$ (see~\cite[Lemma 3]{scharlemann4}). Compress $H$ along $E$ to
obtain a solid torus $H_1$, with boundary $S_1$. If $\partial D$ is
not essential in $S_1$, then $\partial D$ bounds a disk $D'$ also in
$H$, so that the union of $D$ and $D'$ in a 2-sphere that splits a
suitable split spine of $H$.  If $\partial D$ is essential in $S_1$,
then $H_1$ is unknotted and also in this case it is easy to conclude
that $H$ admits a split-spine.
\end{proof}

\smallskip

\begin{remark}\label{tangled}{\rm
Let us say that a (f8)-spine $\Gamma$ is {\it tangled} if there does
not exist any Whitehead move $\tilde \Gamma \to \Gamma$ such that
$\tilde \Gamma $ is a split (hc)--spine. As far as we understand it,
the statement of Proposition 4.1 of \cite{Ghuman2} is obtained from
the statement  of Proposition \ref{1-K-char} above, just by
replacing the first point with: {\it ``There exists a tangled
(f8)-spine of $H$.''}

\smallskip

Clearly this is a strictly weaker hypothesis: being $(2)_S$-knotted is
equivalent to require that {\it every} (f8)-spine is tangled. In fact, it is easy to show
that an unknotted $H$ actually admits both tangled and
untangled (f8)-spines (see \emph{e.g.}~Figure~\ref{handcuff2}), and that 
in this case $\compl (H)$ has compressible
boundary indeed. Notice however
that the arguments in~\cite{Ghuman2} should provide a slightly different proof of Proposition~\ref{1-K-char}.} 
\end{remark}

\smallskip

The following Proposition is close in spirit to Remark
\ref{knot-vs-unknot}.

\begin{prop}\label{non 2-K} 
(a) Let $H$ be $(2)_S$-knotted and not $(2)_L$-knotted. 
Then, up to isotopy, there exists a unique (hc)--spine of $H$ with split constituent link.

\smallskip

(b) If $H$ is not $(2)_S$-knotted then it admits a unique split spine,
up to diffeomorphims of $S^3$ that leave $H$ invariant.
\end{prop}
\begin{proof}
(a):
Let $\Gamma_0,\Gamma_1$ 
be (hc)--spines of $H$ with split constituent links,
and for $i=1,2,$ take a separating meridian disk $D_i\subseteq H$ dual 
to the isthmus of $\Gamma_i$. Observe that by adding a $2$--handle 
to $\compl(H)$ along $D_i$ we obtain a reducible manifold, so
Theorem~6.1 in~\cite{taylor} (see also~\cite{scharlemann3}) implies that $D_1$ and $D_2$ may be
isotoped to be disjoint. Since both $D_1$ and $D_2$ separate $H$, this easily implies that  they are parallel in $H$,
and this gives in turn that $\Gamma_1$ is isotopic to $\Gamma_2$ in $H$.

\smallskip

(b): Every split spine of $H$ determines a sphere transversely intersecting
$\partial H$ in a simple closed curve that separates $\partial H$ into two
once--punctured tori.  
In the language of~\cite{tsukui}, such a sphere decomposes
the pair $(S^3,\partial H)$ into its \emph{prime factors}.
Now, the main theorem in~\cite{tsukui} ensures that the pair
$(S^3,\partial H)$ admits a unique decomposition into prime factors,
up to homeomorphism, and this concludes the proof.
\end{proof}

\section {Quandle coloring obstructions}\label{quandle-obs}

\begin{figure}[htbp]
\begin{center}
 \includegraphics[height=2cm]{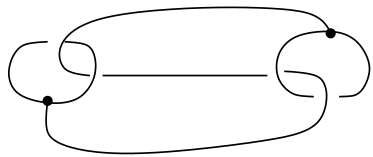}
\caption{\label{handcuff3} A $(3)_S$--knotted $H(\Gamma)$ with trivial
$L_\Gamma$.}
\end{center}
\end{figure}

New invariants of
``links of spatial handlebodies'' have been recently defined
by Ishii 
in~\cite{ishii1, ishii2}. Such invariants are 
based on the analysis of the possible colorings of diagrams, where colors
are intended to belong to
a {\it finite quandle of type $k$} (see below). For example,
by using the simplest instance of these invariants, it has been
remarked in \cite{ishii2} that the handlebody $H(\Gamma)$
corresponding to the spine $\Gamma$ of Figure \ref{handcuff3} is
$(1)_S$-knotted. This example also shows that {\it $(1)_S$-knotting
  does not imply $(1)_L$-knotting}. In fact, using quandle invariants we will show that
  this $H(\Gamma)$ is $(2)_S$--knotted (see Proposition~\ref{1 knotted}),
  so that
  $(2)_S$--knotting does not imply $(1)_L$--knotting. 
  We will
  come back to this example also in Proposition~\ref{alternative}
  and in 
  Proposition~\ref{via-A-obs},
    where we will give two different proofs of the fact that $H(\Gamma)$ is $(3)_S$--knotted.

We are going to show that by means of the same simplest quandle coloring
invariants we can derive more information about our partial order
on the instances of knotting. More precisely in this Section we prove the following:

\begin{prop}\label{K2} 
Let $H$ vary among the (genus 2) spatial handlebodies. Then:
$$\xymatrix{H\ {\rm is}\ (2)_S-{\rm knotted} \ar@{=>}[r]|| & 
 H\ {\rm is}\ (1)_L-{\rm knotted} & ({\rm see\ Proposition\ \ref{1 knotted}}) }
 $$\vspace{-.9cm}\par
$$\xymatrix{H\ {\rm is}\ (1)_L-{\rm knotted} \ar@{=>}[r]|| & 
 H\ {\rm is}\ (2)_S-{\rm knotted} & ({\rm see\ Proposition\ \ref{1Lvs1S:lem}}) }
 $$\vspace{-.9cm}\par
$$\xymatrix{H\ {\rm is}\ (2)_L-{\rm knotted} \ar@{=>}[r]|| & 
 H\ {\rm is}\ (3)_L-{\rm knotted} & ({\rm see\ Proposition\ \ref{2Lvs3L}}) }
 $$\vspace{-.9cm}\par
$$\xymatrix{H\ {\rm is}\ (2)_L-{\rm knotted} \ar@{=>}[r]|| & 
 H\ {\rm is}\ (3)_S-{\rm knotted} & ({\rm see\ Proposition\ \ref{H4}}) }
 $$
\end{prop}

Preliminarily, we need to recall a few facts from  \cite{ishii1,
ishii2} that allow us to compute these invariants in our cases
of interest: either for genus 2 handlebodies or for 2--component links.   
\medskip

\subsection{Quandles: definitions and examples}\label{quandle-recall1}

A {\it quandle} $\X=(X,*)$ is a non-empty set $X$ with a binary
operation that verifies the following axioms.
For every $a,b,c \in X$, we have:
%\smallskip
\begin{enumerate}
\item[(Q1)] $a * a=a$ ;
%\smallskip
\item[(Q2)] $(a * b) * c= (a * c) * (b * c)$ ;
%\smallskip
\item[(Q3)] $S_a(x):= x * a$ defines a bijection on $X$.
\end{enumerate}
\medskip

\noindent{\bf Notation:} for every $a,b \in X$, for every $m\in \N$, set 
$$a*^0b=a,\ a*^1b=a*b,\ a*^2b=(a*b)*b,\ a*^mb=(a*^{m-1}b)*b \ . $$
A quandle $(X,*)$ is {\it of type $k \geq 2$} if $k$ is the minimum
positive integer such that for every $a,b\in X$, we have $a*^kb=a$.
\medskip

\noindent{\bf Dihedral quandles.}
In a sense, the simplest quandles are the
{\it $m$-dihedral quandles} $(R_m , *)$ defined as follows. We 
identify the ring $\Z_m := \Z/m\Z$ with the set $\{0,1,\dots, m-1\}$
of canonical representatives, and it is understood that the operations
act on this concrete set of integer numbers. Then, as a set
$R_m=\Z_m$, while the quandle operation $*$ is defined in terms of the
usual ring operations by: $a*b=2b-a$. It is immediate that {\it
  $(R_m ,*)$ is a finite quandle of type 2}. The name is justified by
the fact that $(R_m ,*)$ can be identified with the set of reflections
of a regular $m$-gon with conjugation as quandle operation.
\medskip

\noindent{\bf Tetrahedral quandle.}
Another important simple example is the {\it tetrahedral quandle}
$(S_4,*)$, where as a set $S_4 = \Z_2[t,t^{-1}]/(t^2+t+1)$, and the
quandle operation is defined in terms of the usual ring operations by:
$a*b = ta+(1-t)b$. It is easy to verify that this is {\it a finite
  quandle of type 3}. 
  %In fact $S_4$ can be identified with the
%4--elements field $\mathcal{F}_4=\{[0],[1], [2], [3]\}$, via the identifications
%$[0]=[0],[1]=[1], [2]=[t], [3]=[t+1]$, and with the quandle operation $*$,
%defined in terms of the field operations by $a*b = [2]a-[2]b $.
\medskip

\noindent{\bf Alexander quandle.}
The above examples belong to the class of so called {\it Alexander
  quandles} $(M,*)$ defined as follows. Let $\Lambda :=
\Z[t,t^{-1}]$. Consider any $\Lambda$-module $M$, as a set, with a
quandle operation $*$ defined (in terms of the usual operations on the
$\Lambda$-module $M$) by $a*b = ta+(1-t)b$.  For every positive
integer $m\geq 2$ and every Laurent polynomial $h(t)\in \Z[t,t^{-1}]$, the module
$M_m= \Z_m[t,t^{-1}]/h(t)$ is an example of Alexander quandle.  
The quandle $(M_m,*)$ is
finite 
if the coefficients of the highest and lowest degree terms are
units in $\Z_m$.  Also 
observe that the dihedral quandle
$(R_m , *)$ introduced above is isomorphic to the Alexander
quandle $(M_m,*)$ associated to the polynomial $h(t)=t+1$. 

\subsection{Quandle coloring invariants.}\label{quandle-recall2}
We are ready to describe the quandle coloring invariants.
Let us fix a finite quandle $\X=(X,*)$ of type $k$.
\smallskip

Let us consider either an (hc)--spine $\Gamma$ of a genus 2 handlebody
$H$, or a 2--component link $L$.  Fix an ordering and an auxiliary
orientation $\omega$ of the two components $K_1, K_2$ of $L_\Gamma$
(resp.~$L$).  

A \emph{$\Z_k$--cycle} on $(\Gamma,\omega)$
(resp.~on $(L,\omega)$) associates to the isthmus of $\Gamma$
the value $0\in\mathbb{Z}_k$, 
and takes an arbitrary value $z_i\in \Z_k$ on
$K_i$ for $i=1,2$.  Therefore, every (hc)--spine and every link supports
exactly $k^2$ different $\Z_k$-cycles, encoded by the couples
$z=(z_1,z_2)$ (the $0$-value associated to the isthmus being understood).
\smallskip
 
Let $\Dd$ be a given diagram of $\Gamma$ (resp.~$L$). Clearly the
orientation $\omega$ and any $\Z_k$-cycle $z=(z_1,z_2)$ descend on
$\Dd$. An {\it arc} of $\Dd$ is an embedded curve in $\Dd$ having as
endpoints either an under-crossing or a vertex of $\Gamma$. 
\smallskip

\begin{figure}[htbp]
\begin{center}
 \includegraphics[height=3cm]{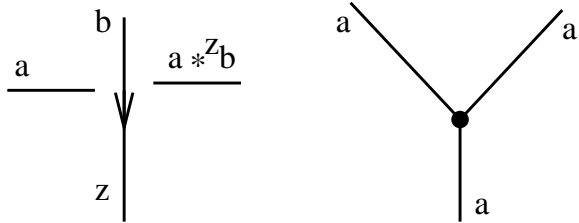}
\caption{\label{quandle} Quandle coloring.}
\end{center}
\end{figure}

By definition, a $\X$-{\it coloring} of $(\Dd,\omega,z)$ assigns to each
arc $e$ of $\Dd$ an element $a(e)\in X$, in such a way that at each
crossing or vertex of $\Dd$ the conditions shown in Figure \ref{quandle}
are satisfied.

Let us denote by $C_\X(\Dd,\omega,z)$ the number of such
$\X$-colorings. Then, by varying the  $\Z_k$- cycle $z$, we obtain a 
{\it non}-ordered $k^2$-uple of positive integers denoted by 
$C_\X(\Dd,\omega)$.
\smallskip

By analyzing the behaviour of $C_\X(\Dd,\omega)$ under the spine moves, it is proved in
\cite{ishii2} that $C_\X(\Dd,\omega)$ does not depend on the choice of
the orientation $\omega$ (so that it makes sense to denote it by
$C_\X(\Dd)$), nor on the choice of the diagram $\Dd$, and not even on
the choice of the (hc)--spine $\Gamma$ of $H$.  It follows that:
\begin{prop}\label{quandle-inv} 
  If $Y$ is either a genus 2 spatial handlebody $H$ or a 2--component
  link $L$, and $\Dd$ is a diagram of either any (hc)--spine $\Gamma$ of
  $H$ or of $L$, then
$$C_\X(Y):= C_\X(\Dd)\in \N^{k^2}/\mathfrak{S}_{k^2} $$
is a well defined isotopy invariant of $Y$ (where $\mathfrak{S}_{k^2}$
is the group of permutations on $k^2$ elements).
\end{prop}
\medskip

\subsection{Specializing to the dihedral quandles}\label{spec-died} 
Let us concentrate our attention on  the dihedral quandle $\X=(R_p,*)$, assuming also for
simplicity that $p$ is an odd prime number.  Then:
\smallskip

\begin{enumerate}
\item
As $\ \X=(R_p,*)$ is of type 2, the orientation $\omega$ is
  unessential in the definition of colorings, so we can forget about it.
%\smallskip
\item
Let $a$ be the number of arcs of a 
diagram $\Dd$ as above. 
Then, for every $\Z_2$-cycle $z$ on $\Dd$, the corresponding set
of colorings is a linear subspace of $\Z_p^{a}$, determined by a
linear equation at each crossing and at each vertex of $\Dd$.  Hence
$C_\X(\Dd,z)$ is a power of $p$, say $p^d$, $d=d(\Dd,z)$, where $d\leq a$.
%\smallskip
\item
For every $\Z_2$-cycle $z$ as above, the corresponding space of
colorings contains the 1--dimensional subspace made by the constant
colorings, so that $d(\Dd,z)\geq 1$.
\end{enumerate}
\smallskip

By using the above remarks we can collect the information
carried by the invariant $C_\X(Y)$ by means of the following invariant
polynomial:

$$\Phi_p(Y)(t):= \Phi_p(\Dd)(t)= \sum_{z} t^{d(D,z)-1} \in \N[t] \ .$$ 

We stress that this notation could be a bit misleading as it could
suggest that the monomials of $\Phi_p(Y)(t)$ are in some way marked by
the $\Z_2$-cycles. This is true for a given diagram, but this
information is lost when we consider the polynomial as an invariant of
$Y$. This corresponds to the fact that $C_\X(Y)$ is a {\it
  non}-ordered 4--uple of positive integers. Clearly $\Phi_p(Y)(t)$
has at most 4 monomials, the sum of its coefficients is equal to $4$,
and its degree is at most $A-1$, where $A$ is the minimal number of
  arcs, when $\Dd$ varies among all the diagrams of $Y$ (when $Y$ is a link) or of all the (hc)--spines of $Y$
  (when $Y$ is a handlebody).
\medskip

\begin{lem}\label{easycase}
For every diagram $\Dd$ of a (hc)--spine
$\Gamma$ of $H$ (resp.~of a 2--component link $L$) we have that 
$d(\Dd,(0,0))=1$ (resp.~$d(\Dd,(0,0))=2$).
\end{lem}
\begin{proof}
Let $\Dd$ be the diagram of a (hc)--spine $\Gamma$ (resp.~of a 2--component link $L$). 
Then, the
assignment of a color to each arc
defines a coloring of $\Dd$ associated to the trivial cocycle
if and only if it is constant (resp.~it is constant on the components of $L$). 
\end{proof}

\subsection{Quandle obstructions}\label{obs}
Let $p$ be an odd prime. From now on, we will consider only 
the dihedral quandle $\X=(R_p ,*)$. 
By using the simplest quandle invariants associated to $\X$ we are
going to determine necessary conditions (``obstructions'') for a given $H$ 
to be $(2)_S$-unknotted, $(1)_L$-unknotted or $(2)_L$-unknotted.
The following Lemma is not strictly necessary to our purposes
(in fact, the statement of Corollary~\ref{L1} regarding $(1)_L$-unknotted
handlebodies
may also be deduced by Lemma~\ref{2-K-obs} below --
see Remark~\ref{L2rem}). However, it
establishes an interesting relation between the quandle invariants of a handlebody
and those of its constituent links. Throughout the whole Subsection, let $H$ be a spatial handlebody.

\begin{lem}\label{HvsL}
Let $\Gamma\in\mathcal{S}(H)$ and 
take the corresponding $L_\Gamma\in\Ll (H)$. If 
$$\Phi_p(L_\Gamma)(t)= t
+ t^{m_1} +t^{m_2}+t^{m_3},\quad m_1\leq m_2\leq m_3\ , $$
%(where the $m_j$'s are not necessarily
%distinct), 
then there exist $n_1,n_2,n_3\in\mathbb{N}$ such that
%$\Phi_p(H)(t)$ is of the form 
$$\Phi_p(H)(t)= 1+t^{n_1}
+t^{n_2}+t^{n_3},\quad n_1\leq n_2\leq n_3\ ,$$
and for $j=1,2,3$ the integer
 $n_j$ is such that either $n_j=m_j$ or
$n_j=m_j-1$.
\end{lem}

\begin{proof}
We can take a
diagram $\Dd$ of $\Gamma$ such that an open neighbourhood of the
isthmus in $\Gamma$ bijectively projects onto its image in $\Dd$.
We can also ask that
this image does not intersect the remaining part of the
diagram. It follows that by removing the interior of the isthmus from
$\Dd$ we get a diagram $\Dd'$ of $L_\Gamma$. Then it is clear that,
for every cycle $z$, the linear system computing $C_\X(\Dd,z)$ is
obtained from the linear system computing $C_\X(\Dd',z)$ just by
adding one equation due to the coloring conditions at vertices. 
Since both these linear systems admit the constant solutions,
for every fixed cycle
the set of solutions corresponding to $\Dd$ 
is an affine subspace 
of codimension $0$ or $1$
of the space
of solutions corresponding to $\Dd'$. 

Together with Lemma~\ref{easycase}, this shows that if
$\Phi_p(L_\Gamma)(t)= t
+ t^{m_1} +t^{m_2}+t^{m_3}$, $m_1\leq m_2\leq m_3$, then
there exist $n_1,n_2,n_3\in\mathbb{N}$, $n_1\leq n_2\leq n_3$,
such that
$\Phi_p(H)(t)= 1+t^{n_1}
+t^{n_2}+t^{n_3}$ and $m_j-1\leq n_{\tau(j)}\leq m_j$, where
$\tau\in\mathfrak{S}_3$ is a permutation. In order to conclude, it is now
sufficient to show that we may assume that $\tau$ is the identity. But 
this is a consequence of the following easy
\smallskip

\noindent{\bf Claim:} Let $a_1\leq \ldots\leq a_n$ and $b_1\leq \ldots\leq b_n$ be 
non--decreasing sequences
of real numbers such that $b_i-1\leq a_{\tau(i)}\leq b_i$ for every $i=1,\ldots,n$,
where $\tau\in\mathfrak{S}_n$ is a permutation. Then, we have
$b_i-1\leq a_i\leq b_i$ for every $i=1,\ldots,n$.

In fact, for every given $i_0$, the assumption implies that the set
$\{i\, |\, a_i\leq b_{i_0}\}$ contains at least $i_0$ elements, so
$a_{i_0}\leq b_{i_0}$. In the very same way one can show that
$a_{i_0}\geq b_{i_0}-1$, whence the conclusion.
\end{proof}

\begin{cor}\label{L1}
Let $L\in \Ll (H)$ be a constituent link of $H$. Then,
$$
\deg \Phi_p (L)-1 \leq \deg \Phi_p (H)\leq \deg \Phi_p (L).
$$
In particular, if $H$ is $(1)_L$-unknotted, then
$\deg \Phi_p (H)\leq 1$.
\end{cor}
\begin{proof}
 The first statement is an immediate consequence of the previous Lemma.
If $H$ is $(1)_L$-unknotted, then the trivial link $L$
belongs to $\Ll (H)$. Since $\Phi_p (L)=4t$, the conclusion follows.
\end{proof}

We now come to the obstruction to being $(2)_S$-unknotted.

\begin{lem}\label{1-k-obs} 
If $H$ is not $(2)_S$-knotted,
then there exist $h_1,h_2\in\mathbb{N}$
such that 
$$ \Phi_p(H)(t)= 1+t^{h_1}+t^{h_2} + t^{h_1+h_2} \ . $$
\end{lem}
\begin{proof} Since $H$ admits a split (hc)--spine $\Gamma$, it has a diagram
$\Dd$ of the form shown in Figure \ref{non-1-K}, where it is
understood that removing the interior of the isthmus one gets diagrams
$\Dd_j$ of the constituent knots $K_j$, $j=1,2$, so that
every box includes a 1--string sub-diagram of the corresponding $\Dd_j$.
The symbol $a$ belongs to $R_p=\{0,1,2,\dots,\ p-1 \}$ and refers to
a portion of a $\X$-coloring of $\Dd$. For every $\Z_2$-cycle
$z=(z_1,z_2)$ on $\Dd$, every $\X$-coloring of $(\Dd,z)$ restricts to a
$\X$-coloring of both the diagrams $(\Dd_j,z_j)$.

\begin{figure}[htbp]
\begin{center}
 \includegraphics[height=3cm]{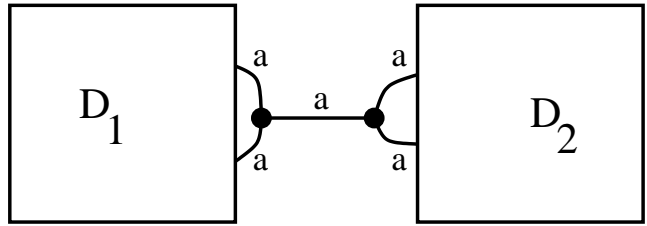}
\caption{\label{non-1-K} A $(2)_S$-unknotted diagram.}
\end{center}
\end{figure} 

For every $a\in R_p$ and $j=1,2$, denote by $n_{j,a}\in \N$ the
number of $\X$-colorings of $(\Dd_j,1)$ that extend the color $a$ near
the vertex of $\Dd$ contained in $\Dd_j$. 
Then we have
$$C_\X(\Dd,(1,0))= \sum_a n_{1,a}$$
$$C_\X(\Dd,(0,1))= \sum_a n_{2,a}$$
$$C_\X(\Dd,(1,1))= \sum_a n_{1,a}n_{2,a} \ .$$

Due to the definition of the dihedral quandle, since $p\neq 2$, for every given $a,c \in R_p$ there exists a unique
$b$ such that $a*b=c$. Moreover,
thanks to the axioms in the definition of quandles, if $b\in\X$ is a fixed color and we replace by $d*b$ 
every color $d$ occurring in a $\X$-coloring of
$(\Dd_j,1)$ that takes the value $a$ near the vertex, then
we get an $\X$-coloring that takes the value $c$ near the vertex. 
As a consequence, for $j=1,2$ there exists $h_j\in\mathbb{N}$ such that 
$n_{j,a}=p^{h_j}$ for every $a\in\X$. 
\smallskip

It follows that we have
$$
\begin{array}{llllll}
 C_\X (\Dd,(0,0))&=&p\ ,\quad &
 C_\X(\Dd,(1,0))&=& p^{h_1+1}\ ,\\
 C_\X(\Dd,(0,1))&=& p^{h_2+1}\ ,\quad
 & C_\X(\Dd,(1,1))&=&p^{h_1+h_2+1}\ ,
 \end{array}
$$
whence the conclusion.
\end{proof}

\begin{lem}\label{2-K-obs} Suppose that $H$ is $(2)_L$-unknotted.
Then there exist $h_1,h_2,h_3\in\mathbb{N}$ such that 
$$ \Phi_p(H)(t)= 1+t^{h_1}+t^{h_2} + t^{h_3}, \qquad h_1\leq h_2\leq h_3, $$
 where $h_3=h_1+h_2$ or $h_3=h_1+h_2+1$.  
Moreover, if $L\in\Ll(H)$ is a split link, then 
$$
\Phi_p (L)(t)=t+t^{h_1+1}+t^{h_2+1}+t^{h_1+h_2+1}\ .
$$
\end{lem}
\begin{proof} Let us take a spine $\Gamma$ of $H$ such that
$L_\Gamma$ is a split
 link. Then there is a diagram $\Dd$ of $\Gamma$ of the form shown in
 Figure \ref{non-2-K}. Here $h=2n+1$ is an odd positive integer and
 the rectangle in the middle represents $h$ parallel strings. All the
 $h+2$ horizontal strings belong to the isthmus. If we remove the
 interior of the isthmus we obtain a split diagram $\Dd_1\cup \Dd_2$ of 
$L_\Gamma=K_1\cup K_2$.

\begin{figure}[htbp]
\begin{center}
 \input{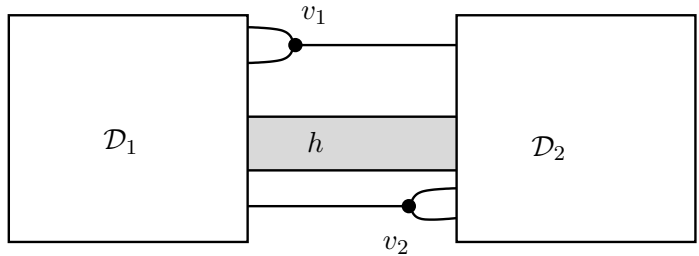}
\caption{\label{non-2-K} A $(2)_L$-unknotted diagram.}
\end{center}
\end{figure}

For $i=1,2$, let $C_i$ be the number of colorings
of $\Dd_i$ corresponding to the non--trivial cycle on $\Dd_i$, and let
$h_1,h_2\in\mathbb{N}$ be such that
$C_i=p^{h_i+1}$. Since any $\mathbb{Z}_2$-cycle on $\Gamma$ vanishes
on the isthmus, we easily get 
$$
\begin{array}{llllll}
 C_\X (\Dd_1\cup\Dd_2,(0,0))&=&p^2\ ,\qquad & C_\X (\Dd,(0,0))&=&p\ ,\\
 C_\X (\Dd_1\cup \Dd_2,(1,0))&=& p^{h_1+2}\ ,\qquad & C_\X(\Dd,(1,0))&=& p^{h_1+1}\ ,\\
C_\X (\Dd_1\cup \Dd_2,(0,1))&=& p^{h_2+2}\ ,\qquad & C_\X(\Dd,(0,1))&=& p^{h_2+1}\ ,\\
 C_\X (\Dd_1\cup \Dd_2,(1,1))&=& p^{h_1+h_2+2}\ .\qquad & & &
 \end{array}
$$
This implies in particular that
$$
\Phi_p (L)(t)=t+t^{h_1+1}+t^{h_2+1}+t^{h_1+h_2+1}\ .
$$

Observe now that 
every $(1,1)$-coloring of $\Dd$ restricts to a $(1,1)$-coloring
of $\Dd_1\cup\Dd_2$. Moreover, once a $(1,1)$-coloring of $\Dd_1\cup\Dd_2$
is fixed, one can try to extend it to a $(1,1)$-coloring of $\Dd$
as follows: the coloring of $\Dd_1$ uniquely determines
the color of the arc of the isthmus starting at the vertex $v_1$;
then, 
following the isthmus from $v_1$ to $v_2$, one assigns to the arcs of the isthmus
the colors uniquely determined by the rules describing the behaviour of colorings
at crossings; finally, one check whether
the color obtained at the arc ending at $v_2$ matches the fixed coloring of $\Dd_2$.
One can express this last condition as a linear equation on the colors of
$\Dd_1\cup \Dd_2$, and this implies in turn that the set of $(1,1)$-colorings
of $\Dd$ admits a bijection with 
a subspace of the $(1,1)$-colorings of $\Dd_1\cup\Dd_2$ having codimension $0$ or $1$. Then
$C_\X(\Dd,(1,1))$ is equal either to $p^{h_1+h_2+2}$ or to $p^{h_1+h_2+1}$,
whence the conclusion.
\end{proof}

\begin{remark}\label{L2rem}
{\rm If $L$ is the trivial link, then $\Phi_p (t)=4t$. Therefore 
 Lemma~\ref{2-K-obs} allows us to refine Corollary~\ref{L1}:
if $H$ is $(1)_L$-unknotted, then we have either
$\Phi_p (H)=4$ or $\Phi_p (H)=3+t$.}
\end{remark}

Our next goal is to use the so obtained obstructions in order to
produce, for example, families of $(2)_S$-knotted
(resp.~$(2)_L$-knotted) handlebodies that are $(2)_L$-unknotted
(resp.~$(3)_L$-unknotted, and even $(3)_S$-unknotted). Having this in mind, it is useful to
introduce and study some elementary tangles that we will combine in
order to get the desired examples.

\subsection{The tangle $E(q)$}\label{tangleE}  
\begin{figure}[htbp]
\begin{center}
\input{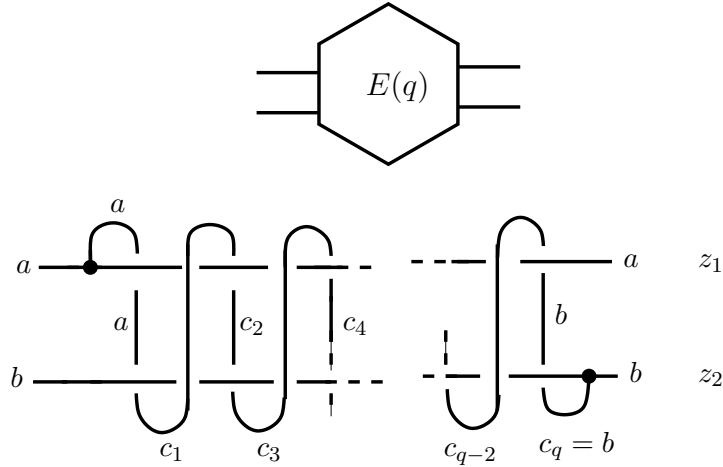}
\caption{\label{spire1} The tangles $E(q)$.}
\end{center}
\end{figure}

Let $q$ be an odd prime, and
consider the tangle $E(q)$ of Figure \ref{spire1}.  Here $q$ is an odd positive integer, the $z_i \in \Z_2$, $i=1,2$, label the horizontal lines and play
the r\^ole of a $\Z_2$-cycle $z$, while the colors $a,b, c_j$ belong
to $R_p$ and refer to a generic $\X$-coloring of this tangle, relative
to the given $z=(z_1,z_2)$. Recall that our assumptions imply that every cycle
vanishes on the isthmus, so for every $z$ any admissible coloring is constant 
on the horizontal lines. As a consequence,  
$a$ (resp.~$b$) are constant along the top (resp.~bottom) line
of the diagram.

The following Lemma computes the number of colorings of the tangle $E(q)$.

\begin{lem}\label{E} 
  For every $p$ denote by $C_p E(q,z,a,b)$ the number of
  $(R_p,*)$-colorings of $E(q)$ relative to $z$ which assume the values
$a$ and $b$ respectively on the top and the bottom line of the diagram.
  Then:
\begin{enumerate}
 \item 
If $q=p$ and $z=(1,1)$, then 
 $C_p E(q,z,a,b)=1$ for every $(a,b)\in R_p^2$.
\item
In all the other cases (\emph{i.e.}~if $p\neq q$ or if $z\neq (1,1)$), then
$C_p E(q,z,a,b)=0$ if $a\neq b$, and $C_p E(q,z,a,b)=1$ if $a=b$
(and in this case, we have only the constant coloring assigning
the color $a=b$ to every arc of the diagram). 
\end{enumerate}
\end{lem}
\begin{proof}
Let us consider only the case $z=(1,1)$, the other cases being easier.
With notation as in Figure~\ref{spire1}, we have $c_1=a*b=2b-a$, and
$c_{2l+1}=(c_{2l-1}*a)*b=c_{2l-1}+2(b-a)$ for every $l=1,\ldots, (q-1)/2$.
We get therefore $c_q=(q-1)(b-a)+2b-a$, and 
the assigned coloring of the horizontal rows extends 
(in a unique way) to the whole diagram if and only if $(q-1)(b-a)+2b-a=b$,
\emph{i.e.}~if and only if $q(b-a)=0$. 
 This equality holds in $\Z_p$ if and only if $p$ divides $q$ or $a=b$,
whence the conclusion.
\end{proof}

\subsection{$(2)_S$-knotting does not imply $(1)_L$-knotting}\label{1Svs1L}
\begin{figure}[htbp]
\begin{center}
\input{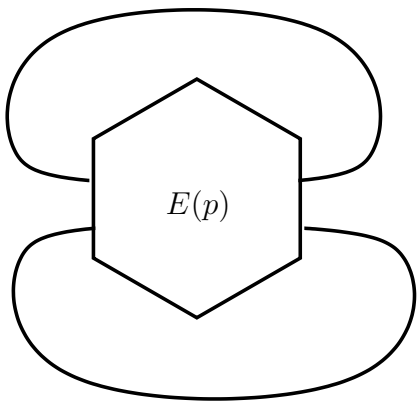}
\caption{\label{1-knotted} The spine $\Gamma_1(p)$. The hexagonal box represents the tangle described in Figure~\ref{spire1}.}
\end{center}
\end{figure}

We are ready to construct the first
pertinent family of examples.  For every odd prime $p$, consider the
(hc)--spines $\Gamma_1(p)$ of Figure \ref{1-knotted}, and set $H_1(p)=H(\Gamma_1(p))$.%qui

\begin{prop}\label{1 knotted}
  For every prime $p$, $H_1(p)$ is $(2)_S$-knotted and 
  $(1)_L$-unknotted.
  Moreover, if $p$ and $p'$ are different
prime numbers, then $H_1(p)$ and $H_1(p')$ are not isotopic.
\end{prop}
\begin{proof} 
Let $p$ be a fixed prime number.
Since the constituent link of $\Gamma_1 (p)$ is trivial,
by definition $H_1(p)$ is not $(1)_L$-knotted. 

An easy application of Lemma~\ref{E} implies that the numbers
of distinct $(R_p,\ast)$--colorings of $\Gamma_1(p)$ with respect to the cycle $z=(z_1,z_2)$
is equal to $p^2$ if $z=(1,1)$ and to $p$ otherwise. This implies that
$$\Phi_p(H_1(p))(t)= 3+t \ . $$
Together with Lemma~\ref{1-k-obs}, this implies that $H_1 (p)$ is 
$(2)_S$--knotted.

Take now a prime number $p'\neq p$. Lemma~\ref{E} now easily implies that
$\Phi_p(H_1(p'))(t)= 4\neq \Phi_p(H_1(p))(t)$,
so $H_1(p')$ is not isotopic to $H_1(p)$.
\end{proof}

\begin{remark}
 {\rm Building respectively on the theory of handlebody patterns developed in Section~\ref{more-bpm}  and on the use
 of Alexander-type invariants, in Propositions~\ref{alternative} and~\ref{via-A-obs} we give two different proofs
 of the stronger fact that $H_1(p)$
is $(3)_S$-knotted for every prime $p$. 
}
\end{remark}

\subsection{The tangle $O(q)$}
\begin{figure}[htbp]
\begin{center}
\input{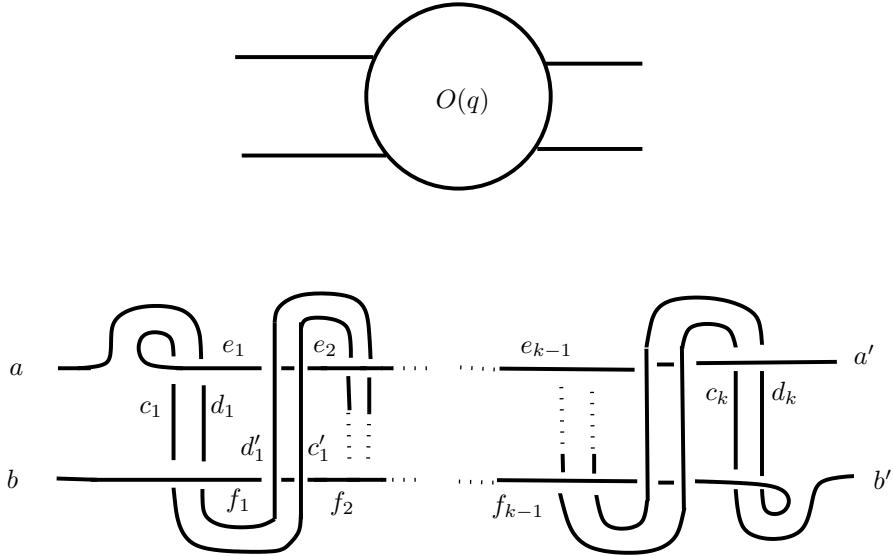}
\caption{\label{spire2} The tangles $O(q)$, where $q=2k-1$.}
\end{center}
\end{figure}

We now 
consider the tangle $O(q)$ of Figure \ref{spire2}, which can be obtained
from $E(q)$ as follows: first, we attach to the band bounded by the two horizontal lines
of $E(q)$ a 2--dimensional 1--handle whose core coincides
with the isthmus; then, we define $O(q)$ to be the boundary of the
so obtained surface. Observe that $O(q)$ is the union of two arcs, one of which
is entering and exiting the diagram on the left, the other on the right.

In what follows, we will be interested in the colorings of $O(q)$ corresponding
to $\Z_2$--cycles which either vanish on both the components of $O(q)$, or
take the value $1$ on both the components of $O(q)$. We denote the corresponding
labelled tangles respectively by $O(q,0)$ and $O(q,1)$.
We also denote by $a,b\in R_p$ a pair of ``input'' colors which are assigned to the arcs
entering and exiting the diagram of $O(q)$ on the left (see Figure~\ref{spire2}).
The following Lemma describes the possible pairs of ``output'' colors
$a',b'\in R_p$ which are associated to the entering/exiting arcs on the right
by a global coloring that extends the ``input'' datum $(a,b)$.

\begin{lem}\label{O} 
For $z\in\{0,1\}$ and $q$ an odd prime,
let $C_p O(q,z,a,b,a',b')$ be the number of
$(R_p,*)$-colorings of $O(q,z)$ that 
extend the given input/output data
$(a,b,a',b')$. 
\begin{enumerate}
 \item Suppose that $p=q$. Then, $C_p O(q,1,a,b,a',b')=p$
if $a=a'$ and $b=b'$, and $C_p O(q,1,a,b,a',b')=0$
if $a\neq a'$ or $b\neq b'$. 
\item
Suppose that $p\neq q$, and let us fix $(a,b)\in R_p^2$. Then
there exists a unique pair $(a',b')\in R_p^2$ such that
$C_p O(q,1,a,b,a',b')=1$. Moreover, if $a=a'$ or $b=b'$ then 
we have $a=a'=b=b'$. 
If $(a'',b'')\neq (a',b')$ is any other pair of colors, then
$C_p O(q,1,a,b,a'',b'')=0$.
\item
Suppose that $z=0$. Then,
%If $z=0$, 
we have $C_p O(q,0,a,b,a',b')=1$
if $a=b$ and $a'=b'$, and $C_p O(q,0,a,b,a',b')=0$
if $a\neq b$ or $a'\neq b'$. 
\end{enumerate}
\end{lem}
\begin{proof}
Point~(3) is obvious, so we concentrate our attention on the colorings
of $O(q,1)$. 

Let us set $k=(q+1)/2$, and label the vertical arcs of the diagram 
by the colors $$c_1,d_1,\ldots,{c_k,d_k}, 
c'_1,d'_1,\ldots,c'_{k-1},d'_{k-1}$$ as in
Figure~\ref{spire2}. Also label by the color $e_i$ (resp.~$f_i$)
the top (resp.~bottom) arc passing over the arcs labelled by $c_i$ and $d_i$.

We set $x=c_1$, and look for the values that the labels introduced above must take in order
to define a coloring of the diagram that extends the given input/output data.
By looking at crossings from the left to the right we obtain
$$
e_1=x,\quad f_1=b,\quad d_1=a\ast x=2x-a\ ,
$$
%Moreover, we have
$$
%\begin{array}{lllllll}
c'_{i}=c_i\ast f_i,\qquad  d'_{i}=d_i\ast f_i, \qquad i=1,\ldots,k-1\ ,
$$
and%qui
$$
e_{i+1}=(e_i\ast d'_i)\ast c'_i,\quad
f_{i+1}=(f_i\ast d'_i)\ast c'_i,\quad
c_{i+1}=c'_i\ast e_{i+1},\quad
d_{i+1}=d'_i\ast e_{i+1}
$$
for every $i=1,\ldots,k-1$.
Therefore, an easy induction shows that 
\begin{equation}\label{colori}
\begin{array}{llllll}
e_k&=&(2k-1)x-(2k-2)a,\quad &
f_k&=&b+(2k-2)(x-a),\\
c_k&=&(6k-5)x-(4k-4)a-(2k-2)b,\quad &
d_k&=&(6k-4)x-(4k-3)a-(2k-2)b.
\end{array}
\end{equation}
Now, the fixed labels $a,b,a',b',x$ extend (uniquely) to a coloring of the whole diagram
if and only if 
$$
\left\{
\begin{array}{l}
e_k=a'\\
f_k=d_k\\
c_k\ast f_k=b'
\end{array}\right. \ .
$$
By~\eqref{colori}, since $q=2k-1$
this linear system is equivalent to the system
\begin{equation}\label{sistema}
\left\{
\begin{array}{l} 
qx-(q-1)a=a'\\
q(2x-a-b)=0\\
-qx+(q+1)b=b'
\end{array}\right. \ .
\end{equation}

Let us assume that $p=q$. Then, the system~\eqref{sistema}
reduces to the conditions $a=a'$, $b=b'$. If these conditions are satisfied,
every choice for $x$ can be (uniquely) extended to a coloring of the diagram, while if
$a\neq a'$ or $b\neq b'$ there do not exist colorings extending
the input/output data $a,a',b,b'$. This proves point~(1).

If $p\neq q$, 
then the second equation of~\eqref{sistema}
implies that $x=(a+b)/2$ (recall that $p$ is odd, so that $2$ is invertible in $\Z_p$). 
Then, looking at the other equations we see that the system admits a (unique)
solution if and only if $2a'=(2-q)a+qb$ and $2b'=-qa+(q+2)b$.
In particular, if $a=a'$ or $b=b'$ we necessarily have $a=a'=b=b'$,
whence the conclusion.
\end{proof}

\subsection{$(1)_L$-knotting does not imply $(2)_S$-knotting}\label{1Lvs1S}
\begin{figure}[htbp]
\begin{center}
\input{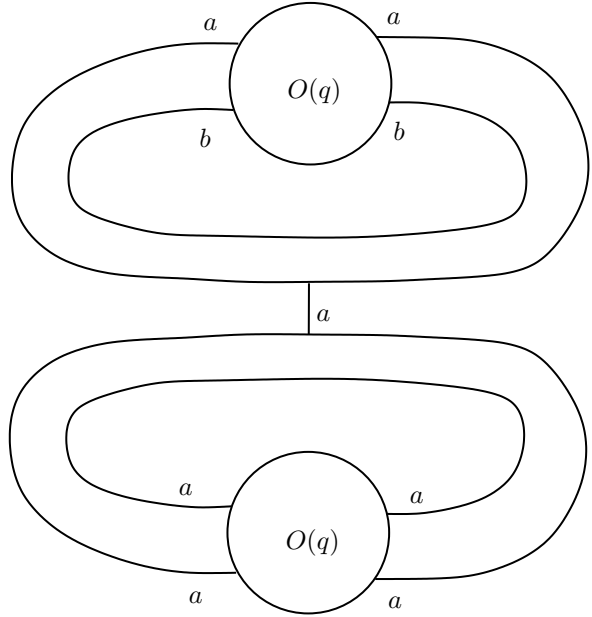}
\caption{\label{gamma1:fig} An $R_p$-coloring of the spine $\Gamma_2(q)$ relative to the cycle $(1,0)$, when $p=q$. The circular boxes
represent tangles as described in Figure~\ref{spire2}.}
\end{center}
\end{figure}

For every odd prime $q$, let us consider the graph $\Gamma_2 (q)$ described in Figure~\ref{gamma1:fig}, and let $H_2(q)=H(\Gamma_2 (q))$. 

We have the following:

\begin{prop}\label{1Lvs1S:lem}
For every prime $p$, $H_2(p)$ is $(1)_L$-knotted and 
  $(2)_S$-unknotted.
  Moreover, if $p$ and $p'$ are different
prime numbers, then $H_2(p)$ and $H_2(p')$ are not isotopic.
\end{prop}
\begin{proof}
Let us fix an odd prime number $p'$. It is clear that 
$\Gamma_2 (p')$ is a split spine of $H_2 (p')$, which is therefore
$(2)_S$-unknotted. By Lemma~\ref{1-k-obs} (and its proof), in order to compute
$\Phi_p (H_2(p'))$ it is sufficient to compute the number
of $z$-colorings of $H_2(p')$ for $z=(1,0)$ and $z=(0,1)$. 

So, let
us suppose that the cycle $z=(1,0)$ assigns the value $1$ (resp.~$0$) to the component
of $L_{\Gamma_2(p')}$ on the top (resp.~on the bottom) of Figure~\ref{gamma1:fig}.
An easy application of Lemma~\ref{O} shows that if $p\neq p'$, then the only
$(1,0)$-colorings of $\Gamma_2(p')$ are the constant ones. The same is true for
$(0,1)$-colorings, so by Lemma~\ref{1-k-obs}  we have
$$
\Phi_p (H_2(p'))=4\qquad {\rm if}\ p\neq p'\ .
$$

Suppose now $p=p'$. By Lemma~\ref{O} it is easily seen that the
knot on the top of Figure~\ref{gamma1:fig} admits exactly $p^3$ colorings
relative to the non-trivial cycle. Moreover,
each such coloring uniquely extends to a $(1,0)$-coloring of the whole diagram
of $\Gamma_2(p)$. By the symmetry of the diagram, the same result holds for $(0,1)$-colorings.
Then, by Lemma~\ref{1-k-obs}  we have
$$
\Phi_p (H_2(p))=1+2t^2+t^4.
$$
By Corollary~\ref{L1}, this implies that $H_2(p)$ is $(1)_L$-knotted.
Moreover, if $p'\neq p$ we have
$\Phi_p (H_2 (p))\neq \Phi_p (H_2(p'))$,
so the spatial handlebodies $H_2(p)$ and $H_2(p')$ are not isotopic.
\end{proof}

\subsection{Quandle colorings of bands}\label{bands:sub}
Our next constructions make an extensive use of links and tangles
obtained by ``doubling'' some given knot or tangle, \emph{i.e.}~by
replacing a knot or a tangle with the boundary of a band representing a fixed
framing on the knot or the tangle. Therefore, it is convenient to 
point out some nice features of quandle colorings of bands.

\begin{figure}[h]
\begin{center}
\input{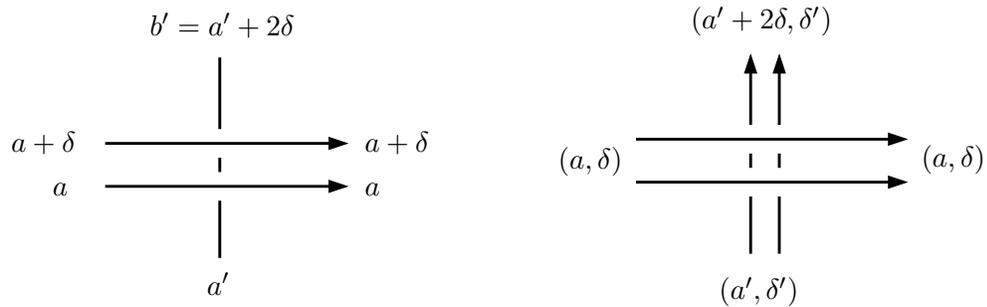}
\caption{\label{bande:fig} 
Quandle colorings of bands.}
\end{center}
\end{figure}

So, let us consider a pair $A$ of parallel arcs in a diagram, labelled with colors $a,b\in R_p$.
After fixing a (coherent) orientation on the arcs, we suppose that 
$a$ (resp.~$b$) is the color of the arc running on the right (resp.~on the left), and
we set $\delta=b-a$. For reasons that will become clear soon, we label 
$A$ by the pair $(a,\delta)\in R_p^2$ (see Figure~\ref{bande:fig}).

Let us now suppose that both the arcs of $A$ are labelled by the non--trivial
$\mathbb{Z}_2$-cycle.
Now, if an arc undercrosses the band $A$ from the right (with color $a'$)
to the left (with color $b'$), it is immediate to observe that in order
to have an admissible coloring the equality $b'=a'+2\delta$ must hold. 
As a consequence, if a band $A'$ undercrosses
$A$ from the right to the left, and  
the portion of $A'$ on the right is labelled by $(a',\delta')$, then the portion
on the left has to be labelled by $(a'+2\delta,\delta')$. In particular,
the parameter $\delta'$ propagates without being affected by crossings.

\subsection{The tangle $\overline{O}(q)$}
\begin{figure}[h]
\begin{center}
\input{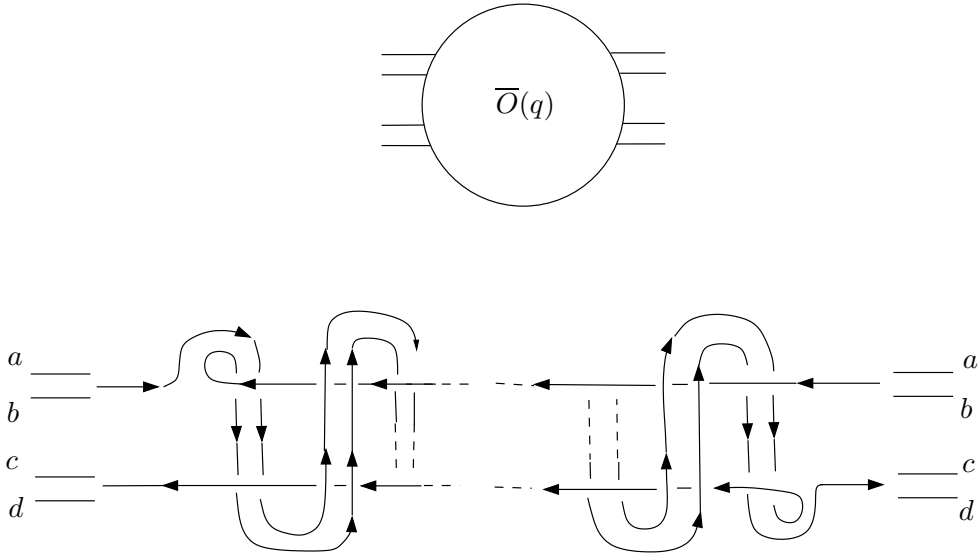}
\caption{\label{spiredoppie:fig} 
The tangle $\overline{O}(q)$.}
\end{center}
\end{figure}

Let us now consider the tangle $\overline{O}(q)$ obtained by replacing
each arc of $O(q)$ with a band, thus obtaining a tangle with four strings. In Figure~\ref{spiredoppie:fig} we represent
$\overline{O}(q)$ by drawing one arrow for each band. We agree that
the $(1,1)$-cycle (resp.~the $(0,0)$-cycle) assigns the value $1$ (resp.~$0$)
to every arc, the
$(1,0)$-cycle (resp.~$(0,1)$-cycle) assigns the value $1$ (resp.~$0$)
on the arcs giving the ``right'' 
boundary of the band (with respect to the orientation given by the arrows)
and the value $0$ (resp.~$1$) on the arcs on the left.

\begin{lem}\label{dO}
  For every $p$ and $z\in\mathbb{Z}_2^2$ let us denote by $C_p \overline{O}(q,z,a,b,c,d)$ the number of
  $(R_p,*)$-colorings of $\overline{O}(p)$ relative to the cycle $z$ 
(in the sense specified above) that extend the coloring described in Figure~\ref{spiredoppie:fig}.
  Then:
\begin{enumerate}
 \item 
If $z=(1,1)$, then $C_p \overline{O}(q,z,a,b,c,d)=1$ if $a=d$ and $b=c$,
and $ C_p \overline{O}(q,z,a,b,c,d)=0$ otherwise.
\item
If $q=p$ and $z=(1,0)$ or $z=(0,1)$, then for every $a\in R_p$ we have $C_p \overline{O}(q,z,a,a,c,d)=p$ 
if $c=d$, and $C_p \overline{O}(q,z,a,a,c,d)=0$ if $c\neq d$.
\item
If $q\neq p$ and $z=(1,0)$ or $z=(0,1)$, then 
for every cycle $z$ we have
$C_p \overline{O}(q,z,a,a,c,d)=1$ if $a=c=d$,
and $C_p \overline{O}(q,z,a,a,c,d)=0$ otherwise.
\end{enumerate}
\end{lem}
\begin{proof}
 (1): The discussion carried out in Subsection~\ref{bands:sub} shows that
a necessary condition for extending the given coloring is that $b-a=c-d$. Let us set
$b-a=c-d=\delta$.  
From the top to the bottom, we may label the bands on the left
by $(b,-\delta)$ and $(c,-\delta)$, and the bands on the right by
$(a,\delta)$ and $(d,\delta)$. Since every band undercrosses itself and the other band
the same number of times (with the same orientation), it follows
that  we need to have $a=d$ and $b=c$. Moreover, it is clear that
if this condition is satisfied, then the coloring extends in a unique way.

\smallskip

(2), (3):  Independently of the cycle $z$, 
it is readily seen that if a coloring assigns
the same color to two parallel arcs bounding a portion of a band, then 
every pair of parallel arcs belonging to that band must have the same color. 
Also observe that this is our case of interest, since we are assigning the color
$a$ on the top arcs both on the right and on the left of the diagram. 
If $z=(1,0)$ or $z=(0,1)$,
this implies in turn that the $z$-colorings of $\overline{O}(q)$ bijectively correspond
to the $(1,1)$-colorings of $O(q)$, \emph{i.e.}~the $z$-colorings
of $\overline{O}(q)$ are exactly the colorings obtained by ``doubling''
a $(1,1)$-coloring of $O(q)$. The conclusion follows now
from Lemma~\ref{O}.
\end{proof}

\subsection{$(2)_L$-knotting does not imply $(3)_L$-knotting}\label{2Lvs3L:sub}
\begin{figure}[h]
\begin{center}
\input{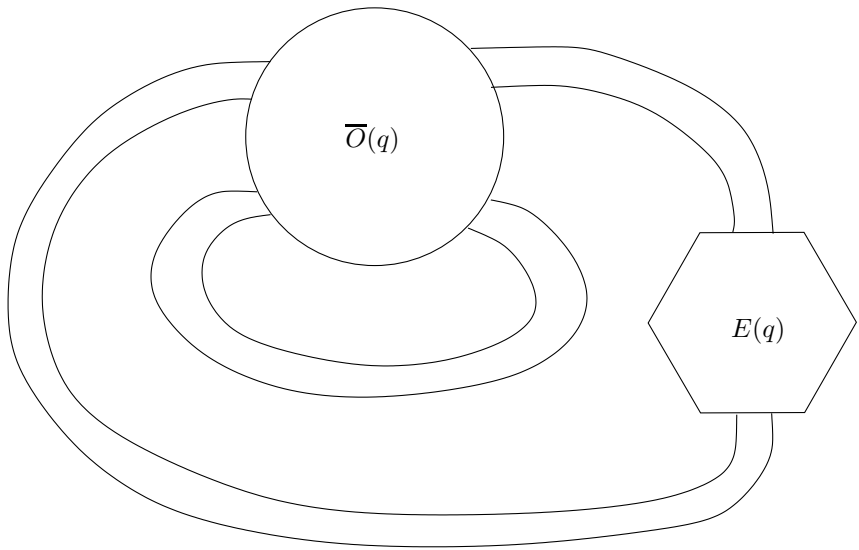}
\caption{\label{gamma3:fig} 
The spine $\Gamma_3(q)$.}
\end{center}
\end{figure}

For every odd prime $q$, 
let us consider the spine $\Gamma_3 (q)$ shown in Figure~\ref{gamma3:fig},
and set $H_3 (q)=H(\Gamma_3(q))$. 
Let $K_1,K_2$ be the constituent knots of $\Gamma_3 (q)$.
For $i=1,2$, 
it is readily seen that the blackboard framing of $K_i$
coincide with the trivial framing,
and this implies that $L_{\Gamma_3 (q)}=K_1 \cup K_2$ is a boundary link (since the linking number of the knots of $K_1$ and $K_2$ is zero, a parallel copy
of a Seifert surface $S_1$ for 
$K_1$ provides a Seifert surface for $K_2$ which is disjoint form $K_1$).

Putting together Lemmas~\ref{E}, \ref{dO}
and~\ref{2-K-obs}
one easily gets the following:

\begin{prop}\label{2Lvs3L}
If $p,q$ are distinct odd primes, then 
$$
\Phi_p (H_3(p))= 1+t+2t^2\ ,\qquad
\Phi_p (H_3(q))=4\ .
$$ 
In particular, $H_3(p)$ is $(2)_L$-knotted and
$(3)_L$-unknotted, and $H_3(p)$ is not isotopic to $H_3(q)$.
\end{prop}

In Subsection~\ref{2L2S} we prove the stronger result
that $(2)_L$-knotting does not imply $(3)_S$-knotting. One may wonder
if this result could be achieved just by replacing the hexagonal box
$E(q)$ with the trivial box $E(1)$ in the construction of $\Gamma_3(q)$.
An easy computation shows that this is not the case. More in general,
let us take a knot $K_1$ with a diagram $\Dd$.
Let us ``double'' $\Dd$ by replacing each arc of $\Dd$ with a band, and let
us add an isthmus in the most trivial possible way, \emph{i.e.}~by adding an arc
which is properly embedded in a small portion of a band, thus getting a (hc)--spine $\Gamma$.

If the blackboard framing defined by $\Dd$ is equal
to the trivial framing of $K$, then $\Gamma$ is a boundary spine of $H(\Gamma)$. However,
it is not difficult to show that for every prime $p$ we have $\Phi_p (H(\Gamma))=2+t^{h_p}$,
where $h_p$ is an integer depending on $p$. In particular, the form of $\Phi_p (H(\Gamma))$
does not allow us to use Lemma~\ref{2-K-obs} in order to conclude
that $H(\Gamma)$ is $(2)_L$-knotted.

In the following Subsections, therefore, we slightly modify our strategy
in order to get the desired handlebodies admitting a boundary spine but
no split constituent link.

One could also wonder if $H_3 (q)$ itself is indeed $(3)_S$-knotted
(\emph{i.e.}~if it does not admit \emph{any} boundary spine). 
We prove that this is the case in Proposition~\ref{alternative}.

\subsection{The tangle $B$}
\begin{figure}[h]
\begin{center}
\input{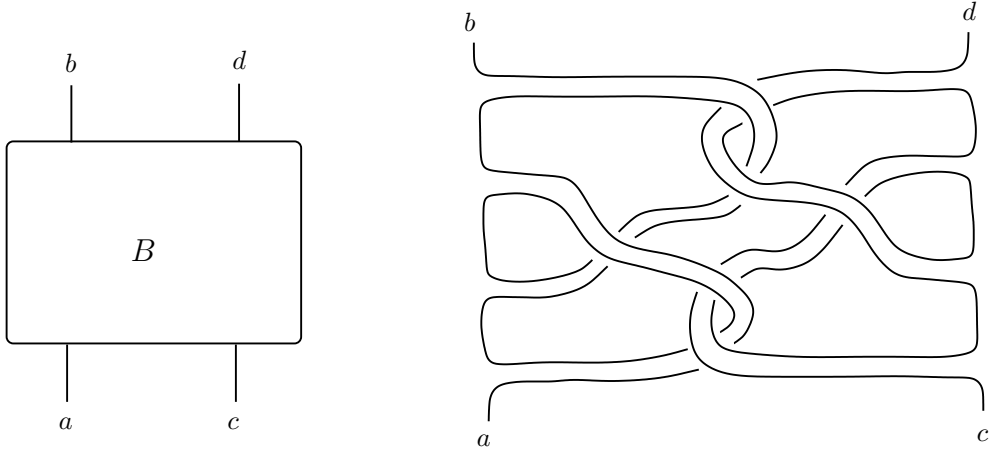}
\caption{\label{bande2:fig} 
The tangle $B$.}
\end{center}
\end{figure}

Let us now consider the tangle $B$ showed in Figure~\ref{bande2:fig}. 
Under the assumption that every arc is labelled with the non-trivial $\mathbb{Z}_2$
cycle,
we would like to compute the number $C_p(a,b,c,d)$ of $R_p$--colorings of $B$ which extend
the colors $a,b,c,d$ assigned on the ``corners'' of the diagram.

\begin{figure}[h]
\begin{center}
\input{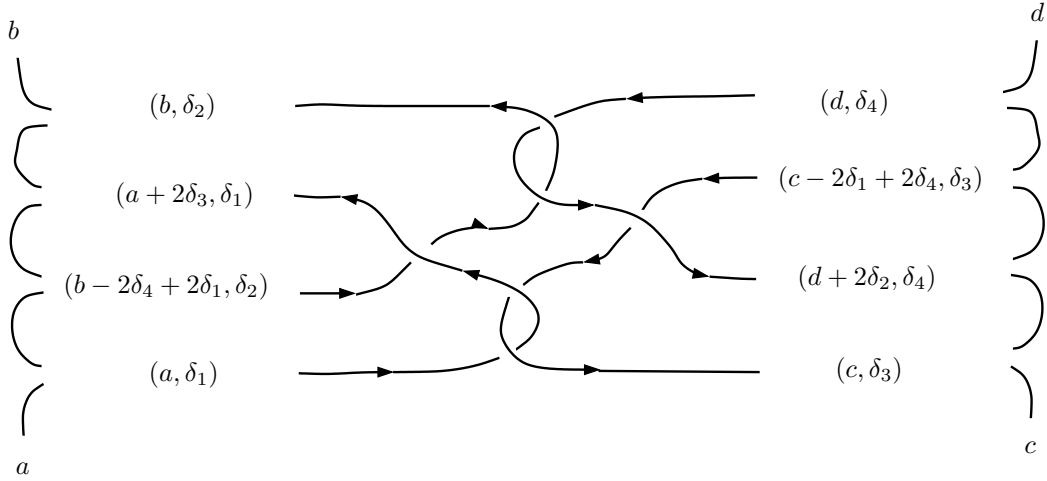}
\caption{\label{bande3:fig} 
Colorings of $B$. Every arrow represents a band.}
\end{center}
\end{figure}

\begin{lem}\label{Blemma}
 We have 
$$
\left\{
\begin{array}{ll}
C_p(a,b,c,d)=p^2\quad &{\rm if}\ a=b,\ c=d\ {\rm and}\ p=3\\
C_p(a,b,c,d)=1\quad &{\rm if}\ a=b,\ c=d\ {\rm and}\ p\neq 3\\
C_p(a,b,c,d)=0\quad & {\rm otherwise}\, .
\end{array}\right.
$$
\end{lem}
\begin{proof}
Let us orient the bands of $B$ as in Figure~\ref{bande3:fig}.
The condition on the colors of the arcs at the corners of $B$ implies
than the bands arriving at the corners of $B$ have to be labelled by the pairs
$(a,\delta_1)$, $(b,\delta_2)$, $(c,\delta_3)$ and $(d,\delta_4)$, where
$\delta_i\in R_p$ for every $i=1,2,3,4$. The discussion above shows that, due to the crossings of the bands, the pairs labelling the bands have to propagate
as described in Figure~\ref{bande3:fig}. Since the arcs at the ends of the bands 
join in pairs as described on the sides of the figure, 
the input coloring $(a,b,c,d)$
can be extended to a coloring of the whole $B$ if and only if the following conditions hold:
$$
\begin{array}{lll}
a+\delta_1=b-2\delta_4+2\delta_1,\ &
b-2\delta_4+2\delta_1+\delta_2=a+2\delta_3+\delta_1,\ &
a+2\delta_3=b+\delta_2,\\
d+\delta_4=c-2\delta_1+2\delta_4,\ &
c-2\delta_1+2\delta_4+\delta_3=d+2\delta_2+\delta_4,\ &
d+2\delta_2=c+\delta_3 .
\end{array}
$$
Such conditions may be rewritten as follows:
$$
\begin{array}{lllllll}
 a-b&=&\delta_1-2\delta_4&=&\delta_1-2\delta_4+\delta_2-2\delta_3&=&\delta_2-2\delta_3\\
d-c&=&\delta_4-2\delta_1&=&\delta_4-2\delta_1+\delta_3-2\delta_2&=&\delta_3-2\delta_2.
\end{array}
$$
 This readily implies that $C_p(a,b,c,d)=0$ whenever $a\neq b$ or $c\neq d$. 

Let us suppose that $a=b$ and $c=d$.
In this case, our conditions are equivalent to the equations
$$
\delta_1=2\delta_4,\quad \delta_4=2\delta_1,\quad \delta_2=2\delta_3,\quad \delta_3=2\delta_2.
$$
If $p\neq 3$, it is readily seen that this implies $\delta_i=0$ for every $i=1,\ldots,4$,
so $C_p (a,a,c,c)=1$. On the other hand, if $p=3$ then the conditions above 
are equivalent to $\delta_1=-\delta_4$ and $\delta_2=-\delta_3$, so the desired
colorings bijectively correspond to the choices of $(\delta_1,\delta_2)\in R_3^2$,
whence the conclusion. 
\end{proof}

\subsection{$(2)_L$-knotting does not imply $(3)_S$-knotting}\label{2L2S}
We are now ready to construct examples of spatial handlebodies
which are $(2)_L$-knotted (\emph{i.e.}~they do not admit a spine
with split constituent link) but $(3)_S$-unknotted
(\emph{i.e.}~they admit a boundary spine).

\begin{figure}[h]
\begin{center}
\input{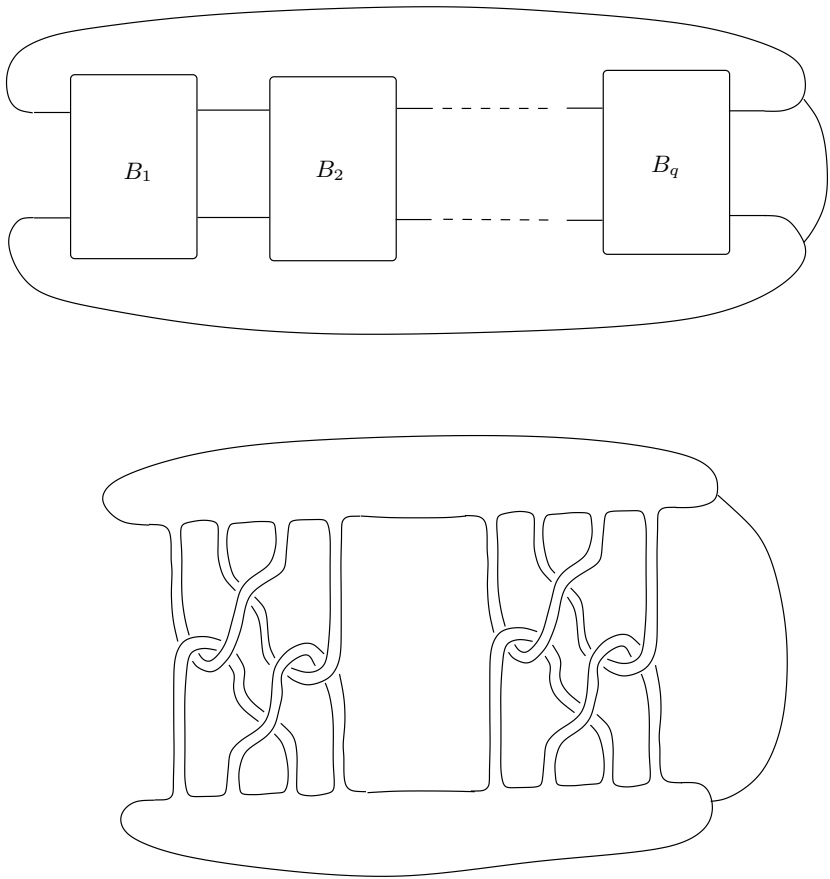}
\caption{\label{bande4:fig} 
On the top: the spine $\Gamma_4 (q)$. Every $B_i$ is a copy of the tangle $B$.
On the bottom: the case $q=2$.}
\end{center}
\end{figure}

For every $q\geq 1$, let $\Gamma_4 (q)$ be the graph described in Figure~\ref{bande4:fig},
and let us set $H_4 (q)=H(\Gamma_4(q))$.
It is obvious from the picture that $\Gamma_4 (q)$ is a boundary spine of
$H_4(q)$, so $H_4 (q)$ is $(3)_S$-unknotted for every $q$. On the other hand, we have the following:

\begin{prop}\label{H4}
For every $q\geq 1$ we have
$$
\Phi_3 (H_4(q))=3+t^{2q}.
$$
In particular,  for every $q\geq 1$ the handlebody $H_4 (q)$ is $(2)_L$-knotted.
Moreover, 
if $q'\neq q$, then the handlebodies
$H_4(q)$ and $H_4(q')$ are not isotopic.
\end{prop}
\begin{proof}
As usual, the only $(0,0)$-colorings of $\Gamma_4 (q)$ are the constant ones.

In order to describe the $(1,0)$-colorings of $\Gamma_4 (q)$, let us
first denote by $T(q)$ (resp.~$B(q)$) the constituent knot
of $\Gamma_4 (q)$ which lies on the top (resp.~on the bottom) of the picture.
It is immediate to observe that $T(q)$ and $B(q)$ are both trivial.
Let us now concentrate on $(1,0)$-colorings of $\Gamma_4 (q)$, where we suppose
for example that the cycle vanishes on $B(q)$.
Under this assumption, 
it is immediate to realize that the $(1,0)$-colorings
of $\Gamma_4(q)$ restrict to $1$-colorings of $T(q)$.
Since such knot is trivial, this implies that every $(1,0)$-coloring
of $\Gamma_4 (q)$ is constant on $T(q)$. 

Together with the discussion in Subsection~\ref{bands:sub}, this also implies
that the colorings of $B(q)$ are not affected by the crossings between the bands of $B(q)$
and the bands of $T(q)$. Then, every $(1,0)$-coloring
of $\Gamma_4 (q)$ restricts to a $0$-coloring (\emph{i.e.}~to a constant coloring)
of $B(q)$. Since the (constant) colors of $T(q)$ and of $B(q)$
have to agree with the color of the isthmus, we can conclude that
the only $(1,0)$-colorings of $\Gamma_4 (q)$ are the constant ones.
The same is true (by the very same argument) also for $(0,1)$-colorings,
so we have already proved that $\Phi_3 (H_4(q))=3+t^\alpha$, 
where $3^{\alpha+1}$ is equal to the number of $(1,1)$-colorings
of $\Gamma_4 (q)$.

Now, let us compute the number of $(1,1)$-colorings of $\Gamma_4(q)$
that induce the color $a\in R_3$ on the isthmus.
It is an immediate consequence of Lemma~\ref{H4} that this number is equal
to the $q$-times product of $C_3 (a,a,a,a)=3^2$. 
Since $a$ can be chosen in $3$ different ways, it readily follows
that the number of $(1,1)$-colorings
of $\Gamma_4 (q)$ is equal to $3^{2q+1}$, so that
$$
\Phi_3 (H_4(q))=3+t^{2q}.
$$

By Lemma~\ref{2-K-obs}, we have that $H_4(q)$ is $(2)_L$-knotted for every $q\geq 1$.
Moreover, if $q\neq q'$ we have $\Phi_3 (H_4(q))\neq \Phi_3 (H_4(q'))$,
and this implies that $H_4(q)$ is not isotopic to $H_4(q')$.
\end{proof}

\section{Handlebody patterns and the maximal free covering}\label{more-bpm} 
The main goal of this Section is to
provide rather handy combinatorial/topological characterizations of the
handlebody complements that admit a $\partial$-connected (resp.~$\partial_R$-connected)
cut system. 
Following Jaco~\cite{jaco}, we observe
in Proposition~\ref{jaco:gen} that the existence 
of a $\partial$-connected or $\partial_R$-connected cut system for 
a handlebody complement $M$ is related to
the way in which $\pi_1 (\partial M)$ sits inside $\pi_1 (M)$. 
Then, in order to study the image of $\pi_1(\partial M)$ into $\pi_1 (M)$ 
we extend 
some techinques coming from the theory of homology boundary links
to the context
of spatial handlebodies.

After recalling the definition of link pattern, we define the analogous notion
of \emph{handlebody pattern}. 
In Proposition~\ref{criterio2}
we exploit this notion for constructing
easily computable obstructions that allow to decide, starting from an explicitely given
cut system for $M$, if
$M$ admits a $\partial$-connected cut system.
 As an application, in Proposition~\ref{alternative} we provide a proof of the fact
that the handlebodies $H_1(p)$ introduced in Subsection~\ref{1Svs1L} are
$(3)_S$--knotted. For $p=3$, this fact was already known to Lambert~\cite{lambert}.
In our opinion, the characterization described in Proposition~\ref{criterio2}
is easier to handle in comparison, for
instance, with the original Lambert's topological treatment of his
example. 

Building on Proposition~\ref{jaco:gen}, we also describe a group--theoretic
obstruction for $M$ to admit a $\partial_R$-connected cut system. This obstruction
is then exploited in Proposition~\ref{nuovaprop} for proving
that there exist $(4)_L$-unknotted handlebodies whose complement
does not admit any $\partial_R$-connected cut system. As a consequence,
the notions of $(4)_L$-knotting and
$(3)_L$-knotting are not equivalent.

Jaco's obstruction for $M$ to admit a $\partial$-connected cut system
admits a nice topological
interpretation in terms of the 
\emph{maximal free covering} $\widetilde{M}_\omega$ of $M$. 
At the end of the Section we introduce such a covering, and
we prove 
that the boundary of $\widetilde{M}_\omega$ is connected if and only if
$M$ admits a $\partial$-connected cut system. 
We also show that the study of the first homology groups of $\widetilde{M}_\omega$
and $\partial\widetilde{M}_\omega$ provides other
obstructions for $M$ to admit $\partial$-connected (or~$\partial_R$-connected)
cut systems.

\subsection{Cut systems and epimorphisms of the fundamental group}\label{cuts:sub}
Let $M$ be the complement in $S^3$ either of a genus--2 handlebody or of
a $2$--component link, set $G=\pi_1 (M)$ and 
let $F_2=F(t_1,t_2)$ be the free group on two generators $t_1,t_2$.
Recall that a necessary and sufficient condition for $M$ to admit
a cut system is that there exists an epimorphism $\varphi\colon \pi_1 (M)\to F_2$. 
More precisely, if $\calS=\{S_1,S_2\}$ is a cut system for $M$, then
we can fix a basepoint $x_0\in M\setminus \calS$, and 
define an epimorphism $\varphi\colon G=\pi_1 (M,x_0)\to
F_2$ in such a way that, if $g\in G$ is represented by a loop
disjoint from $S_2$ (resp.~$S_1$) and intersecting positively $S_1$ (resp.~$S_2$) in exactly one
point, then $\varphi (g)=t_1$ (resp.~$\varphi (g)=t_2$). We say that such 
a $\varphi$ is \emph{associated to $\calS$}.

The following result is due to Stallings~\cite{stall} and shows that, up to post--compositions with automorphisms
of $F_2$, there exists a unique epimorphism from $G$ to $F_2$.
Define $G_1=[G,G]$, $G_{n+1}= [G_n,G]$, 
$G_\omega = \bigcap_n G_n$.  

\begin{teo}[\cite{stall}]\label{stallings}
  Suppose $G=\pi_1 (M,x_0)$ is the fundamental group of $M$,
  and let $\varphi\colon G\to F_2$ be any epimorphism.  Then $\ker
  \varphi=G_\omega$.
\end{teo}

Henceforth we tacitly make the assumption that every connected
component of the boundary of a given cut system is {\it essential},
\emph{i.e.}~it
does not bound a disk on $\partial M$. Of course, every
cut system for $M$ can be compressed in order to satisfy this
requirement.

\subsection{Link patterns}
Keeping notations from Subsection~\ref{cuts:sub}, let us now specialize
to the case when $M=\compl (L)$ is the complement of an (ordered and oriented)
homology boundary link $L$. 

As elements of $G$, the meridians of $L$ are defined only up to conjugacy.
For $i=1,2$,
let $\gamma_i\in G$ be a representative of the $i$--th meridian of $L$, and
set $w_i=\varphi (\gamma_i)$. Then $w_i$ is well--defined  up to conjugacy in
$F_2$.  Adding to $M$ two 2-handles along the meridians of $L$ we obtain a space
homeomorphic to $S^2\times [0,1]$, which is simply connected, 
so $w_1,w_2$ normally generate $G$
(\emph{i.e.}~there do not exist proper normal subgroups of $G$ containing 
$w_1$ and $w_2$). 
Also recall that any epimorphism
$\psi\colon G\to F_2$ 
is obtained from $\varphi$ by post--composition
with an automorphism of $F_2$. 

An old result by Nielsen (see \emph{e.g.}~\cite{MKS} for a proof)
ensures that every $n$--uple of generators of the free group $F_n$ of rank $n$
is in fact a set of free generators of $F_n$. Such an $n$--uple is called
a \emph{base} of $F_n$.
The following definition is taken from~\cite{CL}. 

\begin{defi}
{\rm
A \emph{link pattern} is a pair $(w_1,w_2)\in F_2\times F_2$ such that 
$w_1$ and $w_2$ normally generate $F_2$.
The pattern $(w_1,w_2)$ is realized by the link $L$ if there exist
an epimorphism $\varphi\colon \pi_1 (\compl(L))\to F_2$ and a choice
of meridians $\gamma_1,\gamma_2$ such that $\varphi (\gamma_i)=w_i$ for $i=1,2$.
Two link patterns $(w_1,w_2)$ and $(w_1',w_2')$ are \emph{equivalent} if 
there exist $h_1,h_2\in F_2$  
and  $\alpha\in {\rm Aut} (F_2)$ such that
$w'_i=h_i\alpha (w_i)h_i^{-1}$. 
A pattern is \emph{trivial} if it is equivalent to a base
of $F_2$.}  
\end{defi}

The discussion above shows that to any homology boundary link there is associated
a well-defined equivalence class of link patterns. Moreover, it is proved in~\cite{CL} that
every pattern is realized by a homology boundary link (see also~\cite{bellis}
for an explicit construction).

\begin{remark}\label{base:rem}
{\rm Of course every two bases of $F_2$ are equivalent as link patterns, but
if $(t_1,t_2)$ is a base of $F_2$, then 
the pair $(t_1,wt_2w^{-1})$, while being trivial as a pattern,
is not necessarily a base of $F_2$.
More precisely, let us show that 
$(t_1,wt_2 w^{-1})$ is a base of $F_2$ if and only if $w=t_1^nt_2^m$ 
for some $m,n\in\mathbb{Z}$ (we will need this result later). 

Of course, if $w=t_1^nt_2^m$, then $t_1$ and $wt_2 w^{-1}=t_1^nt_2t_1^{-n}$
generate the whole $F_2$, so $(t_1,wt_2 w^{-1})$ is a base of $F_2$.
On the other hand, 
let us suppose that 
$(t_1,wt_2 w^{-1})$ is a base of $F_2$, and
let us choose $n\in\mathbb{Z}$ in such a way that
$w'=t_1^{-n}w$ either is the identity,
or is represented by a reduced word starting with the symbol
$t_2$ or $t_2^{-1}$. Observe that $(t_1,w't_2(w')^{-1})$ is also a base of $F_2$,
so there exists an element $R(a,b)$ in the free group over two generators $F(a,b)$ such
than $R(t_1,w't_2(w')^{-1})=t_2$.
Let now $w''$ be the reduced word representing $w't_2(w')^{-1}$. Then it is easily seen
that in any product of the form $t_1^{\pm 1} (w'')^{\pm 1}$ there cannot be cancellations.
It is easily seen that this forces $R(a,b)=b$, whence $w''=t_2$, and $w'=t_2^m$ for some
$m$. We have therefore $ w=t_1^nt_2^m$, as claimed.} 
\end{remark}

The following result, which was already observed by Smythe in~\cite{Smythe},
characterizes in terms of patterns those homology boundary links which are in fact boundary links.

\begin{prop}\label{pattern1}
 A homology boundary link is a boundary link if and only if its associated
link patterns are trivial.
\end{prop}

Recall that an element $w\in F_2$ is \emph{primitive} if it is an element 
of a base of $F_2$. There is an extensive literature about primitive elements
in the free group on $n$ generators, and a particular interest has been devoted
to the case of rank two. The following Lemma provides an useful
characterization of trivial link patterns:

\begin{lem}\label{free1:lemma}
Let $F(t_1,t_2)=F_2$ be the free group on two generators $t_1,t_2$. Then:
\begin{enumerate}
\item 
A link pattern $(w_1,w_2)\in F_2\times F_2$ 
is trivial if and only 
if $w_1$ and $w_2$ are both primitives. 
\item
Suppose that $$w=t_1^{\alpha_1} t_2^{\beta_1}t_1^{\alpha_2}t_2^{\beta_2}\ldots t_1^{\alpha_m}t_2^{\beta_m}$$ 
is a cyclically reduced word representing a primitive element, where $\alpha_i\neq 0$,
$\beta_i\neq 0$ for every $i=1,\ldots,m$. Then, all the $\alpha_i$'s share the same sign,
and all the $\beta_i$'s share the same sign.
\end{enumerate}
\end{lem}
\begin{proof}
Since $(w_1,w_2)$ normally generate $F_2$, they project
onto a base of $F_2/[F_2,F_2]=\mathbb{Z}^2$. Therefore, point~(1)
follows from~\cite[page 167]{hibook}.

Point~(2) dates back to Nielsen~\cite{nielsen} (see \emph{e.g.}~\cite{wei},
\cite{met}, \cite{os}
for alternative proofs).
\end{proof}

\subsection{Jaco's characterization of handlebody complements
admitting $\partial$-connected cut systems}\label{nuovenote}
Let us now consider the case when $M$ is the complement
of a spatial handlebody $H$
such that $\cut (M)=2$.
Let $\Ss=\{S_1,S_2\}$ be any cut system of $M$, and set $G=\pi_1(M,x_0)$, where $x_0$ is a basepoint such that
$x_0\in \partial M\setminus (S_1\cup S_2)$.
Let $F_2=F(t_1,t_2)$ be the free group on two generators $t_1,t_2$, and
let $\varphi\colon G\to
F_2$ be the epimorphism associated to $\calS$.

Let us set $G_\partial =\pi_1 (\partial M,x_0)$ and denote
by 
$i\colon \partial M\to M$ the inclusion. The following result 
is proved
in~\cite[Theorems 2 and 3]{jaco}, and in some sense extends Proposition~\ref{pattern1}
to the case of spatial handlebodies.

\begin{prop}\label{jaco:prop}
The manifold $M$ admits a
 $(M\to W)$--boundary--preserving--map (or, equivalently, a $\partial$-connected
cut system) if and only if
$\varphi(i_\ast(G_\partial))=F_2$, \emph{i.e.}~if and only if
$i_\ast(G_\partial)$ surjectively projects onto $G/G_\omega$.
\end{prop}

In what follows we  show how Proposition~\ref{jaco:prop} can be exploited
for proving that the handlebodies $H_1 (p)$ 
introduced in Subsection~\ref{1Svs1L} are
$(3)_S$--knotted. Moreover, in Proposition~\ref{jaco:gen}
we extend Proposition~\ref{jaco:prop} in order to obtain an obstruction for a 
handlebody complement to admit a $\partial_R$-connected cut system.

\subsection{Cut systems and longitudes}\label{longitude}
Suppose that $X= \compl (L)$ is the complement of a 2--component homology
boundary link and let $\calS$ be a cut system for $X$. The boundary components
of the surfaces of $\calS$ belong to two families of parallel curves, one on each component
of $\partial X$. It is well--known that the isotopy classes of such curves on $\partial X$ 
do not depend on the particular cut system $\calS$ (see also
Subsection~\ref{longitude2}). These isotopy classes define the
\emph{longitudes} of $L$.  

We would like to extend this notion
to the case when $M=\compl (H)$ is the complement
of a handlebody $H$ such that $\cut (M)=2$. It turns out that 
in this case the definition of longitudes is less obvious, and longitudes are in fact
no more independent of the choice of a cut system (see Remark~\ref{nocommon}).

\begin{lem}\label{long:lemma}
  Let $\calS$ be a cut system for $M$.  There exist three disjoint
  essential simple closed curves $\ell_1,\ell_2,\ell_3$ on $\partial
  M$ such that the following conditions hold:
\begin{enumerate}
\item each component of the reduced boundary of $\calS$ is parallel to
  $\ell_1$, $\ell_2$ or $\ell_3$;
\item
the union $\ell_1\cup \ell_2$ is not separating in $\partial M$; 
\item
$\ell_1$ and $\ell_2$ (provided with some orientations) give a basis
of $\ker \left(i_\ast\colon H_1 (\partial M)\to H_1 (M)\right)$;
\item if $\partial \Ss \setminus \partial_R \Ss$ is non--empty, then
  each component of the reduced boundary of $\calS$ is parallel to
  $\ell_1$ or $\ell_2$, and each further component of $\partial \Ss$ 
is parallel to $\ell_3$.
\end{enumerate}
\end{lem}
\begin{proof}
  Since on $\partial M$ there exist at most three non--parallel
  non--trivial unoriented simple loops, there exist three loops
  $\ell_1,\ell_2,\ell_3$ such that every (unoriented) loop in
  $\partial\calS$ is parallel  to some $\ell_i$,
  $i=1,2,3$. Moreover, we may suppose that $\ell_1\cup \ell_2$ does
  not separate $\partial M$, so that $\ell_1,\ell_2,\ell_3$ may be
  oriented in such a way that either $[\ell_3]=[\ell_1]+[\ell_2]$ (if
  $\ell_3$ does not separate $\partial M$) or $[\ell_3]=0$ (if
  $\ell_3$ separates $\partial M$). In any case the submodule of
  $H_1 (\partial M)$ generated by the homology classes $[\partial
  S_1]\in H_1 (\partial M)$ and $[\partial S_2]\in H_1 (\partial M)$
  is contained in $\langle [\ell_1],[\ell_2]\rangle$. It is proved for
  example in~\cite{RRR} that $\langle [\partial S_1],[\partial
  S_2]\rangle$ has rank 2, is equal to $\ker i_\ast$ and is not a
  proper finite--index submodule of any submodule of $H_1 (\partial
  M)$. These facts easily imply (3), and (4) is obvious.
\end{proof}

\begin{figure}
 \begin{center}
  \input{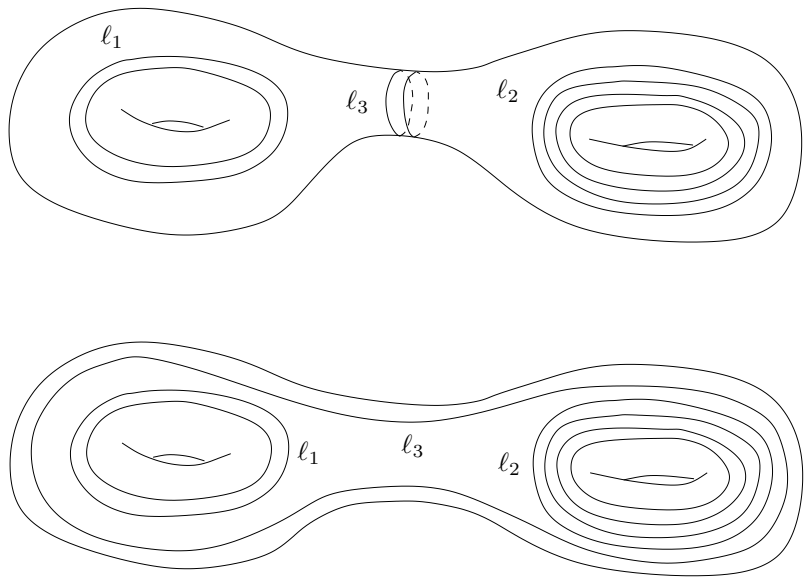}
\caption{The boundary components of a cut system $\calS$. On the top: $\ell_3$ separates $\partial M$, so $\calS$ is good. On the bottom: $\ell_3$ does not separate $\partial M$, and $\calS$ is not good.}\label{cutsystems:fig}
 \end{center}
\end{figure}

\begin{defi}
Let $\calS$, $\ell_1$ and $\ell_2$ be as in the statement of Lemma~\ref{long:lemma}.
Then we say that (the isotopy classes of) $\ell_1,\ell_2$ are a pair of \emph{longitudes} associated
to $\calS$.  
\end{defi}

If $\partial \calS\setminus \partial_R\calS\neq \emptyset$, then the longitudes
of $\calS$ are uniquely determined by $\calS$ (see the top of Figure~\ref{cutsystems:fig}). On the contrary, if the
connected components of $\partial \calS$ are divided into three
non--isotopic families of parallel curves each of which is
non--separating, then $\calS$ defines three pairs of longitudes (see the bottom
of Figure~\ref{cutsystems:fig}). In Lemma~\ref{good} we show how to get rid of this ambiguity, which anyway is not really relevant to our purposes.

\medskip

\begin{defi}
 {\rm
We say that the a cut system $\calS$ is \emph{good} if the components
of $\partial_R \calS$ fall into two isotopy classes
of curves on $\partial M$ (and in this case such classes define the unique pair of longitudes
associated to $\calS$). 
%It is \emph{very good} if it is good and satisfies
%$\partial\calS=\partial_R \calS$. Equivalently, $\calS$ is very good if and only
%if the components
%of $\partial \calS$ fall into two isotopy classes
%of curves on $\partial M$.
}
\end{defi}

\begin{lem}\label{good}
Let 
$\ell_1,\ell_2$ be a fixed pair of longitudes of the cut system $\calS$. Then
$M$ admits a good cut system 
$\calS'$ with longitudes $\ell_1,\ell_2$.
If $\calS$ is $\partial_R$-connected, we may set $\calS'=\calS$.
\end{lem}
\begin{proof}
If every connected component of $\calS$ is parallel to $\ell_1$ or
  to $\ell_2$, then we may set $\calS'=\calS$, and we are done.
  Otherwise, let $\ell_3$ be the loop defined in
  Lemma~\ref{long:lemma}. 
If $\ell_3$ is separating
 (and this is the case,
in particular, if $\calS$ is $\partial_R$-connected), then we may set $\calS'=\calS$.
Otherwise, $\ell_3$ separates $\partial
  H\setminus (\ell_1\cup \ell_2)$ in two pair of pants $Y_1$ and
  $Y_2$ with $\partial Y_1=\partial Y_2=\ell_1\cup\ell_2\cup\ell_3$.  
If $\partial \calS$ has $n$
  components parallel to $\ell_3$, we define $\calS'$ by replacing
  small neighbourhoods in $\calS$ of such components with $n$ parallel
  copies of $Y_1$ (or of $Y_2$). 
\end{proof}

\subsection{Handlebody patterns}\label{handlepat}
Suppose now that $\calS$ is a good cut system for $M$, and fix
an ordering and an auxiliary orientation on the longitudes $\ell_1,\ell_2$
of $\calS$. 

We now define three elements $w_0 (\calS)$, $w_1 (\calS)$, $w_2 (\calS)$
of $F_2=F(t_1,t_2)$ as follows. For $i=1,2$, let $W_i$ be the non--annular component of $\partial M\setminus \partial \calS$ whose closure contains (a loop isotopic to) $\ell_i$ 
(if $\partial\calS=\partial_R\calS$, then $W_1=W_2$ is homeomophic to
a $4$--punctured sphere, otherwise both $W_1$ and $W_2$ are homeomorphic to
$3$--punctured spheres).
Take basepoints 
$x_i\in W_i$, $i=1,2$, and recall that every component of $\partial\calS$
inherits a well-defined orientation induced by the orientations of $S_1$, $S_2$.

\begin{figure}
 \begin{center}
  \input{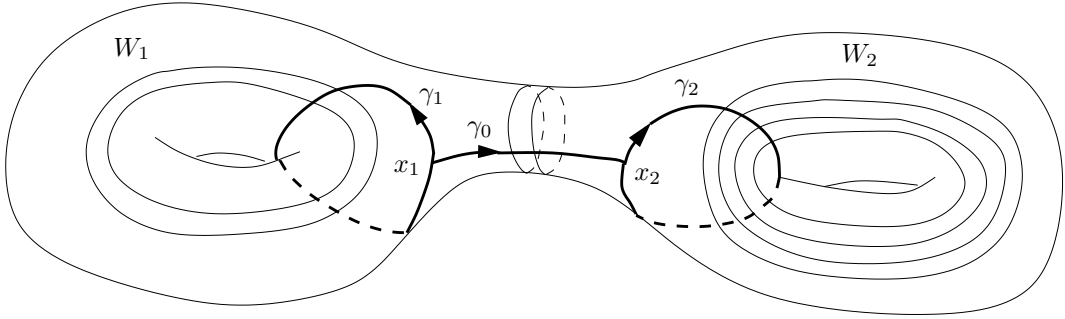}
\caption{The definition of handlebody pattern.}\label{tree2:fig}
 \end{center}
\end{figure}

For $i=1,2$ we fix a simple oriented loop $\gamma_i$ on $\partial M$
such that the following conditions hold (see Figure~\ref{tree2:fig}):
\begin{itemize}
\item 
$\gamma_i$ is based at $x_i$ and transverse to $\partial \calS$; 
\item
$\gamma_i$ is disjoint from every separating
connected component of $\calS$;
\item
$\gamma_i$ is disjoint from every component of $\partial \calS$
isotopic to $\ell_j$ for $j\neq i$;
\item 
$\gamma_i$ positively intersects a representative of $\ell_i$ exactly in one point;
\item
$\gamma_i$ transversely intersects every component of $\partial\calS$
isotopic to $\ell_i$ exactly in one point.
\end{itemize}
Starting and ending in $x_1$, we now follow $\gamma_1$
 and write down a letter $t_i$ (resp.~$t_i^{-1}$)
every time $\gamma_1$ positively (resp.~negatively) intersects a component of 
$\partial S_i$, thus obtaining a (not necessarily reduced)
word $\widetilde{w}_1(\calS)$ that represents the element
$w_1 (\calS)\in F_2$. The word $\widetilde{w}_2(\calS)$ and the element $w_2(\calS)$ are obtained by applying the very same procedure
to $\gamma_2$. Finally, we take an arc $\gamma_0$ starting at $x_1$, ending at $x_2$,
and transversely intersecting each component of $\partial\calS\setminus \partial_R\calS$
exactly in ont point,
and we denote by $\widetilde{w}_0(\calS)$ the word obtained by associating the element
$t_i$ (resp.~$t_i^{-1}$) to every positive (resp.~negative) intersection of
$\gamma_0$ with $\partial S_i$. We also denote by $w_0 (\calS)$ the element
of $F_2$ associated to $\widetilde{w}_0 (\calS)$.
Observe that $w_i (\calS)$ depends both on the orientations
of $S_1$, $S_2$ (which is part of the datum $\calS$), and on the fixed
auxiliary orientations
on $\ell_1$, $\ell_2$. However, our notation forgets 
about this last dependence, since it is not relevant to our purposes.

Our next aim is to describe as explicitely as possible the relations
between the topological properties of the cut system $\calS$
and the algebraic properties of the epimorphism $\varphi$ associated to $\calS$. 
These last properties are encoded by the $w_i(\calS)$'s, while the 
words $\widetilde{w}_i(\calS)$'s keep track of the actual 
components of $\partial\calS$. However, we will see in Lemma~\ref{simpli} below
that the $w_i(\calS)$'s encode in some sense
the ``essential'' information about $\partial\calS$.
\smallskip

Let us now come back to the notation of Subsection~\ref{nuovenote},
suppose that $\calS$ is good, and set $x_0=x_1$, \emph{i.e.}~let
$G=\pi_1 (M,x_1)$ and $G_\partial =\pi_1(\partial M,x_1)$.

\begin{lem}\label{generatori}
Let $w_0 (\calS)$, $w_1(\calS)$, $w_2(\calS)$
be as above. Then:
\begin{enumerate}
\item 
 $\varphi (i_\ast(G_\partial))\subseteq F_2$ is generated
by the elements $w_1 (\calS)$ and $w_0(\calS)w_2(\calS)w_0(\calS)^{-1}$.  
\item
The pair $(w_1(\calS),w_2(\calS))$ normally generates $F_2$.
\item
The kernel of the map $\varphi\circ i_\ast\colon G_\partial\to F_2$
is normally generated by any pair of loops 
which are freely homotopic to the longitudes of $\calS$.
\end{enumerate}
\end{lem}
\begin{proof}
(1): For $i=1,2$, let
$l_1\subseteq \partial M$ be a loop based at $x_i$ which is isotopic to $\ell_i$
and disjoint from $\partial \calS$. Then, the group $G_\partial$ is generated by the elements $$
g_1=[\gamma_1], \quad g_2=[\gamma_0\ast \gamma_2\ast \gamma_0^{-1}],\quad n_1=[l_1],\quad n_2=[\gamma_0\ast
l_2\ast \gamma_0^{-1}]\ .
$$
Since $\varphi$ is associated to $\calS$, we have
\begin{equation}\label{kerphi}
\varphi (i_\ast(g_1))=w_1 (\calS),\
\varphi (i_\ast(g_2))=w_0(\calS)w_2(\calS)w_0(\calS)^{-1},\ 
\varphi (i_\ast(n_1))=\varphi(i_\ast(n_2))=1,
\end{equation}
 whence point~(1).

(2): By point~(1), it is sufficient to prove that $i_\ast(G_\partial)$ normally generates $G$.
Let $m_1,m_2$ be simple closed loops of $\partial M$ which bound
disjoint compressing disks $D_1,D_2$ in the handlebody $H$ in such a way that
$D_1\cup D_2$ does not separate $H$. 
The smallest normal subgroup $N$ of $G$ containing $i_\ast(G_\partial)$ also contains
the (conjugacy classes of) 
two elements $h_1,h_2$ which are freely homotopic to $m_1$ and $m_2$ in $\partial M$.
But the fundamental group of $M\cup D_1\cup D_2$ is obviously trivial, so an easy application of Seifert--Van Kampen Theorem implies that the smallest normal subgroup of $G$
containing $h_1,h_2$ coincides with $G$. We have \emph{a fortiori} $N=G$, whence
point~(2).

(3): 
Take $g\in G_\partial$,
and let $g_1,g_2,n_1,n_2$ be the generators of $G_\partial$ introduced 
in the proof of point~(1). There exists a word 
$R(a,b,c,d)$ in the symbols $a^{\pm 1},b^{\pm 1}, c^{\pm 1}, d^{\pm 1}$ such that $g=R(g_1,g_2,n_1,n_2)$. Since $\varphi(i_\ast(n_i))=1$ for $i=1,2$, the element
$\varphi(i_\ast(g))\in F_2$ is represented
by the word $\overline{R}(\varphi(i_\ast(g_1)), \varphi(i_\ast(g_2)))$,
where $\overline{R}(a,b)=R(a,b,\emptyset, \emptyset)$ is obtained from $R$ by replacing each occurence
of $c^{\pm 1}$ and $d^{\pm 1}$ with the empty word.

It follows from point~(2)
that the elements $\varphi(i_\ast(g_1))$, $\varphi(i_\ast(g_2))$ freely generate
a rank--2 subgroup of $F_2$. This implies in turn that 
$\varphi(i_\ast(g))=1$ if and only if
the word $\overline{R}$
represents the trivial element of $F(a,b)$. This last condition is in turn equivalent
to the fact that $g$ belongs to the subgroup of $G_\partial$
normally generated by $n_1,n_2$, whence the conclusion.
\end{proof}

We are now ready to define the notion of handlebody pattern.

\begin{defi}\label{handle:patt}
{\rm A \emph{handlebody pattern} is a triple $(w_0,w_1,w_2)\in F_2\times F_2\times F_2$ such that
$w_1,w_2$ normally generate $F_2$ (\emph{i.e.}~$(w_1,w_2)$ is a link pattern). 
If $M$, $\calS$ and $\ell_1,\ell_2$ are as above, then we say that
$(w_0(\calS),w_1(\calS),w_2(\calS))$ is the pattern associated to
$\calS$. We also say that $(w_0(\calS),w_1(\calS),w_2(\calS))$
is realized by $M$ (or by $H$).

A pattern $(w_0,w_1,w_2)$ is \emph{trivial} if the pair $(w_1, w_0w_2w_0^{-1})$
is a base of $F_2$ (by Remark~\ref{base:rem}, this condition is strictly 
stronger than the condition that
$(w_1, w_2)$ is trivial as a link pattern).
} 
\end{defi}

\begin{remark}
 {\rm 
With some effort it is possible to define an equivalence relation on the set
of handlebody patterns, in such a way that a fixed 
$M$ uniquely defines an equivalence class of handlebody patterns.
Such equivalence relation is a bit more complicated than the one
defined on link patterns, and since we won't need to exploit the notion
of equivalent handlebody patterns, we are not discussing it here (however,
it is maybe worth mentioning that Lemma~\ref{simpli} below shows for
example that, with respect to this relation, the handlebody pattern $(w_0,w_1,w_2)$ should be equivalent to $(1,w_1,w_0w_2w_0^{-1})$). 

Moreover, 
putting together Lemma~\ref{translation2} with
the fact that every link pattern is realized by a 
homology boundary link 
it can be easily proved that every
handlebody pattern is realized by a spatial handlebody. 
}
\end{remark}

The following Lemma shows that patterns encode the relevant information
about the topology of the boundary of cut systems.

\begin{lem}\label{simpli}
Let $\calS$ be a good cut system for $M$ and let $g_1,g_2$ be elements in $F_2$. 
Then $M$ admits a good cut system $\calS'$ satisfying the following conditions:
\begin{enumerate}
 \item 
$\calS'$ has the same longitudes as $\calS$.
\item
$(w_0(\calS'),w_1(\calS'),w_2(\calS'))=(g_1w_0(\calS)g_2^{-1}, g_1w_1(\calS)g_1^{-1},
g_2w_2(\calS)g_2^{-1})$.
\item
for $i=1,2,3$, the word $\widetilde{w}_i(\calS')$ is reduced
(if $w_0(\calS')=1$, then it is understood that $\widetilde{w}_0 (\calS')$
is the empty word, \emph{i.e.}~that $\partial\calS=\partial_R\calS$).
\end{enumerate}
\end{lem}
\begin{proof}
 Let $\ell_1,\ell_2$ be the longitudes of $\calS=\{S_1,S_2\}$, and let 
us denote by $m_1,\ldots,m_l$ the components
of $\partial\calS\setminus \partial_R \calS$, and by $W_1$ (resp.~$W_2$)
the non--annular component of $\partial M\setminus \partial\calS$
whose boundary contains a loop isotopic to $\ell_1$ (resp.~$\ell_2$).  
If $l\geq 1$, we order the $m_i$'s
in such a way that $m_1$ (resp.~$m_l$) is a boundary component of $W_1$ (resp.~$W_2$),
and $m_i,m_{i+1}$ bound an annulus in $\partial M\setminus\partial\calS$
for $i=1,\ldots,m-1$.

Let us first consider the case
when $g_1=t_i^{\pm 1}$ and $g_2=1$. Let us fix an embedded arc $\alpha\colon [0,1]\to M$
satisfying the following properties: $\alpha(0)\in S_i$, $\alpha (1)\in W_1$
and $\alpha (t)\in M\setminus (\partial M\cup S_1\cup S_2)$ for every $t\in (0,1)$
(such an arc exists because $M\setminus (S_1\cup S_2)$ is connected).
Let $Z=D^2\times [0,1]\subseteq M$ be a $1$--handle satisfying the following
conditions (see Figure~\ref{costruzione:fig}): 
\begin{itemize}
\item 
$\alpha$ is the core of $Z$; 
\item 
$Z\cap (S_1\cup S_2)=Z\cap S_1=D^2\times \{0\}$
is a regular neighbourood of $\alpha (0)$ in $S_i\setminus \partial S_i$; 
\item
$Z\cap \partial M=Z\times \{1\}$
is a regular neighbourhood of $\alpha (1)$ in $W_1\setminus \partial W_1$. 
\end{itemize}

\begin{figure}
 \begin{center}
  \input{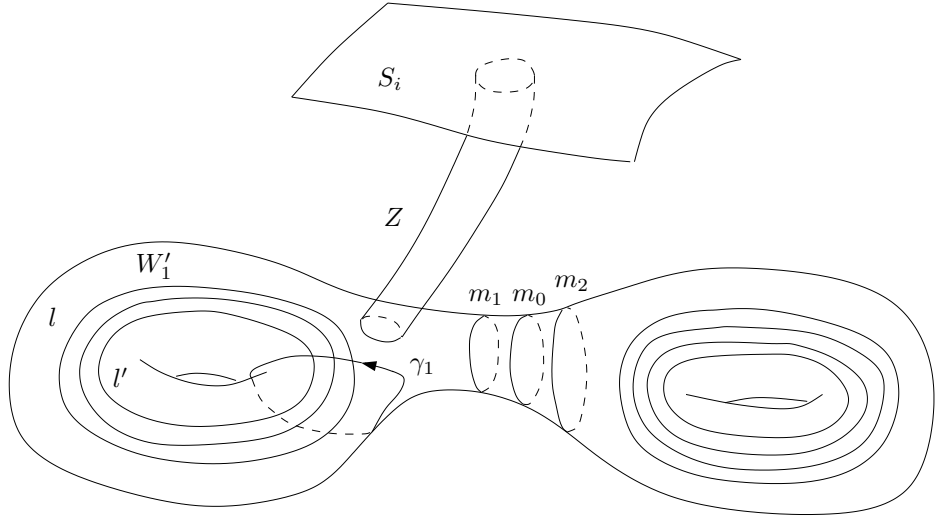}
\caption{The construction described in the proof
of Lemma~\ref{simpli}.}\label{costruzione:fig}
 \end{center}
\end{figure}

Let now $\gamma_0$ and $\gamma_1$ be the loop and the arc of $\partial M$
entering in the definition of the handlebody pattern associated to $\calS$. Let
$m_0\subseteq W_1$ be a simple loop parallel to $m_1$ which 
is disjoint from $\gamma_1$ and transversely intersects $\gamma_0$ exactly in one point, and let $l,l'\subseteq W_1$ be
two loops isotopic to the components of $\partial W_1\setminus m_1$.
We also assume that $\gamma_1$ transversely meets $l$ and $l'$ exactly in one point
and in this order. Finally, we denote by $W'_1\subseteq W_1$ the pair of pants bounded
by $m_0,l,l'$. Up to shrinking $Z$, we may suppose that $D^2\times \{0,1\}$
is contained in the internal part of $W'_1$. 
Let now $S'_i$ be defined by 
$$
S'_1=\left( S_1\cup (\partial D^2\times [0,1])\cup W'_1\right)\setminus 
\left({\rm int}(D^2)\times \{0,1\}\right)\ ,
$$
and endow $S'_i$ with the orientation induced by $S_i$.
It is clear that we may push ${\rm int}(S'_i)\cap \partial M$ slightly inside $M$ thus obtaining
a new cut system $\calS'=\{S'_i,S_j\}$ for $M$. By construction, as a set we have that
$\partial\calS'$ is obtained from $\calS$ by adding 
$m_0,l$ and $l'$. 
It is easily seen that if $\gamma_0$ intersects
$m_0$ positively (resp.~negatively), then $\gamma_1$ intersects $l$ positively (resp.~negatively) and $l'$ negatively (resp.~positively). 
Therefore, we get
$$
w_0(\calS')=t_i^{\pm 1}w_0(\calS),\quad 
w_1(\calS')=t_i^{\pm 1}w_1(\calS)t_i^{\mp 1},\quad
w_2(\calS')=w_2(\calS) \ .
$$
In order to get the desired exponent for the added factors $t_i^{\pm 1}$, it is sufficient to
replace $\alpha$, if necessary, with an arc enjoying the very same properties
as $\alpha$, but exiting from $S_i$ on the opposite side with respect to $\alpha$. 

The very same proof also shows how $\calS$ can be modified in order to obtain 
a cut system $\calS'$ such that
$$
w_0(\calS')=w_0(\calS)t_i^{\pm 1},\quad 
w_1(\calS')=w_1(\calS),\quad 
w_2(\calS')=t_i^{\mp 1}w_2(\calS)t_i^{\pm 1}\ .
$$

After repeating the construction just described a finite number of times, 
we obtain
a cut system $\calS'$ satisying conditions~(1) and~(2).

In order to conclude, it is now sufficient to show that we may replace $\calS'$
with a cut system having the same longitudes and the same pattern, and satisfying the additional property 
that its associated words are all 
reduced.
So, let us suppose that a string of the form $t_j^{\pm 1}t_j^{\mp 1}$
appears in $\widetilde{w}_i(\calS')$ for some $i=0,1,2$. Then, the components
$c,c'$ of $\partial \calS'$ corresponding to the symbols $t_j^{\pm 1}$, $t_j^{\mp 1}$
belong to the same surface $S_j$ of $\calS$, bound an annulus $A$
in $\partial M\setminus \partial \calS$,
and inherit from $S_j$ opposite
orientations. This implies that 
the surface $\widehat{S}_j$ obtained by cutting from $S_j$ 
small neighbourhoods of $c$ and $c'$ and adding the annulus $A'$ obtained by pushing
$A$ slightly inside $M$ is orientable, and may be given, therefore, the orientation
induced by $S_j$. Moreover, $\widehat{S}_j$ is obviously disjoint
from $S_2$, so we may modify $\calS'$ by replacing
$S_j$ with $\widehat{S}_j$. This operation has the effect of
cancelling out the string $t_j^{\pm 1}t_j^{\mp 1}$ from $\widetilde{w}_i (\calS')$.
After a finite number of operations of this type we end up with a cut system
satisfying all the properties of the statement.
\end{proof}

For later purposes we point out the following easy Corollary of the previous Lemma:

\begin{cor}\label{verygood}
 Let $\calS$ be a good cut system for $M$.
Then $M$ admits a good cut system $\calS'$ having the same longitudes
as $\calS$ and satisfying the additional property that
$\partial\calS=\partial_R\calS$.
\end{cor}
\begin{proof}
It is sufficient to apply Lemma~\ref{simpli} to the case
$g_1=w_0(\calS)^{-1}$, $g_2=1$.
\end{proof}

\subsection{$(4)_L$-knotting is equivalent to $(4)_S$-knotting}
As an application of Lemma~\ref{simpli} we obtain the following:

\begin{prop}\label{L4S4}
A spatial handlebody $H$ is $(4)_L$-knotted if and only if it is $(4)_S$-knotted.
\end{prop}
\begin{proof}
Let us prove that if $H$ is $(4)_L$-unknotted, then it is $(4)_S$-unknotted, the other implication being trivial.
So, let $\Gamma$ be a spine of $H$ such that the link $L_\Gamma$ is homology boundary,
and let $S_1,S_2$ be a pair of generalized Seifert surfaces for $L_\Gamma$.
We may suppose that $S_i$ is transverse to the isthmus of $\Gamma$ for $i=1,2$.
Then, up to shrinking $H$ onto a smaller neighbourhood of $\Gamma$,
the surfaces $S_1\cap M$, $S_1\cap M$ define a good cut system $\calS$ for $M$
such that the components of $\partial\calS\setminus\partial_R\calS$ bijectively
correspond to the points where $S_1\cup S_2$ intersects the isthmus
of $\Gamma$. 

By Corollary~\ref{verygood}, we may replace $\calS$ with a cut system $\calS'$
having the same longitudes as $\calS$ and such that
$\partial\calS=\partial_R\calS$.
Now it is not difficult to realize that one can add some annuli to $\calS'$ in order to obtain
a pair of disjoint generalized Seifert surfaces for $L_\Gamma$ whose interiors do not intersect the isthmus
of $\Gamma$.
\end{proof}

\subsection{Patterns and obstructions}
We are now ready to exploit patterns in order to decide about 
the existence
of $\partial$-connected or $\partial_R$-connected cut systems for $M$.

\begin{prop}\label{jaco:gen}
Suppose $\calS$ is any good cut system for $M$ with associated pattern
$(w_0,w_1,w_2)$. Then:
\begin{enumerate}
 \item 
$M$ admits a $\partial$-connected cut system if and only if 
$(w_0,w_1,w_2)$ is trivial;
\item
if $M$ admits a $\partial_R$-connected cut system, then
there exists a trivial link pattern whose elements are contained in
the subgroup of $F_2$ generated by 
$w_1$ and $w_0w_2w_0^{-1}$.
\end{enumerate}
\end{prop}
\begin{proof}
Point~(1) is an immediate consequence of 
Proposition~\ref{jaco:prop} and Lemma~\ref{generatori}, so we only have
to prove point~(2).

By the very definition of associated pattern, if $\calS'$ is a 
 $\partial_R$-connected cut system for $M$, then, up to suitably choosing
the ordering and the orientations of the longitudes of $\calS'$, we may assume
that $w_1(\calS')=t_1$, $w_2(\calS')=t_2$. By Lemma~\ref{generatori}
(applied with respect to $\calS'$),
this implies that $\varphi(i_\ast(G_\partial))$ contains the link pattern
$(t_1,w_0(\calS')t_2w_0(\calS')^{-1})$, which is obviously trivial.  
But by Lemma~\ref{generatori} again (now applied to the cut system $\calS$), the group $\varphi(i_\ast(G_\partial))$ is generated by
$w_1$ and $w_0w_2w_0^{-1}$, whence the conclusion.
\end{proof}

\subsection{Patterns of $(4)_L$-unknotted handlebodies}\label{4sub}
Let $H$ be  
a $(4)_L$-unknotted
spatial handlebody and let $\Gamma\in\mathcal{S}(H)$ be a spine
of $H$ such that $L_\Gamma=K_1\cup K_2$ is a homology boundary link.
Also denote by $S_1,S_2$ a pair of disjoint generalized Seifert surfaces
for the knots $K_1,K_2$, and by
$\alpha$ the isthmus of $\Gamma$. Up to isotopy, we may suppose that
$\alpha$ transversely intersects $S_1\cup S_2$ in a finite number of points
$x_1,\ldots,x_l$. Moreover, we may order $x_1,\ldots,x_l$ in such a way that
they appear consecutively along $\alpha$ when running from $K_1$ to $K_2$,
and label each $x_i$ with the letter $t_j$ (resp.~$t_j^{-1}$) 
if $x_i$ belongs to $S_j$ and
$\alpha$ intersects $S_j$ at $x_i$ positively (resp.~negatively).
We define the element
$$
\I (S_1,S_2,\Gamma)\in F_2
$$
as the product of the labels of $x_1,\ldots,x_l$.

Up to isotopy, we may assume that $H$ transversely intersects $S_1\cup S_2$
in some annuli (each of which has one boundary component on $K_1\cup K_2$)
and in a collection of meridian disks that separate $H$ and bijectively correspond
to the points of intersection between $\alpha$ and $S_1\cup S_2$.
Let us set $\calS=\{S_1\cap M, S_2\cap M\}$, where as usual 
$M=\compl (H)$, and fix on $S_i\cap M$ the   
orientation induced by $S_i$.

The following Lemma is an immediate consequence of our definitions:

\begin{lem}\label{translation}
We have 
$$
w_0(\calS)=\I(S_1,S_2,\Gamma)\ .
$$
\end{lem}

\begin{lem}\label{translation2}
Let $L$ be a homology boundary link, and suppose that $(w_1,w_2)$ is a
pattern realized by $L$. Then, for every $w_0\in F_2$ there exists 
a spatial handlebody $H$ having $L$ as constituent link
and realizing the pattern $(w_0,w_1,w_2)$.
\end{lem}
\begin{proof}
Let us set $X=\compl (L)$, and denote by $K_1$ and $K_2$ the knots
such that $L=K_1\cup K_2$. By the very definition of associated pattern, 
we may choose disjoint generalized Seifert surfaces $S_1,S_2$ for $K_1,K_2$,
a basepoint $x_0\in X\setminus (S_1\cup S_2)$, and two elements
$m_1,m_2\in  \pi_1(X,x_0)$ representing the meridians of $K_1,K_2$ in such a way that
$\varphi (m_1)=w_1$, $\varphi(m_2)=w_2$, where 
$\varphi\colon \pi_1(X,x_0)\to F_2$ is the epimorphism associated
to $S_1\cup S_2$. 

Let us denote by $\partial_i X$ the component
of $\partial X$ corresponding to $K_i$, and choose a basepoint
$x_i\in \partial_i X\setminus (\partial S_1\cup\partial S_2)$.
For $i=1,2$, let $\overline{m}_i\subseteq \partial_i X$
be a simple loop based at $x_i$ and representing a meridian of $K_i$, 
and choose 
a simple arc $\alpha_i\subseteq X\setminus (S_1\cup S_2)$ 
joining $x_0$ to $x_i$. Let $m'_i\in \pi_1 (X,x_0)$
be the element represented by the loop
$\alpha_i\ast \overline{m}_i\ast \alpha_i^{-1}$. 
Then there exists $h_i\in F_2$ such that
$\varphi (m'_i)=h_i\varphi(m_i)h_i^{-1}=h_iw_ih_i^{-1}$.

Let
now
$\beta$ be the (homotopy class of a) loop in $\pi_1 (X,x_0)$ such that
$\varphi (\beta)=h_1w_0h_2^{-1}$ and let 
$\alpha'=\alpha_1^{-1}\ast\beta\ast\alpha_2$. Let also $\alpha\subseteq X$
be a simple arc satisfying the following properties: it is properly embedded in $X$
with endpoints $x_1,x_2$; it is homotopic in $X$ relative to its endpoints
to the arc $\alpha'$; it transversely intersects
$S_1\cup S_2$ in a finite number of points. Finally, let us set $\Gamma=K_1\cup K_2
\cup \alpha$, $H=N(\Gamma)$ and $M=\compl(H)$. 
Also denote by $\calS=\{S_1\cap M,S_2\cap M\}$ the cut system
of $M$ obtained as above from $S_1\cup S_2$.

Our construction now implies that
$$
w_0(\calS)=\I(S_1,S_2,\Gamma)=h_1w_0h_2^{-1},\quad
w_1(\calS)=\varphi(m_1')=h_1w_1h_1^{-1},\quad
w_2(\calS)=\varphi(m_2')=h_2w_2h_2^{-1}\ .
$$
Lemma~\ref{simpli} now implies that $M$ admits a cut system $\calS'$
such that $w_i(\calS')=w_i$ for $i=0,1,2$, whence the conclusion.
\end{proof}

\subsection{Patterns of $(3)_L$-unknotted handlebodies}
Building on  Proposition~\ref{jaco:gen} we are now able to 
describe an effective 
algorithm which decides if a $(3)_L$-unknotted handlebody 
admits a $(M\to
W)$--boundary--preserving--map (or, equivalently, a $\partial$-connected cut system). 

\begin{prop}\label{criterio2}
Let $H$ be a $(3)_L$-unknotted
spatial handlebody and let $\Gamma$ be a spine
of $H$ such that $L_\Gamma=K_1\cup K_2$ is a boundary link
with Seifert surfaces $S_1,S_2$. We also set as usual $M=\compl(H)$.
Then:
\begin{enumerate}
 \item 
If $\I(S_1,S_2,\Gamma)=t_1^nt_2^m$ for some $n,m\in\mathbb{Z}$, then
$\Gamma$ is a boundary spine for $H$. In particular, $H$ is $(3)_S$-unknotted
and $M$ admits a $\partial$-connected cut system.
\item
Otherwise, $M$ does not admit any $\partial$-connected cut system. In particular,
$H$ is $(3)_S$-knotted.
\end{enumerate}
\end{prop}
\begin{proof}
Let $\calS$ be the cut system of $M$ obtained from $S_1\cup S_2$ as described 
above,
and let $\ell_1,\ell_2$ be the longitudes
of $\calS$, oriented in such a way that $\ell_i\subseteq \partial S_i$  inherits 
the orientation induced by $S_i$.  By Lemma~\ref{translation}, the manifold
$M$ realizes the pattern  $(\I(S_1,S_2,\Gamma), t_1,t_2)$.

Let us now suppose that $M$ admits a $\partial$-connected cut system. 
By Proposition~\ref{jaco:gen} we have that
$(\I(S_1,S_2,\Gamma), t_1,t_2)$ is trivial, and by Remark~\ref{base:rem}
this implies in turn that
$\I(S_1,S_2,\Gamma)=t_1^nt_2^m$ for some
$m,n\in \mathbb{Z}$. We have thus proved point~(2).

In order to conclude, it is sufficient to show that 
if $\I(S_1,S_2,\Gamma)=t_1^nt_2^m$ for some $n,m\in\mathbb{Z}$, then
$\Gamma$ is a boundary spine for $H$.

In fact,
under the assumption $\I(S_1,S_2,\Gamma)=t_1^nt_2^m$ we will show that Lemma~\ref{simpli}
provides an explicit procedure which replaces
$S_1,S_2$ with a pair of Seifert surfaces for $K_1,K_2$ 
whose internal parts are disjoint from the isthmus $\alpha$ of $\Gamma$. 

As mentioned above, we have
$w_0(\calS)=\I(S_1,S_2,\Gamma)=t_1^nt_2^m$, $w_1(\calS)=t_1$, $w_2(\calS)=t_2$. 
We may now apply Lemma~\ref{simpli} to the case
$g_1=t_1^{-n}$, $g_2=t_2^{-m}$, thus obtaining a cut system
$\calS'=\{S'_1,S'_2\}$ such that $\widetilde{w}_1(\calS')=t_1$, $\widetilde{w}_2(\calS')=t_2$
and $\widetilde{w}_0(\calS')$ is the empty word. It follows that
$\calS'$ is $\partial$-connected. Moreover,
since the longitudes $\ell_1,\ell_2$ of $\calS'$ coincide with those of $\calS$,
we may obtain the desired Seifert surfaces for $K_1$ and $K_2$
just by adding to $S_i'$
an annulus $A_i\subseteq H$ bounded by $K_i\cup \ell_i$ for $i=1,2$.
\end{proof}

The following results are immediate consequences of Proposition~\ref{criterio2}. By Theorem~\ref{ex<=>in:teo}
proved in Section~\ref{Ext-vs-intr} below, in Corollary~\ref{casopart} the hypothesis that
$H$ is $(3)_L$-unknotted is superfluous. Note however that the proof
of Theorem~\ref{ex<=>in:teo} relies on very deep results, such as the genus--2 Poincar\'e conjecture
(see Remark~\ref{poincrem}).

\begin{cor}\label{casopart}
Suppose $H$ is a $(3)_L$-unknotted spatial handlebody. Then,
$H$ is $(3)_S$-knotted if and only if $M=\compl (H)$ does not admit any $\partial$-connected cut system.
\end{cor}

\begin{cor}
Suppose $H$ is a $(3)_S$-unknotted spatial handlebody, and let $\Gamma\in \Ss (H)$ be a 
spine of $H$ such that $L_\Gamma$ is a boundary link. Then, $\Gamma$ is a boundary spine of $H$.
\end{cor}

\subsection{$(1)_L$-knotting does not imply $(3)_S$-knotting}\label{newlambert}
Let us come back to the examples described in Subsection~\ref{1Svs1L}
(see Figure~\ref{1-knotted})
and in Subsection~\ref{2Lvs3L:sub} (see Figure~\ref{gamma3:fig}).

The constituent link $L$ of
the spine $\Gamma_1 (p)$ is trivial, and it admits an obvious pair of disjoint
Seifert surfaces $S_1,S_2$ given by the disks bounded by the components of $L$
and lying on the blackboard plane. It is readily seen that
we may choose orientations in such a way that 
$$
\I(S_1,S_2,\Gamma_1(p))=t_1 (t_2t_1)^{\frac{p-1}{2}}t_2\ .
$$ 
In a similar way,
the isthmus of $\Gamma_3 (p)$ intersects the union of the obvious
Seifert surfaces of $L_{\Gamma_3 (p)}$ in such a way that 
$$
\I(S_1,S_2,\Gamma_3(p))=t_1 (t_2t_1)^{\frac{p-1}{2}}t_2\ .
$$ 
Therefore, the criterion described in
Proposition~\ref{criterio2} immediately implies the following:

\begin{prop}\label{alternative}
For every odd prime $p$, the manifolds $M_1 (p)$ and $M_3 (p)$ introduced
in Subsections~\ref{1Svs1L} and~\ref{2Lvs3L:sub} do not admit any $\partial$-connected cut system.
In particular, the handlebodies $H_1 (p)$ and $H_3(p)$ are
$(3)_S$-knotted.  Since $H_1 (p)$ is clearly $(1)_L$-unknotted, it follows that
$(3)_S$-knotting does not imply $(1)_L$-knotting.
\end{prop}

We give a different proof of the fact that $H_1(p)$ is $(3)_S$-knotted in
Proposition~\ref{via-A-obs}.

\subsection{$(3)_L$-knotting does not imply $(4)_L$-knotting}
In this Subsection we exploit Proposition~\ref{jaco:gen}
for constructing examples of $(4)_L$-unknotted handlebodies
which do not admit any constituent boundary link.
More precisely,
we prove the following:

\begin{prop}\label{nuovaprop}
Let $L$ be a homology boundary link which is not a boundary link. Then,
$L$ is a constituent link of a handlebody $H$ whose complement $M=\compl (H)$
does not admit any $\partial_R$-connected cut system.
In particular:
\begin{enumerate}
\item
 $H$ is $(3)_L$-knotted and $(4)_L$-unknotted. 
\item
$M$ admits a cut system, but does not admit any
$\partial_R$-connected cut system.
\end{enumerate}
\end{prop}
\begin{proof}
%Let us set $X=\compl (L)$,  
%let $\calS=\{S_1,S_2\}$ be a cut system for $X$
%and
%fix a basepoint $x_0\in X\setminus \calS$. If $K_1$ and $K_2$ are the knots
%such that $L=K_1\cup K_2$, we denote by $\partial_i X$ the component
%of $\partial X$ corresponding to $K_i$, and we fix
%simple arcs $\alpha_1, \alpha_2$ in $X$ satisfying the following properties:
%for $i=1,2$, the arc $\alpha_i$ is disjoint from $S_1\cup S_2$ and joins 
%$x_0$ to $\partial_i X$;
%the intersection $\alpha_1\cap \alpha_2$ reduces to the basepoint $x_0$. 
%We also denote by $x_i$ the endpoint of $\alpha_i$ lying on $\partial_i X$.

%Let $\varphi\colon \pi_1 (X,x_0)\to F(t_1,t_2)=F_2$ the epimorphism associated
%to $\calS$.
%For $i=1,2$, 
%let $m_i\subseteq \partial_i X$ be a meridian loop based at $x_i$,
%and set $w_i=\varphi(\alpha_i\ast m_i\ast \alpha_i^{-1})$. The pair
%$(w_1,w_2)$ is a link pattern realized by $L$. 

Let $(w_1,w_2)$ be a pattern realized by $L$.
Recall that $w_1$ and $w_2$ project onto a basis of $F_2/[F_2,F_2]$, so 
if we had $w_1=t_1^k$ and $w_2=t_2^h$ for some $h,k\in\mathbb{Z}$, then
we
would get $k=\pm 1$, $h=\pm 1$, and $(w_1,w_2)$ would be trivial as a link pattern,
against our assumption that $L$ is not a boundary link.
Let us assume that $w_1\neq t_1^{\pm h}$, the case when
$w_2\neq t_2^{\pm h}$ being similar.
%Therefore, we may assume for example that the reduced word representing $w_1$
%contains both the letter $t_1^{\pm 1}$ and the letter $t_2^{\pm 1}$
%(the case when $w_2$ contains both letters is analogous).
Then there exists $n\in\mathbb{N}$ such that the following conditions hold:
\begin{enumerate}
 \item 
the reduced word representing $t_1^{n} w_1 t_1^{-n}$ is given by
$t_1^{k_1} z_1 t_1^{-h_1}$, where $h_1\geq 1$, $k_1\geq 1$ (so
$z_1$ is not empty, and contains the symbol $t_2^{\pm 1}$);
\item
the reduced word representing $t_2^n w_2 t_2^{-n}$ is either equal
to $t_2^{l_2}$ for some $l_2\in\mathbb{Z}$, or to 
$t_2^{k_2} z_2 t_2^{-h_2}$, where $h_2\geq 1$, $k_2\geq 1$ (and in this case
$z_2$ is not empty, and contains the symbol $t_1^{\pm 1}$).
\end{enumerate}

By Lemma~\ref{translation2}, there exists a handlebody $H$ with constituent link $L$
and associated pattern $(t_1^{-n}t_2^{n}, w_1,w_2)$. 
Let us set $w'_1=t_1^nw_1t_1^{-n}$, $w'_2=t_2^nw_2t_2^{-n}$. 
By Lemma~\ref{simpli},
the handlebody $H$ also realizes the pattern
$(1,w_1', w_2')$.
Let us set $M=\compl (H)$.
By Proposition~\ref{jaco:gen}--(2),
in order to conclude it is sufficient to show that there does not exist a trivial link
pattern whose elements are contained in the subgroup $J$ of $F_2$ generated by
$w_1'$ and $w_2'$. 
%The conjugation by $t_1^{n}$ defines
%an automorphism of $F_2$ which takes this subgroup into the subgroup $J$
%generated by 
%$$
%w_1'=t_1^n w_1t_1^{-n},\qquad w_2'=t_2^n w_2t_2^{-n}\ ,
%$$
%so it is sufficient to show that $J$ does not contain two elements which define
%a trivial link pattern.

We first look at which elements of $J$ can be primitive in $F(t_1,t_2)$.
So, let us take a non--trivial element $R\in F(a,b)$ such that $R(w'_1,w'_2)$
is primitive in $F(t_1,t_2)$. Let $R'$ be a cyclically reduced conjugate
of $R$, and observe that $R'(w_1',w_2')$ is also primitive.
Let us first assume that both $a$ and $b$ appear in $R'$ (\emph{i.e.}~that
$R'$ is not of the form $R'=a^h$ or $R'=b^h$ for some $h\in\mathbb{Z}$).
Then, points~(1) and~(2) above easily imply that the reduced word 
representing $R'(w_1',w_2')$ is ciclically reduced, and contains the symbol $t_1$
both with positive and with negative exponents. The criterion described in
Lemma~\ref{free1:lemma}--(2) now implies that $R'(w_1',w_2')$ cannot be primitive
in $F_2$.

Therefore, if $R(w_1',w_2')$ is primitive, then
we have either $R'(w_1',w_2')=t_1^{n}w_1^h t_1^{-n}$ or
$R'(w_1',w_2')=t_2^{n}w_2^h t_2^{-n}$ for some $h\in\mathbb{Z}$.
Since primitives of $F_2$ project onto indivisible elements
of $F_2/[F_2,F_2]$, this forces $h=\pm 1$. Since $R(w_1',w_2')$ is conjugate
(in $J$) to $R'(w_1',w_2')$, we may conclude that every element
of $J$ which is primitive in $F_2$ is conjugated (in $F_2$) either to $w_1^{\pm 1}$
or to $w_2^{\pm 1}$. Since the elements of a pattern project onto a base
of $F_2/[F_2,F_2]$, this implies that if $(z_1,z_2)$ is a link pattern 
whose elements are contained in $J$, then $(z_1,z_2)$ is equivalent
to $(w_1^{\pm 1},w_2^{\pm 1})$. Since $L$ is not a boundary link, the pattern
$(w_1^{\pm 1},w_2^{\pm 1})$ is not trivial, and this implies in turn that $(z_1,z_2)$
is not trivial. We have eventually showed that there does not exist any trivial link pattern
whose elements are contained in $J$, and this concludes the proof of the Proposition.
\end{proof}

\subsection{An example}
We now describe an explicit example of a $(3)_L$-knotted handlebody admitting
a constituent homology boundary link.

\begin{figure}
 \begin{center}
  \input{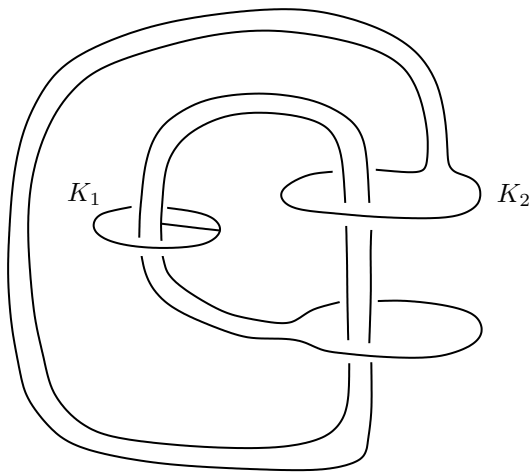}
\caption{The spine of a $(4)_L$-unknotted handlebody which is $(3)_L$-knotted.}\label{hom1:fig}
 \end{center}
\end{figure}

\begin{figure}
 \begin{center}
  \input{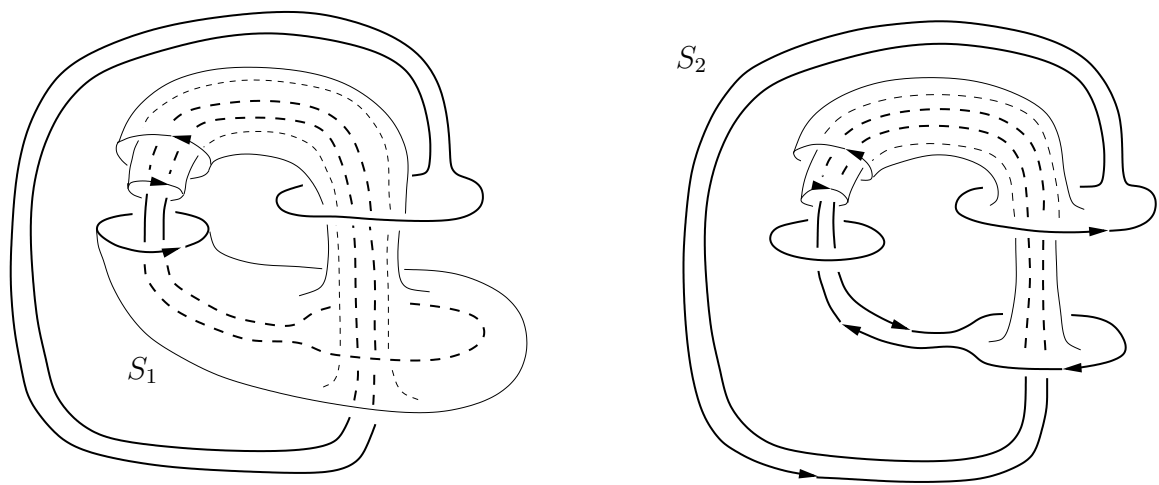}
\caption{The generalized Seifert surfaces $S_1$ and $S_2$. In order to get a clearer
picture, for $i=1,2$ we have cut from $S_i$ small annular neighbourhoods 
of two components of $\partial S_i$.}\label{hom2:fig}
 \end{center}
\end{figure}

Let $\Gamma$ be the spine described in Figure~\ref{hom1:fig}, and set
$H=H(\Gamma)$, $M=\compl(H)$. The constituent link $L_\Gamma=K_1\cup K_2$
was introduced by Cochran and Orr in~\cite{CO2,CO1} and provides the first
example of a homology boundary link which is not concordant to a boundary link
(in particular, it is not a boundary link).
Following~\cite[Diagram 1.4]{tesi},
in Figure~\ref{hom2:fig} we have drawn a pair $S_1,S_2$ of generalized Seifert surfaces for
$K_1\cup K_2$. The three boundary components of $S_1$ all lie on $K_1$ (and exactly two of them have the same orientation as $K_1$), while $S_2$ has one boundary component on $K_2$
and two boundary components (having opposite orientations) on $K_1$.

Let us denote by $\calS=\{S_1\cap M,S_2\cap M\}$ the cut system of $M$ defined by $S_1$ and
$S_2$ (see Subsection~\ref{4sub}). Up to isotopy, we may suppose that
the isthmus of $\Gamma$ is disjoint from ${\rm int}(S_1)\cup {\rm int}(S_2)$. It is now an easy exercise to show that
$$
w_0(\calS)=\I(S_1,S_2,\Gamma)=1,\quad
w_1(\calS)=t_1t_2t_1^{-1}t_2^{-1}t_1,\quad
w_2(\calS)=t_2\ .
$$
Since in each product of the form $w_1(\calS)^{\pm 1}w_2(\calS)^{\pm 1}$
there cannot be cancellations, just as in the proof of Proposition~\ref{nuovaprop}
we can conclude that $M$ does not admit any $\partial_R$-connected cut system.
In particular, the handelbody $H$ is $(3)_L$-knotted, while being
obviously $(4)_L$-unknotted.
%qui

\subsection{The maximal free covering}\label{maximalfree:sub}
Let $H$ be a spatial handlebody, set as usual $M=\compl(H)$
and suppose that $\cut(M)=2$.
Under this assumption, we are now going to describe
the \emph{maximal free covering} $\widetilde{M}_\omega$
of $M$. As mentioned at the beginning of the Section, 
we will discuss some aspects of the topology of $\widetilde{M}_\omega$
which are related to the knotting level of $H$.

Let $\Ss=\{S_1,S_2\}$ be any cut system of $M$, and let 
$V$ be the manifold with
boundary obtained by cutting $M$ along $\calS$.  
Then, the boundary of
$V$ consists of some ``horizontal'' boundary region (given by
$\partial M\cap \partial V$) and some ``vertical'' boundary region,
coming from the cuts along $S_1$ and $S_2$.  More precisely, recall that $S_i$
is oriented, and call $S_i^+$, $S_i^-$ the vertical
components of $\partial V$ associated to $S_i$, $i=1,2$, 
in such a way that at any point of $S_i$ a positive basis of the tangent space of $S_i$ is completed to a positive
basis of the tangent space of $M$ by adding a vector pointing towards $S_i^+$.

Let us set $G=\pi_1(M,x_0)$, where $x_0$ is a basepoint such that
$x_0\in\partial M\setminus ( S_1\cup S_2)$, 
let $F_2=F(t_1,t_2)$ be the free group on two generators,
and denote by $\varphi\colon G\to
F_2$ the epimorphism associated to $\calS$.

Let $p_\omega\colon \widetilde M_\omega\to M$ be the covering
associated to $\ker \varphi$. The space $\widetilde M_\omega$ admits
an easy topological description as a tree of spaces (called \emph{pieces}), where each piece is homeomorphic
to $V$. 
If $\{V_h\}_{h\in F_2}$ is a
countable family of copies of $V$ indexed by the elements of $F_2$,
then $\widetilde M_\omega$ is homeomorphic to the quotient of the
disjoint union $\bigsqcup_{h\in F_2} V_h$ by the equivalence relation generated by
 $$
V_h\ni x \ \sim \ y \in V_{h'} \quad \Longleftrightarrow \quad
h'=h_i h,\ x\in S_i^-\subseteq V_h,\, y\in S_i^+\subseteq V_{h'},\
{\rm and}\ x=y\ {\rm in}\ M\ ,
$$
where we  now consider $M$ as the space obtained from any $V_h$ by gluing
in pairs its vertical boundary components.

The group of the covering automorphisms of $\widetilde{M}_\omega$ is isomorphic to $F_2$, and
for every $h_0, h\in F_2$ the covering translation associated to $h_0$
translates $V_h$ onto its copy $V_{h_0 h}$.

The notation $\widetilde M_\omega$ is justified by 
Theorem~\ref{stallings},
which states that $\ker\varphi=G_\omega$, and implies
that the topology of $\widetilde M_\omega$ does not
depend on the particular epimorphism $\varphi$, nor on the chosen
cut system $\calS$, and is therefore intrinsically associated to $M$
(whence to $H$). The covering
  $p_\omega\colon \widetilde
M_\omega\to M$ is called the {\it maximal free
  covering} of $M$ because
the group of the covering automorphisms of $\widetilde{M}_\omega$ is 
isomorphic to the maximal free quotient of $G$.

\smallskip

The proof of the following Lemma is elementary, and it is left to the reader.

\begin{lem}\label{connect:lemma}
  Suppose $X\subseteq M$ is path--connected, choose a base point
  $x_0\in X\subseteq M$ and denote by $i\colon X\to M$ the
  inclusion. Let $p\colon (\widehat{M},\widehat{x}_0)\to (M,x_0)$ be a
  regular covering, and set $\widehat{X}=p^{-1} (X)$. Then
  $\widehat{X}$ is path--connected if and only if $p_\ast (\pi_1
  (\widehat{M}, \widehat{x_0}))\cdot i_\ast(\pi_1 (X,x_0))$ coincides
  with the whole group $\pi_1 (M,x_0)$.
\end{lem}
%\begin{proof}
 % The set $\widehat X$ is path--connected if and only if every point
  %  $\widehat{x}_1\in p^{-1} (x_0)$ can be connected to $\widehat{x}_0$
 % by a path in $\widehat{X}$.

  %Suppose first that $p_\ast (\pi_1 (\widehat{M}, \widehat{x}_0))\cdot
  %i_\ast(\pi_1 (X,x_0))=\pi_1 (M,x_0)$, and take a path
  %$\widehat{\gamma}$ joining $\widehat{x}_0$ to $\widehat{x}_1$ in
  %$\widehat M$.  Our hypothesis implies that the path
  %$p\circ\widehat\gamma$ is homotopic to the concatenation of loops
  %$\alpha$ and $\beta$ based at $x_0$, where the homotopy class of
  %$\alpha$ lies in $p_\ast (\pi_1 (\widehat M,\widehat{x}_0))$, and
  %the image of $\beta$ is contained in $X$.  Now, the path
  %$\alpha\ast\beta$ lifts to a path starting at $\widehat{x}_0$ and
  %ending at $\widehat{x}_1$ which is the concatenation of a loop based
  %at $\widehat{x}_0$ and a path contained in $\widehat X$.  It follows
  %that $\widehat X$ is path--connected.

%  Suppose now that $\widehat X$ is path--connected, and let $\gamma$
 % be a loop in $M$ based at $x_0$. If $\widehat{\gamma}$ is the lift
  %of $\gamma$ starting at $\widehat{x}_0$, then there exists a path
  %$\widehat{\beta}$ starting at $\widehat{x}_0$ and contained in
  %$\widehat X$ such that the concatenation
  %$\widehat{\alpha}=\widehat{\gamma}\ast\widehat{\beta}^{-1}$ is
  %defined. Then, $\gamma$ is homotopic to
  %$(p\circ\widehat{\alpha})\ast (p\circ\widehat{\beta})$, whence the
  %conclusion.
%\end{proof}
%\cvd

We now apply the previous Lemma to the case  we are interested in.
The following result provides an interesting relation between the
topology of $\widetilde{M}_\omega$ and the  knotting level of $H$.

\begin{prop}\label{bordorivestimento}
%(1) The subspace $\partial \widetilde{M}$ of $\widetilde M$
%is path connected.
The subspace $\partial\widetilde{M}_\omega$ of
$\widetilde{M}_\omega$ is path connected if and only if $M$ admits a
$(M\to W)$--boundary--preserving--map.
\end{prop}
\begin{proof}
%We have the commutative diagram
%$$\xymatrix {
%\pi_1 (\partial M,x_0)\ar[r]^{i_\ast} \ar[d]^\theta & 
%\pi_1 (M,x_0)\ar[d]^\theta\\
%H_1 (\partial M) \ar[r]^{i_\ast} & H_1(M) },
%$$
%where the vertical arrows are given by the Hurewicz epimorphisms.
%Since $i_\ast\colon H_1(\partial M)\to H_1 (M)$ is surjective and the kernel
%of the homomorphism $\theta\colon \pi_1 (M,x_0)\to H_1 (M)$ is exactly 
%the commutator subgroup of $\pi_1 (M,x_0)$ we have
%$$
%[\pi_1 (M,x_0),\pi_1 (M,x_0)]\cdot i_\ast (\pi_1 (\partial M,x_0))=\pi_1 (M,x_0),
%$$
%so Lemma~\ref{connect:lemma} implies that $\partial\widetilde{M}$ is
%connected.
%
%\medskip
%
There is an elementary topological proof of the fact that if $M$
admits a $\partial$-connected cut system $\calS$, then
$\partial\widetilde{M}_\omega$ is connected.  In fact, 
in this case
the horizontal boundary of
$V=M\setminus \calS$ is connected. Therefore, the boundary of
$\widetilde{M}_\omega$ is obtained by gluing connected spaces
following a tree--like pattern, and is therefore connected.  

To get
both implications it is sufficient to observe that Lemma~\ref{connect:lemma}
implies that $\partial \widetilde{M}_\omega$ is connected if and only if
$G=i_\ast(G_\partial) \cdot G_\omega$, and Proposition~\ref{jaco:prop} 
ensures that this last condition is equivalent
to the fact that $M$ admits a
$(M\to W)$--boundary--preserving--map.
\end{proof}

\subsection{Lifts of longitudes}\label{longitude2}
If $X= \compl (L)$ is the complement of a 2--component homology
boundary link and $\calS$ is a cut system of $X$, we can construct in the same way as above 
the maximal free covering $\widetilde{X}_\omega$ of $X$.
Of course, since $\partial X$ is disconnected, the space
$\partial \widetilde{X}_\omega$ cannot be
connected. Also observe that every connected component 
of $\partial\widetilde X_\omega$ is obtained by gluing to each other an infinite number of annuli, and is therefore homeomorphic to
an annulus.

Recall that
a \emph{longitude} of $L$ is (the isotopy class in $\partial X$ of) 
a connected component
of $\partial \Ss$. 
Every 
longitude lifts to a loop in $\partial\widetilde
X_\omega$ which generates the first homology 
group of the annular component of $\partial\widetilde
X_\omega$ where it lies. Since two non--trivial simple loops on an annulus are isotopic, 
this readily implies the already mentioned fact that
longitudes do not depend on the fixed cut system $\Ss$. Moreover,
since $S_1$ and $S_2$ also lift to $\widetilde{X}_\omega$, if $\partial S_1$ and $\partial S_2$
are connected then every lift of a longitude bounds in $\widetilde{X}_\omega$, and is 
therefore null--homologous in $\widetilde{X}_\omega$. 
We can summarize this brief discussion in the following well-known:

\begin{lem}\label{longitudini}
The $\Z$--module $H_1 (\partial \widetilde{X}_\omega)$ is generated
by the lifts of the longitudes of $\partial X$. If $L$ is a boundary
link, then  $\widetilde{i}_*(H_1 (\partial \widetilde{X}_\omega))=0$, where 
$\widetilde{i}_*$ is induced
by the inclusion $\widetilde{i}\colon \partial\widetilde{X}_\omega\to \widetilde{X}_\omega$.  
\end{lem}

In Propositions~\ref{comb2} and~\ref{comb3} below we
extend Lemma~\ref{longitudini}
to the case when $M=\compl (H)$ is the complement
of a handlebody $H$ such that $\cut (M)=2$, thus obtaining
some more obstructions for $M$ to admit $\partial$--connected or $\partial_R$--connected
boundary. Notice that
some difficulties arise due to the fact that 
in this case the definition of longitudes is less obvious, and longitudes are in fact
no more independent of the choice of a cut system.

\smallskip 

If $\calS$ is a good cut system for $M$, let us define $L(\calS)\subseteq
H_1 (\partial\widetilde{M}_\omega)$ as the submodule generated by the
lifts to $\partial\widetilde{M}_\omega$ of the longitudes of $\calS$.
We also denote by $\varphi\colon \pi_1(M,x_0)\to F_2$ the epimorphism induced by $\calS$,
where $x_0$ is a basepoint in $\partial M\setminus (S_1\cup S_2)$.

\begin{prop}\label{comb2}
We have $L(\calS)=H_1 (\partial \widetilde{M}_\omega)$.
In particular, $L(\calS)$ does not depend on $\calS$.
\end{prop}
\begin{proof}
Take a basepoint $\widetilde{x}_0\in p_\omega^{-1} (x_0)\subseteq \partial \widetilde{M}_\omega$, and let $\partial_0 \widetilde{M}_\omega$ be the
connected component of $\partial\widetilde{M}_\omega$ containing
$\widetilde{x}_0$. Since the maximal free covering is regular,
it is sufficient to show that $H_1(\partial_0 \widetilde{M}_\omega)$
is generated by those lifts of longitudes of $\calS$
that lie on $\partial_0 \widetilde{M}_\omega$.

So, let us take $z\in H_1(\partial_0 \widetilde{M}_\omega)$. By Hurewicz's
Theorem we may suppose that $z$ is represented by (the class of) a loop
$\widetilde{\gamma}\in \pi_1 (\partial_0 \widetilde{M}_\omega,\widetilde{x}_0)$.
Let now $i\colon \partial M\to M$, $\widetilde{i}_0\colon \partial_0 \widetilde{M}_\omega
\to \widetilde{M}_\omega$ be the inclusions, let us denote
the restriction of $p_\omega$ to $\partial_0 \widetilde{M}_\omega$
simply by $p_\omega$, and set $\gamma=(p_\omega)_\ast (\widetilde{\gamma})\in
\pi_1(\partial M,x_0)$. Since $p_\omega\circ\widetilde{i}_0=i\circ p_\omega$ and
$(p_\omega)_\ast (\pi_1 (\widetilde{M}_\omega,\widetilde{x}_0))=\ker\varphi$, we have
$\gamma\in \ker\varphi\circ i_\ast$. 

Let now $\gamma_1,\gamma_2$ be elements of $\pi_1 (\partial M,x_0)$ 
which are represented by simple loops isotopic to the longitudes of $\calS$,
and recall from Lemma~\ref{generatori}--(3) that $\ker \varphi\circ i_\ast$
coincides with 
the smallest normal subgroup  of $\pi_1 (\partial M,x_0)$ containing $\gamma_1,\gamma_2$.
Therefore, $\gamma$ is a product of conjugates of $\gamma_1,\gamma_2$ in $\pi_1(\partial M,x_0)$, so $\widetilde{\gamma}$ is a product of loops each of which is homologous
to the lift of a longitude. This implies that $z\in L(\calS)$, whence the conclusion. 
\end{proof}

\subsection{The image of $H_1(\partial\widetilde{M}_\omega)$ in $H_1(\widetilde{M}_\omega)$
as an obstruction}
Recall that the group of the covering automorphisms of $\widetilde{M}_\omega$
is isomorphic to the free group $F_2$. Therefore, both $H_1 (\partial \widetilde{M}_\omega)$
and $H_1(\widetilde{M}_\omega)$ admit a natural structure of $\mathbb{Z} F_2$--module,
where $\mathbb{Z} F_2$ is the group ring of $F_2$. Moreover, if
$\widetilde{i}\colon \partial\widetilde{M}_\omega\to \widetilde{M}_\omega$ is the inclusion, then
$\widetilde{i}_\ast\colon  H_1(\partial\widetilde{M}_\omega)\to H_1(\widetilde{M}_\omega)$
is a homomorphism of $\mathbb{Z} F_2$--modules, so 
$\widetilde{i}_\ast (H_1(\partial\widetilde{M}_\omega))$ is a  
$\mathbb{Z} F_2$--submodule of $H_1(\widetilde{M}_\omega)$.

\smallskip

We are now ready to point out a further topological obstruction to the existence
of $\partial$-connected and $\partial_R$-connected cut systems for $M$.

\begin{prop}\label{comb3}
The following facts hold:
\begin{enumerate}
\item
 If $M$ admits a $\partial$-connected cut system, then
$\widetilde{i}_\ast (H_1(\partial\widetilde{M}_\omega))=\{0\}$.
\item 
If $M$ admits a $\partial_R$-connected cut system, then
$\widetilde{i}_\ast (H_1(\partial\widetilde{M}_\omega))$ is a cyclic $\mathbb{Z} F_2$--module.
\end{enumerate}
\end{prop}
\begin{proof}
(1):  By Proposition~\ref{comb2}, $H_1 (\partial \widetilde{M}_\omega)$ is
  generated by the lifts of the longitudes.  But, since $\partial
  S_i$ is connected for $i=1,2$, the longitudes bound in
  $\widetilde{M}_\omega$, whence the conclusion.
  
  \smallskip
  
  (2): Let $\calS$ be a $\partial_R$-connected cut system for $M$ with longitudes $\ell_1,\ell_2$. 
  By point~(1), we may assume that $\partial\calS\setminus \partial_R\calS$ is a non--empty collection
  of simple loops. Let $\ell_3$ be one of these loops, and observe that $\ell_3$ lifts to a loop
  $\widetilde{\ell}_3$ in $\partial\widetilde{M}_\omega$. For $i=1,2$, since $\calS$ is $\partial_R$--connected and
  $S_i$ lifts to $\widetilde{M}_\omega$, every lift of $\ell_i$ is homologous in $\widetilde{M}_\omega$
  to a sum of parallel copies of translates of $\widetilde{\ell}_3$ (such a sum is empty if $\partial S_i$ is connected).
  This implies that every lift of a longitude of $\calS$ lies in the cyclic $\mathbb{Z} F_2$--submodule of $H_1 (\widetilde{M}_\omega)$
  generated by $\widetilde{\ell}_3$. The conclusion now follows from Proposition~\ref{comb2}.
\end{proof}

\begin{cor}\label{3kn-4kn}
Let $H$ be a spatial handlebody, set $M=\compl(H)$ and suppose that
$\cut (M)=2$. Then:
\begin{enumerate}
\item
If $H$ is $(3)_S$-unknotted, then $\widetilde{i}_\ast (H_1(\partial\widetilde{M}_\omega))=\{0\}$.
\item
If $H$ is $(3)_L$-unknotted, then $\widetilde{i}_\ast (H_1(\partial\widetilde{M}_\omega))$ is cyclic (as a $\mathbb{Z}F_2$--module).
\end{enumerate}
\end{cor}

\begin{remark}
{\rm Alexander--type obstructions, which will be described in
Section~\ref{Alex-obs}, arise from the analysis of the \emph{maximal abelian covering}
$\widetilde M$ of $M$. One may wonder if the arguments developed in this Section could 
take place in that (more classical) context, but this does not seem the case.} 

{\rm For example, an easy application of Lemma~\ref{connect:lemma} implies that $\partial\widetilde{M}$ is always connected (even if
$\cut (M)=1$), so that the maximal abelian covering cannot provide obstructions as the one described in Proposition~\ref{bordorivestimento}.}

{\rm  Moreover, while the last results of this Section are inspired by analogous results for links, it turns out that the maximal abelian covering of
a handlebody complement (having maximal cut number) displays properties quite dissimilar from the ones of maximal abelian coverings
of (homology boundary) links.
For instance, while Lemma~\ref{longitudini} (which concerns links) also holds when the maximal
free covering is replaced by the maximal abelian one, even when $H$ is unknotted the image
of $H_1(\partial \widetilde{M})$ into $H_1(\widetilde{M})$ does not vanish, so  Proposition~\ref{comb3}
and Corollary~\ref{3kn-4kn} do not admit analogous statements if $\widetilde{M}_\omega$ is replaced by $\widetilde M$.

In the case of (the complement of) homology boundary links, the structure of the first homology group
of the maximal free covering as a $\mathbb{Z}F_2$--module is studied in detail in~\cite{hibook}.}
\end{remark}

\section{Extrinsic vs intrinsic levels of knotting}\label{Ext-vs-intr}
In this Section we investigate the implications that 
the existence of $\partial$-connected or $\partial_R$-connected
cut systems for $M=\compl (H)$ has on the knotting level of $H$.
Since the existence of such cut systems is clearly an intrinsic
property of $M$, this issue mainly concerns the relations between
the intrinsic and the extrinsic properties of handlebody complements.

More precisely, recall that if $H$ is $(3)_S$-unknotted (resp.~$(3)_L$-unknotted),
then $M$ admits a $\partial$-connected (resp.~$\partial_R$-connected) cut system.
We now face the question whether also the converse implications are true.

\medskip

\begin{defi}\label{H-separate}{ \rm Let $M=\compl(H)$ as usual.
A cut system  $\Ss=\{S_1,S_2\}$
of $M$ is said {\it $H$-separated} if there exists a simple
 essential curve $\ell$ on $\partial M$, such that:
\begin{enumerate}
\item $\ell$ separates $\partial M$;

\item $\ell$ does not intersect $S_1\cup S_2$;

\item $\ell$ bounds a compressing 2--disk in $H$. 
\end{enumerate}
 }
\end{defi}

It is not difficult to show that the map that associates to any spine $\Gamma$ of $H$ the boundary of the meridian disk dual to
the isthmus of $\Gamma$ establishes a bijection between
the isotopy classes (in $H$) of (hc)--spines of $H$ and
the isotopy classes (in $\partial H$) of the loops $\ell$ satisfying properties~(1) and (3) described in Definition~\ref{H-separate}
above. Building on this remark, it is not difficult to prove the following:
 
\begin{lem} Let $M=\compl (H)$ be as usual. Then:
\begin{enumerate}
\item[(a)] $H$ is $(3)_S$-unknotted if and only if $M$ admits a $H$-separated
$\partial$-connected cut system.
\item[(b)]  $H$ is $(3)_L$-unknotted if and only $M$ admits a $H$-separated
$\partial_R$-connected cut system.
\item[(c)]  $H$ is $(4)_L$-unknotted if and only if $M$ admits a $H$-separated
cut system.
\end{enumerate}
\end{lem}

On the other hand we can prove the following remarkable equivalence
between extrinsic and intrinsic properties, which was stated
as Theorem~\ref{ex<=>in:teo} in Section~\ref{cut}. 

\begin{teo}\label{ex<=>in} $M=\compl(H)$ admits a 
  $\partial$-connected cut system if and only if $H$ is not
  $(3)_S$-knotted (equivalently, $M$ admits a $H$-separated
  $\partial$-connected cut system).
\end{teo}
\begin{proof} By using a $(M\to W)$--boundary--preserving--map $f: M\to W$ we
can construct a degree--1 map $g: S^3 \to N:= H\cup_fW$. It was
remarked in \cite{Bing} that such a 3--manifold $N$ is a homotopy
sphere. As the Poincar\'e conjecture holds true, then $N$ is
homeomorphic to $S^3$ and is endowed by construction with a Heegaard
splitting (of genus 2). Since every Heegaard splitting of the sphere is trivial,
$W$ admits a $H$-separated
$\partial$-connected cut system (actually made by two 2--disks), say
$\Ss$. We can put $g$ transverse to $\Ss$, without modifying it on a
neighbourhood of $H$. Then the pull-back of $\Ss$ via $g$ provides
the required $H$-separated cut system of $M$.
\end{proof}

\begin{remark}\label{poincrem} {\rm 
By~\cite{birhil}, every 3-manifold with a Heegaard splitting of
genus two is a two-sheeted cyclic branched cover of $S^3$ branched over a knot or link. 
This reconduces the validity of the Poincar\'e conjecture for manifolds
of genus at most two to
the positive solution to the Smith conjecture~\cite{smith}.
However,
    a statement similar to Theorem~\ref{ex<=>in} holds for $H$ of arbitrary genus (see below),
    and in such a generality it is very close to be equivalent to the
    full Poincar\'e conjecture.  
This remark strongly suggests
that the issue of characterizing the relations
    between intrinsic and extrinsic properties of spatial handlebodies
    definitely involves deep results in $3$--dimensional topology.}
\end{remark}

\begin{remark} \label{nocommon}{\rm In general a given $\partial$-connected cut system
    is {\it not} $H$-separated. For example let $H$ be unknotted. Then also
    $M=\compl(H)$ is a handlebody. Let us take a (hc)--spine $\Gamma$ of $M$ such
    that at least one component of $L_\Gamma$ is non--trivial (see \emph{e.g.}~Figure~\ref{tunnel}). We claim that
        the compression disks $D_1,D_2$ of $M$ dual to the components of $L_\Gamma$
    form a $\partial$-connected cut system of $M$ which is not
    $H$-separated. In fact, suppose that $D_3'\subseteq H$ is a separating meridian disk whose boundary
    is disjoint from $\partial D_1\cup \partial D_2$. Now, cutting $M$ along $D_1\cup D_2$ we obtain a ball,
    so $\partial D_3'$ bounds a meridian disk  
    $D_3\subseteq M$ separating $M$ and disjoint from $D_1\cup D_2$. The $2$--sphere $D_3\cup D_3'$ is a reducing sphere
    for the Heegaard splitting $S^3=M\cup H$, and the components of $L_\Gamma$ appear now as cores of tori of a genus--1 Heegaard
    splittings of $S^3$. This implies that both the components of $L_\Gamma$ are unknotted, a contradiction.

Also observe that if both components of $L_\Gamma$ are unknotted, then for $i=1,2$ the loop $\partial D_i\subseteq \partial M=\partial H$
cannot transversely intersect in exactly one point the boundary $\partial D'_i$ of a meridian disk $D'_i\subseteq H$ (otherwise,
by compressing $M$ along $D_i$ one would get a genus--1 Heegaard splitting of $S^3$, so the component of $L_\Gamma$ dual
to $D_j$, $j\neq i$, would be unknotted). This readily implies that the standard pair of longitudes of $M$ (\emph{i.e.}~the pair associated
to the compressing disks dual to the constituent link of an unknotted spine for $M$)
and the pair of longitudes associated to the cut system $\{D_1,D_2\}$ have no element in common.
    }
\end{remark}

\section{Alexander module obstructions}\label{Alex-obs}
In this Section we will recognize obstructions having a much more
classical flavour as they are based on the {\it elementary
  determinantal ideal} $E_2(G)$ derived from any presentation of the
{\it Alexander module} $A(G)$ of the fundamental group $G$ of $M=
\compl(H)$. Note that every invariant arising in this way is forced
to detect only {\it intrinsic}
features of $M$, \emph{i.e.}~only properties that do not depend
on the realization of $M$ as
a cube--with--holes. 
The main result 
proved in 
this Section is
%Proposition~\ref{K3} below and 
Theorem~\ref{nonempty}, which states that 
there exists an infinite family 
of handlebodies $\{H_i\}_{i\in I}$ such that every $M_i=\compl (H_i)$ has
  cut number equal to 1 (hence $H_i$ is $(4)_L$-knotted), and
$M_i$ is not homeomorphic to $M_j$ for $i\neq j$.
However, this Section is mainly devoted to discuss in detail how Alexander invariants can provide obstructions
to the existence both of generic and of $\partial$--connected cut systems: such obstructions are described
 in Proposition~\ref{main-ob}, and applied in Propositions~\ref{via-A-obs}
and~\ref{via-A-obs2}. More precisely, in Proposition~\ref{via-A-obs} we give a different proof
of Proposition~\ref{alternative}, which asserts that
the handlebodies $H_3(p)$ introduced in  Subsection~\ref{1Svs1L} are $(3)_S$-knotted,
while in Proposition~\ref{via-A-obs2} we prove that the complement of Kinoshita's graph (see Figure~\ref{kino} below)
has cut number equal to one.

\subsection{A short account about the existing literature}
\label{ex-lit}
\begin{figure}[htbp]
\begin{center}
 \includegraphics[height=5cm]{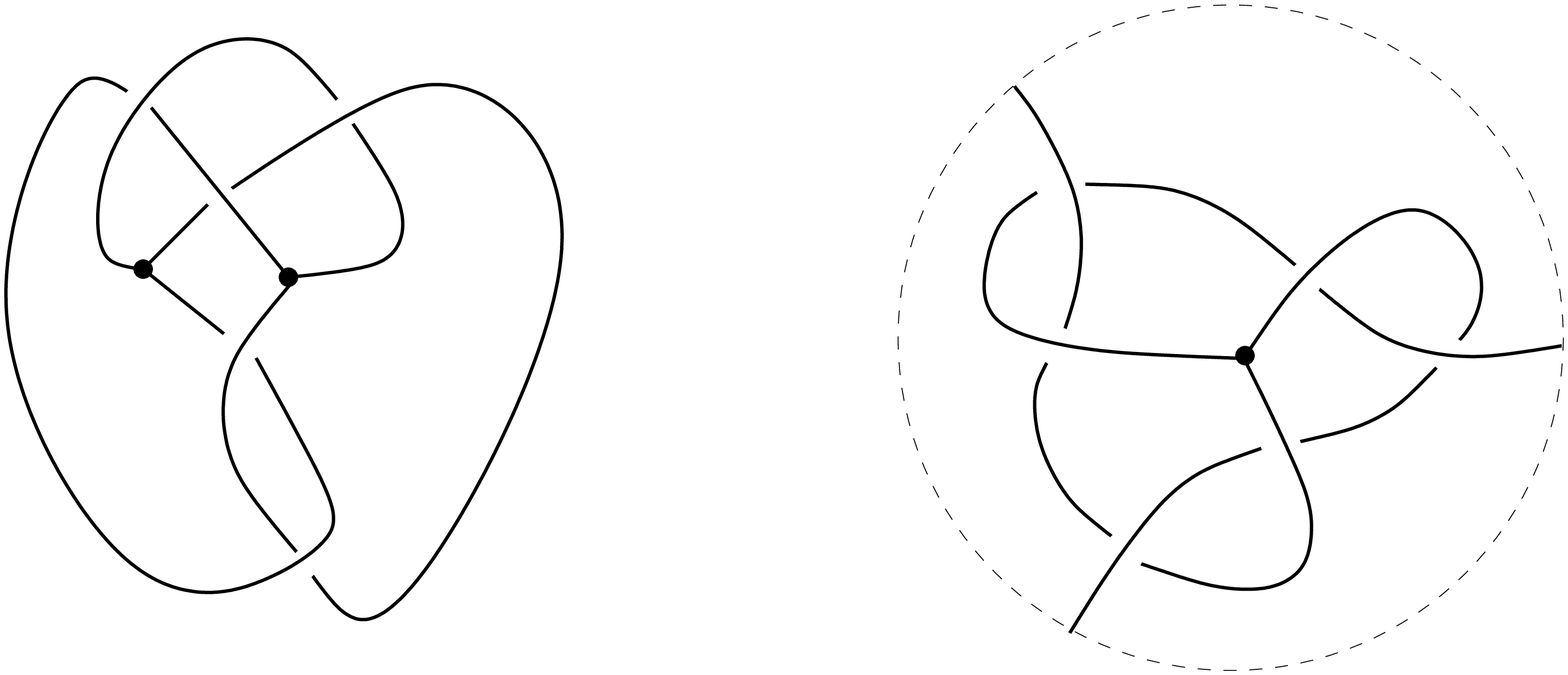}
%\vspace{-1em}
\caption{\label{kino} On the left, Kinoshita's $\theta$-graph. On the right,
Thurston's knotted wye: the second vertex is understood to be at infinity, in such a way
that the complement of a regular neighbourhood of the graph coincides with
the complement in a ball of a regular neighbourhood of the tangle here described.}
\end{center}
\end{figure}

The graph $\Gamma_K$ of Figure \ref{kino}  is the so-called
{\it Kinoshita $\theta$-graph}~\cite{kino}. It is the spine of
the spatial handlebody $H_K$, whose complement $\compl (H_K)$ will be denoted
by $M_K$. 

\smallskip

In \cite{kino} Kinoshita introduced some elementary ideals
$E_d(\Gamma,z)$ associated to any presentation of the Alexander module
of the fundamental group of $S^3\setminus \Gamma$. These ideals turn out to be
isotopy invariants for the couple $(\Gamma, z)$, where $\Gamma$ is a
spatial {\it graph} (not necessarily of genus 2) endowed with a given
$\Z$-cycle $z$. By means of these invariants he proved for example
that ``his'' graph $\Gamma_K$ is knotted.

The Kinoshita $\theta$-graph has several nice properties and since
\cite{kino} does occur in several papers in order to test, for
example, whether or not certain invariants are able to distinguish it
from the unknotted $\theta$-graph. An interesting property of $\Gamma_K$ is that
all its constituent knots are unknotted (one says that it is a
``minimally knotted'' graph), so the unknotting criterium of
\cite{scharlemann1} applies and one can conclude that:

\smallskip

{\it The Kinoshita's $\theta$-graph $\Gamma_K$ is knotted if and only
  if the handlebody $H_K$ is $(1)_S$-knotted}.  \medskip

This is a rather exceptional behaviour, because we know that in
general the knotting of a given spine does not imply the
$(1)_S$-knotting of the associated handlebody. On the other hand, as observed in~\cite{ushi},
Kinoshita's $\theta$-graph is isotopic to the ``knotted wye'' graph introduced
by Thurston in~\cite[Example 3.3.12]{thu}, where it is also shown that the  
manifold
$M_K=\compl(H_K)$ admits a hyperbolic structure with geodesic boundary 
(in fact, it turns out that $M_K$ is the hyperbolic $3$-manifold with geodesic boundary
of smallest volume~\cite{koji} and smallest complexity~\cite{FMP}).
This implies
that the boundary of $M_K$ is incompressible, hence Proposition \ref{1-K-char}
does apply, and $H_K$ is (at least) $(2)_S$-knotted. For these reasons
it is natural to look for the true level of knotting of $H_K$.
\smallskip

The first example of a genus 2 cube--with--holes $M=\compl (H)$ not
admitting any $(M\to W)$--boundary--preserving--map is due to Lambert
\cite{lambert}. In fact, Lambert's example coincides with the handlebody $H$ associated to the
spine of Figure \ref{handcuff3}, equivalently $H=H_1(3)$. The
proof in \cite{lambert} is of topological nature, based on an accurate
analysis of such a specific example.

\smallskip

The first example of a genus 2 cube--with--holes $M=\compl (H)$ 
having cut number equal to 1 is due to Jaco \cite{jaco}. 
The discussion of this example exploits
the following topological obstruction to be of corank 2: {\it If $M$
  has corank equal to 2, then every map $f: M\to S^1\times S^1$ is
homotopic to a non--surjective map}.  This obstruction looks (at least to
us) not so handy in order to discuss other examples.

\smallskip

In \cite{suzuki} the author remarked that Kinoshita's invariants (and
some variations of them) can be used to face questions concerning the
cut number of $M$ or the existence of $(M\to W)$--boundary--preserving--maps.
As applications he gave a different proof of the fact that Lambert's
example does not admit any $(M\to W)$--boundary--preserving--map, and proved that $\cut(M_K)=1$. In \cite{suzuki}
Suzuki adopts Kinoshita's set up, in terms of spatial graphs endowed with
$\Z$-cycles.
The fact that
one is actually working up to spine moves is somehow
implicit. Moreover, sometimes the proofs simply refer
to different sources strewn in the literature, where formally
analogous statements had been previously achieved in the classical
case of (boundary as well as homology boundary) links. For these
reasons we have preferred to provide below an essentially
self--contained, detailed account about these ``Alexander module
obstructions'', adopting as much as possible an intrinsic, geometric
approach.  By the way, we will point out the analogies but also some
remarkable differences between the case of genus 2
cubes--with--holes and the case of links. As usual we will limit to
deal with the case of present interest, although the discussion can be
generalized. A detailed and comprehensive account on Alexander modules of groups and spaces
is given in~\cite{hibook} (which is mainly concerned with links). Several algebraic results we are discussing here are proved in~\cite{hibook}
in greater generality. However, in order to make our treatment as elementary as possible, when this does not imply a big waste of space,
we provide proofs for the statements that are relevant to our purposes.

%\smallskip

%We will actually provide two different proofs of 
%Proposition~\ref{via-A-obs}. One is based on such Alexander
%obstructions. The other is based on a suitable topological
%characterization of the $M$'s admitting $\partial$-connected cut
%systems among the ones that admit a generic cut system. This
%characterization somehow refines the one established in \cite{jaco} and,
%in our opinion, is eventually easier to handle in comparison, for
%instance, with the original Lambert's topological treatment of his
%example. This matter will be developed in Section \ref{more-bpm}.

\smallskip

We stress that it is quite remarkable that {\it easily} computable
obstructions are able to recognize in some case, like $M_K$, that the
cut number is equal to 1. It is known that starting with the input of a
{\it finite presentation} of a group $G$, the determination of its
corank can be done {\it in principle} by means of an algorithm (see
\cite{mak}). This is an important conceptual fact, however the time of
execution a such a generic algorithm grows too fast with the input
complexity, so this is not of practical utility, also dealing with
rather simple examples.  In some cases (as already remarked in
\cite{Stallings}) one can associate to the finite presentation of $G$
some pertinent 3--dimensional (triangulated) manifold, and try to
exploit geometric/topological tools in order to simplify the
determination of the corank. Note that in our situation the problem is
3--dimensional from the very beginning, and one can try to use for
example the theory of {\it normal surfaces} in order to detect the
potential cut systems (if any). To this respect $M_K$ should appear
rather promising as it admits a very simple minimal triangulation as
well as simple presentations of the fundamental group. So we have
tried for a while to treat $M_K$ along this way. However, there is a
complication due to the fact that the theory of normal surfaces {\it
  with boundary} deals with many more elementary local configurations
with respect to the closed case. For example it is not hard to realize
in this way, {\it by bare hands}, that $M_K$ is ``small'', \emph{i.e.}~that it
does not contain any incompressible and not boundary parallel
 closed surface; on
the other hand, the computation of the cut number becomes rather
demanding, reasonably should need a computer aid, and eventually we
have preferred to turn towards more handy obstructions.
   
\subsection{The Alexander module}
Let us denote by $M$ either $M=\compl (H)$ or $M=\compl (L)$, $H$
being as usual a genus $2$ spatial handlebody, $L$ being a
$m$--component link (in fact we will refer mostly to the cases $m=1,\
2$). $G$ denotes the fundamental group of $M$, and $K$ the abelianized
of $G$, \emph{i.e.}~$K=G/[G,G]$.  We denote the operation of the
abelian group $K$ multiplicatively. By Hurewicz Theorem, there exists
a {\it canonical} isomorphism $K\cong H_1 (M)$. We  denote by
$\Lambda=\Z K$ the group ring of $K$.

Let $p^M_A\colon \widetilde{M}\to M$ be the {\it maximal abelian
  covering} of $M$, that is the one associated to $[G,G]$. Fix a
base point $x_0\in M$, and set $\widetilde{M}^0=(p_A^M)^{-1} (x_0)$.
The group of covering automorphisms of $\widetilde{M}$ is canonically
isomorphic to $K$, and acts on the pair
$(\widetilde{M};\widetilde{M}^0)$.  Hence $H_1
(\widetilde{M};\widetilde{M}^0)$ admits a natural structure of
$\Lambda$-module, denoted by $A(M)$, which is by definition the {\it
  Alexander module of $M$}. There are two important homomorphisms defined
on $\Lambda$: 

\smallskip

\begin{enumerate}

\item {\it (Augmentation map)} $\epsilon: \Lambda \to \Z$, which
  sends every element $\sum_i m_i k_i$, where $m_i\in\Z$ and $k_i\in
  K$, to the integer $\sum_i m_i$. Its kernel $I$ is called the {\it
    augmentation ideal} of $\Lambda$.

\item {\it (Canonical involution)} $\sigma: \Lambda \to \Lambda$,
  which sends every element $\sum_i m_i k_i$ to $\sum_i m_i k_i^{-1}$
  (the fact that $\sigma$ is a homomorphism relies on the fact that
  $K$ is abelian).
\end{enumerate}

\subsection{Finite presentations via free differential calculus}
\label{free-calc}
Let us briefly recall the definition of Fox derivative
\cite{Fox2}. Suppose $F=F(x_1,\ldots,x_n)$ is a free group on $n$
generators. For $j=1,\ldots, n$, the Fox derivative
$$
\partial_j \colon \Z F \to \Z F
$$
is defined as the unique $\Z$--linear map such that 
$$
{\partial_j (w\cdot w')}={\partial_j w}+w\cdot {\partial_j w'}\qquad
\forall w,w'\in F
$$
and
$$
{\partial_j x_j}=1,\quad  {\partial_j x_j^{-1}}=-x_j^{-1},\quad 
{\partial_j x_i^{\pm 1}}=0 \quad \forall i\neq j .
$$

Let $A$ be a $\Lambda$--module; $A$ is {\it finitely presented} if
it is  isomorphic, as a $\Lambda$--module, to the quotient
$\Lambda^n/\langle r_1,\ldots,r_s\rangle$, where each $r_i$ in an
element of $\Lambda^n$, and $\langle r_1,\ldots,r_s\rangle$ is the
$\Lambda$--module generated by the $r_i$'s.  If $B$ is the $(s\times
n)$--matrix whose rows are given by the $r_i$'s, then we say that $B$
is a {\it presentation matrix} for $A$.
\medskip

Let $M$, $G$ be as above. A proof of the following result can be found
for example in~\cite{lickorish}:

\begin{teo}\label{pres}
  Let $\langle x_1,\ldots,x_n \, |\, r_1,\ldots,r_s\rangle$ be a
  presentation of $G$, and let $j\colon F(x_1,\ldots,x_n)\to G$ and
  $k\colon G\to G/[G,G]=K$ be the natural projections.  Then a matrix
  presentation of $A(M)$ is provided by the matrix $(B)_{il}$, where
$$
B_{il}=k\left(j\left({\partial_l r_i}\right)\right),\qquad i=1,\ldots,
s,\ l=1,\ldots,n .
$$ 
\end{teo}

In particular, this implies that, up to isomorphism, $A(M):=A(G)$ only
depends on $G$, and it is finitely presented.

\subsection{Elementary ideals}\label{elementary-ideal}
We briefly recall here the definition of elementary ideal of a finitely presented
$\Lambda$--module. The interested reader can find a detailed account
on this issue for instance in~\cite{hibook}.

Let $A$ be a finitely presented $\Lambda$--module with a given
$(s\times n)$ presentation matrix $B$ as above. For every $d\in\matN$,
let us now define the ideal $E_d (B)\subseteq \Lambda$ as follows: 
\begin{itemize}
\item
If
$n-d>s$, then $E_d (B)=0$; 
\item
If $0< n-d\leq s$, then $E_d (B)$ is the ideal generated by the
determinants of the $(n-d)\times (n-d)$ minors of $B$;
\item
If $n-d\leq 0$, then $E_d (B)=\Lambda$.
\end{itemize}
We have therefore
$E_d(B)\subseteq E_{d+1} (B)$ for every $d\in\matN$.
\smallskip

A well-known purely algebraic result on presentations of modules implies
that $E_d (B)$ depends in fact only on $A$, so that for every $d\in\matN$
the {\it Alexander elementary ideal} of $A$
$$
E_d (A):=E_d (B)
$$
is well--defined.  The $k$--th {\it Alexander principal ideal
  $\Pp_k(A)$} is the smallest principal ideal containing $E_k (A)$.
Since $\Lambda$ is a UFD, a generator $\Delta_k (A)$
of $\Pp_k(A)$, which is usually called the \emph{$k$-th Alexander polynomial} of $A$) is simply the greatest common divisor of any set of generators for $E_k
(A)$.  

\medskip

We can apply these definitions to $A(M)$ so that the elementary ideals
$E_k(M)$ as well as $\Pp_k(M)$ and $\Delta_k (M)$ are well defined
topological invariants of $M$ (actually depending only on the
fundamental group $G$). In the case of a link $L$  we also write 
$A(L),\dots, \Delta_k(L)$, instead of $A(M),\dots,\Delta_k (M)$ .

The following Lemma is useful to study these invariants. Its proof
only relies on elementary algebraic arguments involving the very definition
of Alexander ideals.

\begin{lem}\label{esercizio}
\begin{enumerate}
\item Suppose $A$ is a finitely presented $\Lambda$--module. Then for
  every $d,k\in\matN$ we have
$$
E_{d+k} (A\oplus \Lambda^k)=E_d (A).
$$
\item
Suppose 
$$
\xymatrix{
0 \ar[r] & A_1 \ar[r]^i & A_2 \ar[r]^\pi & A_3 \ar[r] & 0 }
$$
is an exact sequence of $\Lambda$--modules. Then
$$
E_d (A_1) E_{d'} (A_3)\subseteq E_{d+d'} (A_2)\quad {\rm for\ every}\
d,d'\in\matN .
$$ 
\item
Suppose $B$ is a square presentation matrix for $A$
of order $n$, and
let 
$${\rm Ann} (A)=\{\gamma\in\Lambda \, |\, \gamma (a)=0\ \forall a\in A\}
\subseteq \Lambda$$ be the annihilator ideal of $A$.  Then we have
$$
E_0 (A)\subseteq {\rm Ann} (A).
$$
\end{enumerate}
\end{lem}
\begin{proof}
  Point~(1) is an easy consequence of the fact that a presentation
  matrix for $A\oplus\Lambda^k$ is obtained by adding $k$ null columns
  to a presentation matrix for $A$.

  Points~(2) and~(3) are stated respectively as Theorem~3.12--(1) and Theorem~3.1--(1)
in~\cite{hibook}.
\end{proof}

\subsection{Polynomial ideals}\label{polynomials} 
$K$ is {\it non}--canonically isomorphic to $\mathbb{Z}^m$ (where
$m=2$ when $M=\compl (H)$, and $m$ equals the number of components in
the case of links). Let $\langle t_1,\dots, t_m\rangle$ be the
multiplicative abelian group freely generated by the symbols
$t_1,\dots , t_m$.  Choosing an isomorphism $K\cong \langle t_1,\dots
t_m\rangle$ is equivalent to choosing a basis for ${\rm Hom}(H_1
(M),\Z)=H^1 (M)$: in fact, since $K=H_1 (M)$ is torsion--free of rank
$m$, then to each basis $l_1,\dots, l_m$ of $H^1(M)$ there is
associated the isomorphism sending every $z\in H_1 (M)$ to the
monomial $t_1^{l_1 (z)}\cdots t_m^{l_m (z)}$, and every isomorphism
$K\cong \langle t_1,\dots, t_m\rangle$ arises in this way.  Every such
an isomorphism canonically extends to a ring isomorphism between
$\Lambda$ and $\Z[t_1^{\pm 1},\dots , t_m^{\pm 1}]$.

\smallskip

Alexander-Lefschetz duality provides a canonical isomorphism $H^1
(M)\cong H_1 (H) \cong H_1 (\Ll)$, where $\Ll$ is either some $L_\Gamma
\in \Ll(H)$ or $L$ itself: every integral cycle in $H_1 (\Ll)$ defines
a cohomology class in $H^1(M)$ via the linking number.  Therefore, if
$c_1,\dots , c_m$ is a basis of $H_1 (\Ll)$, then the map $z\mapsto
t_1^{\lk (c_1,z)}\cdots t_m^{\lk (c_m,z)}$ provides an isomorphism
between $K$ and $\langle t_1,\dots, t_m\rangle$. If the components
$K_1,\dots, K_m$ of $\Ll$ are {\it ordered and oriented}, then we can
select the {\it distinguished basis} of $H_1 (\Ll)$ such that 
$c_j= \sum_j \delta_{ij}[K_j]$. Hence in the case of a ordered and oriented
link $L$ we have a distinguished isomorphisms between $K$ and
$\langle t_1,\dots, t_m \rangle$. On the other hand, one can prove the following:

\begin{prop} If $M=\compl(H)$ then, by varying the (hc)--spine $\Gamma$
as well as the ordering and the orientation on $L_\Gamma$,
every  isomorphism between $K$ and $\langle t_1,t_2\rangle$ arises from
the  distinguished basis of some $L_\Gamma$.
\end{prop}
\begin{proof}
Of course, it is sufficient to show that every basis of $H_1 (H)$ is represented by the
(ordered and oriented) components of the constituent link of some spine of $H$. However, it is well--known 
that the group of homeomorphisms of $H$ into itself acts transitively on the set of bases of $H_1 (H)$
(see \emph{e.g.}~\cite[Lemma 2.2]{birman}). Therefore, if $\Gamma$ is any spine of $H$ and $\mathcal{B}$
is a fixed basis of $H_1 (H)$, we may find a homeomorphism $\varphi\colon H\to H$ whose induced map in homology
takes to $\mathcal{B}$ the basis associated to $L_\Gamma$. This implies that $\mathcal{B}$
is the distinguished basis associated to the spine $\varphi (\Gamma)$ of $H$. 
\end{proof}

Summing up, in the case of an {\it ordered and oriented} link $L$ (in
particular when $L$ is an oriented knot) there is a canonical
identification of $E_k(M)$, $\Pp_k(M)$ as {\it polynomial ideals} of
$\Z[t_1^{\pm 1},\dots , t_m^{\pm 1}]$, also denoted by $\Lambda$, and
$\Delta_k(M)$ is called the {\it $k$--th Alexander polynomial} of
$L$.  In the case of $M=\compl (H)$ such polynomial invariants are
well defined only up to the natural action of $SL(2,\Z)$ on $\Z[t_1^{\pm
  1}, t_2^{\pm 1}]$.

\begin{remark}\label{confronto}{\rm 
    Let $M=\compl (H(\Gamma))$, $G$ be as usual. Let $c\in H_1
    (\Gamma)$ be a non--trivial primitive homology class. The map
    $z\mapsto t^{\lk (c,z)}$ defines a surjective homomorphism
    $\alpha_c\colon H_1 (M)\to \langle t\rangle \cong \Z$, which
    induces in turn a ring homomophism $\alpha_c\colon \Lambda\to
    \Z[t,t^{-1}]$. Since $\alpha_c$ is sujective, if $B$ is a
    presentation matrix for $A(M)$ and $\alpha_c (B)$ is the matrix
    obtained by applying $\alpha_c$ to every coefficient of $B$, then
    $\alpha_c (E_d (M))=E_d (\alpha_c (B))$ for every $d\in\matN$ (see
    also~\cite{kino,suzuki}). Let $p\colon \widehat{M}\to M$ be the
    covering associated to $\ker \alpha_c\circ\pi$, where $\pi\colon
    G\to H_1 (M)$ is the Hurewicz epimorphism.  Then, $\widehat{M}$ is
    an infinite cyclic covering whose automorphism group is
    canonically isomorphic to $\langle t\rangle$. A slight variation
    of Theorem~\ref{pres} implies that $\alpha_c (B)$ is a
    presentation matrix of $H_1 (\widehat{M},p^{-1} (x_0))$ as a
    $\Z[t,t^{-1}]$--module, so the ideals $\alpha_c (E_d (M))$ are
    related to the topology of $\widehat{M}$.  Starting from this
    consideration, in~\cite{cimasoni} (which is concerned with links)
    it is shown how one can deduce results about the ideals $E_d (M)$
    starting from the study of all the infinite cyclic coverings of
    $M$. Note however that the polynomial $\alpha_c (\Delta_d (M))$
    may be different from the generator $\Delta^c_d (M)$ of the
    smallest principal ideal containing $\alpha_c (E_d (M))$. More
    precisely, it is obvious that $\alpha_c (\Delta_d (M))$ divides
    $\Delta^c_d (M)$; however, if $E_d (M)$ is not principal, then it
    may happen that $\Delta^c_d (M)$ does not divide $\alpha_c
    (\Delta_d (M))$. For example, in the case when $G$ is the
    fundamental group of the complement of a $m$-component (oriented)
    link $L$, $m\geq 2$, it is not difficult to show that if $c$ is
    the cycle given by the sum of all the components of $L$, then
    $\Delta^c_d (L)=(t-1)\cdot \alpha_c (\Delta_d (L))$
    (see~\cite[Proposition 2.1]{cimasoni}).}
\end{remark}

\subsection{The absolute module of  $\widetilde M$} \label{A-vs-AA} 
    
Instead of considering the module $A(M)$, one may consider the module
$AA(M)=H_1 (\widetilde M)$ (the \emph{absolute} Alexander module of
$M$), as also $AA(M)$ admits a natural structure of
$\Lambda$--module. The long exact sequence of the pair
$(\widetilde{M},\widetilde{M}^0)$ in homology provides a sequence of maps
of $\Lambda$--modules which can be eventually written in the form

$$
\xymatrix{
0 \ar[r] &AA(M)\ar[r] &A(M)\ar[r] &\Lambda \ar[r]^\varepsilon& \Z\ar[r] &0
}
$$

or equivalently in the form

\begin{equation}\label{cromwell}
\xymatrix{
0\ar[r] &AA(M)\ar[r]& A(M)\ar[r] &I \ar[r] &0.
}
\end{equation} 

Note that under any fixed isomorphism $\Lambda\cong \Z[t_1^{\pm
  1},\dots, t_m^{\pm 1}]$ as in the previous Section, the augmentation
ideal $I$ corresponds to $(t_1-1,\dots , t_m-1)$.  The sequence~\eqref{cromwell}
is usually known as the \emph{Crowell sequence} for $A(M)$.
\medskip

\subsection{The case of knots.} When $M=\compl (L)$, $L$ being an oriented
{\it knot}, since $I$ is a principal ideal, then
$I\cong\Lambda=\Z[t^{\pm 1}]$, so the sequence~\eqref{cromwell}
splits, and we have
\begin{equation}\label{knots:eq}
A(L)\cong AA(L) \oplus \Lambda .
\end{equation}

Putting together equation~\eqref{knots:eq} and
Lemma~\ref{esercizio}--(1) we obtain the following equalities
between polynomial ideals
$$
E_{k+1} (A(L))=E_k (AA(L))\quad {\rm for\ every}\ k\geq 0 .
$$
Moreover, it is known that any Wirtinger presentation of $G$ has
deficiency one (see Definition~\ref{defin}), 
and defines therefore a presentation matrix $B$ for
$A(L)$ having $n$ rows and $n+1$ columns. Another property of
Wirtinger presentations implies that the determinant of the square
matrix obtained by omitting the $i$--th column of $B$ does not depend
(up to units in $\Z[t^{\pm 1}]$) on $i$, so $E_1 (A(L))$ is
principal. The generator of
$$
E_1 (A(L))=E_0 (AA(L))
$$
is the classical \emph{Alexander polynomial} of the knot.

\medskip

In our cases of interest $M=\compl (H)$ or $M=\compl (L)$, $L$ being a
link, the Crowell sequence does not necessarily split, so the
relation between $A(M)$ and $AA(M)$ is less direct. We will come back
later to this issue.

\subsection{Alexander obstructions}\label{A-ob}
In this Section we focus on $M=\compl (H)$.

\begin{defi}\label{intrinsic-prop} {\rm Let $J$ be an ideal of
    $\Lambda$.  

\begin{enumerate}

\item We say that $J$ is
      \emph{unitary} if $\varepsilon (J)=\Z$, where $\varepsilon
      \colon \Lambda\to\Z$ is the augmentation map introduced above.

\item We say $J$ is \emph{symmetric} if $\sigma (J)=J$, where $\sigma$
  is the canonical involution introduced above.

\end{enumerate}
}
\end{defi}

The proof of the following easy Lemma is omitted.

\begin{lem}\label{in-coord} Let us fix any isomorphism
  $\Lambda \cong \Z[t_1^{\pm 1},t_2^{\pm 1}]$ (associated to the
  distinguished basis of some $L_\Gamma$ as in Section
  \ref{polynomials}). Assume that the ideal $J$ is principal
  generated by the polynomial $f$.  Then: 
\begin{enumerate}
\item  $J$ is unitary if and only if $f(1,1)=\pm 1$.
\item
 $J$ is symmetric if and only if $f(t_1,t_2)=\pm t_1^{a_1}
  t_2^{a_2} f(t_1^{-1},t_2^{-1})$ for some $a_1,a_2\in\mathbb{Z}$.
\end{enumerate}
\end{lem}
\smallskip

The following Proposition collects the most important properties of the
elementary ideals of $M$, including those providing the promised ``obstructions''.

\begin{prop} \label{main-ob} Let $M=\compl(H)$, $H$ being a genus 2
spatial handlebody. Then
\begin{enumerate}
\item $E_0(M)=E_1(M)=0$.

\item $E_2(M)$ is unitary. For every $L_\Gamma \in \Ll(H)$, we have $E_2(M)
  \subseteq E_2(L_\Gamma)$.

\item $E_2(M)=E_1(AA(M))$.

\item If  $\crk(\pi_1(M))=2$, then  $E_2(M)$ is principal.

\item If $M$ admits a $\partial$-connected cut system (equivalently,
  $M$ admits a $(M\to W)$--boundary--preserving--map), then $E_2(M)$ is
  symmetric.
\end{enumerate}
\end{prop}

We devote most of the Section to the proof of this Proposition. 
More precisely, points~(1) and (2) are proved in Subsection~\ref{1proof}, point~(3)
in Subsection~\ref{2proof}, point~(4) in Proposition~\ref{principal2} (see Subsection~\ref{3proof}), and point~(5) in Corollary~\ref{suzuki:cor} (see Subsection~\ref{4proof}). 

\begin{remark}
{\rm
To our purposes, the most relevant results described in Proposition~\ref{main-ob} are points~(4) and~(5).
Since $\pi_1 (M)$ has deficiency 2 (see Lemma~\ref{def}), 
point~(4) can be deduced from~\cite[Theorem 4.3]{hibook}, which implies in fact the stronger
result that $E_2(M)$ is principal if and only if $\pi_1 (M)$ admits an epimorphism onto
$F_2/F_2''$, where $F_2=[F'_2,F'_2]$ and $F'_2=[F_2,F_2]$.

Moreover, in~\cite{suzuki} point~(5) is described as a consequence of the results of~\cite{guti}, which -- 
however -- are concerned only with links. 

We have thus decided to include here a detailed account on how the argument in~\cite{guti}
can be adapted to the case of handlebody complements. In order to achieve this, it is necessary to introduce a machinery
which allows us to give a self--contained proof also of point~(4) without a too big waste of space. 
}
\end{remark}

We begin by pointing out
some analogies and differences with respect to the case of links.
\smallskip 

\begin{enumerate}
\item[(a)] In the case of a link $L$, it is proved
in~\cite{traldi1} that for every $d\in\matN$ there exist natural
numbers $k,h$ such that
$$
E_d (AA(L))\cdot I^k\subseteq E_{d+1} (A(L)),\qquad E_{d+1}
(A(L))\cdot I^h\subseteq E_d (AA(L)).
$$
Since the only principal ideal containing $I$ is the whole ring
$\Lambda$, this implies that for every $d\in\matN$ we have the
equality of Alexander polynomials
$$
\Delta_d (AA(L))=\Delta_{d+1} (A(L)).
$$ 
\smallskip
\item[(b)] Let $L$ be any link.
The following facts are proved in~\cite{torresfox}:
\begin{itemize}
 \item 
If $L$ is a knot, then all the Alexander ideals
$E_d (L)$, $d\in\matN$
 are symmetric;
\item
If $L$ is a knot, then the Alexander polynomial $\Delta_1 (L)$
satisfies
$$
\Delta_1 (t)=t^n \Delta_1 (t^{-1}),
$$ 
where $n$ is even;
\item If $L$ has $m\geq 2$ components, then the \emph{first} Alexander
  polynomial of $L$ satisfies
$$
\Delta_1 (L)(t_1^{-1},\ldots,t_m^{-1})=-t_1^{k_1-1}\ldots t_m^{k_m-1}
\Delta_1 (t_1,\ldots,t_m),
$$
where $k_i$ has the same parity of the sum of the linking numbers of $L_i$
with all the other components of $L$.
\end{itemize}
A proof of these facts using Seifert surfaces can be found 
in~\cite{cimasoni}.

\smallskip

\item[(c)] If $L$ is a 2--component {\it homology boundary} link, then $E_2(L)$
is not necessarily principal. According to \cite{HILLM}, such an ideal is principal
if and only if 
the image of $H_1(\partial  \widetilde{M})$ 
in  $H_1(\widetilde{M})$ vanishes. This last condition holds whenever $L$ is a {\it boundary} link. These results are in sharp contrast with point~(3) of Proposition~\ref{main-ob} above.
\end{enumerate}

\subsection{Proofs of points~(1) and (2) of Proposition~\ref{main-ob}.}\label{1proof}
We begin with the following:
\begin{defn}\label{defin}
  The \emph{deficiency} of a finite presentation $\langle S\, |\,
  R\rangle$ of a group $G$ is equal to the difference between the
  number $|S|$ of generators and the number $|R|$ of relations of the
  presentation.  The deficiency of a finitely presented group $G$
  is the maximal deficiency of finite presentations of $G$ (note
  that the deficiency of a group may be negative).
\end{defn}

The following easy result is proved \emph{e.g.}~in \cite[Theorem 7]{kino}:

\begin{lem}\label{def}
The deficiency of the fundamental group $G$ of $M=\compl (H)$ is $2$.
\end{lem}

\smallskip

Putting together Lemma~\ref{def} and Theorem~\ref{pres} we deduce that
$A(M)$ admits a $n\times (n+2)$ presentation matrix. Clearly this
implies that $ E_0 (M)=E_1 (M)=0 $.
\medskip

Let us now show that $E_2(M)$ is unitary.  It is not difficult to show
that if $\langle x_1,\ldots,x_g | r_1,\ldots,r_s\rangle$ is a
presentation of $G$ and $a\colon F(x_1,\ldots,x_g)\to \{1\}$ is the
trivial homomorphism, then $a(\partial r_i/\partial x_j)\in
\Z[\{1\}]=\Z$ computes the sum of the exponents of $x_j$ in the word
$r_i$. It readily follows that if $B$ is a presentation matrix for
$A(M)$, then applying $\varepsilon$ to every element of $B$ we obtain
a presentation matrix for $G/[G,G]\cong \Z^2$. 
Since the presentation of $A(M)$ can be chosen of deficiency
2, say of the form $n\times (n+2)$, this implies that the GCD of the
minors $n\times n$ of $\varepsilon (B)$ has to be 1. Then $\varepsilon
(E_2 (M))=\Z$.

\smallskip

If $\Gamma$ is a spine of $H$ with constituent link $L_\Gamma$, then
$\pi_1 (S^3\setminus L_\Gamma)$ is obtained by adding to any presentation
of $\pi_1 (M)$ a relation representing the boundary of a $2$--handle dual to the isthmus
of $\Gamma$. It follows that  
a 
presentation matrix for $A(L_\Gamma)$ is obtained by adding a row 
to a
presentation matrix of $A(M)$, and this readily implies that $E_2(M) \subseteq E_2( L_\Gamma)$.

\medskip

These facts already say that $E_2 (M)$ is the first non--vanishing
Alexander ideal of $M$. It is therefore not surprising that it encodes
several geometric properties of $M$.

\subsection{A relation between the ideals of $A(M)$ and of $AA(M)$.}\label{2proof}
Let us prove now that $E_2 (M)=E_1 (AA(M))$. This is 
a consequence of the more general results 
proved
in~\cite{traldi1,traldi2}. 
For the sake of completeness, 
we describe here the proof in the case we are interested in.  
\smallskip

We begin by
computing the ideal $E_1 (I)$. We fix an identification $\Lambda\cong
\Z[t_1^{\pm 1}, t_2^{\pm 1}]$. Let us consider the exact sequence
\begin{equation}\label{I}
\xymatrix{0 \ar[r] & \Lambda\ar[r]^{\alpha_1} & \Lambda^2 \ar[r]^{\alpha_2} & 
I \ar[r] & 0},
\end{equation}
where $\alpha_1 (p)=((t_2-1)p, (1-t_1)p)$ and $\alpha_2
(p,q)=(t_1-1)p +(t_2-1)q$.  The presentation matrix for $I$ related to
this exact sequence is given by $(t_2-1\ 1-t_1)$, so $E_1 (I)=I$.

Let us now show that $E_2 (M)\subseteq E_1 (AA(M))$. Since $\Lambda$
and $\Lambda^2$ are free $\Lambda$--modules, we can arrange the short
exact sequence~\eqref{I} and the Crowell exact sequence for $A(M)$ in
the commutative diagram
$$
\xymatrix{ 0 \ar[r] & \Lambda\ar[r]^{\alpha_1} \ar[d]^\beta &
  \Lambda^2
  \ar[r]^{\alpha_2} \ar[d]^\gamma & I \ar[r] \ar[d]^{\rm Id} & 0\\
  0 \ar[r] & AA(M) \ar[r]^{\varphi} & A(M) \ar[r]^{\psi} & I \ar[r] &
  0 }
$$
It is now easy to check that the sequence
$$
\xymatrix{ 0 \ar[r] & \Lambda\ar[r]^-{k} & AA(M)\oplus \Lambda^2
  \ar[r]^-{h} & A(M) \ar[r] & 0 }
$$
is exact, where $k(p)=(-\beta (p), \alpha_1 (p))$ and $h(a,(p,q))=\varphi
(a)+\gamma (p,q)$.  We may now apply Lemma~\ref{esercizio}, thus
obtaining
$$E_2 (A(M))=E_2(A(M))\cdot E_1 (\Lambda)\subseteq E_3 (AA(M)\oplus \Lambda^2)
=E_1 (AA(M)).$$

In order to show the opposite inclusion, let us introduce the
following notation: if $J,J'$ are ideals of $\Lambda$, then we set
$$
(J : J')=\{\lambda\in\Lambda\, |\, \lambda\cdot J'\subseteq J\}.
$$

Since $E_2 (M)$ is unitary we have $E_2 (M)+I=\Lambda$. As a
consequence, for every ideal $J$ of $\Lambda$ we have
$$
E_2 (M)+J=(E_2 (M)+J)\cdot (E_2 (M)+ I)\subseteq E_2 (M)+J\cdot I
\subseteq E_2 (M)+J,
$$
whence 
$E_2 (M)+J=E_2 (M)+J\cdot I$. Applying this equality to the case
$J=(E_2 (M) : I)$ we obtain
$$E_2 (M)+ (E_2 (M) : I) = E_2 (M)+(E_2 (M) : I)\cdot I\subseteq E_2 (M)+ 
E_2 (M)=E_2 (M),$$
whence
\begin{equation}\label{bah}
E_2(M)=(E_2 (M) : I).
\end{equation} 

Now, by applying Lemma~\ref{esercizio}--(2) to the Crowell sequence
for $A(M)$ we have $E_1 (AA(M))\cdot I=E_1 (AA(M))\cdot E_1
(I)\subseteq E_2 (M)$. Together with equation~\eqref{bah}, this
finally implies
$$
E_1 (AA(M))\subseteq ((E_1 (AA(M))\cdot I) : 
I)\subseteq (E_2 (M) : I)=E_2 (M).
$$

\subsection{The case when $\cut (M)=2$}\label{3proof}
Assume now that $\cut (M)=2$.  Let $Y$ be the figure-eight
graph $S^1\vee S^1$ and denote by $y_0$ the singular point of $Y$. We
also choose a base point $x_0\in M$.  Then $G=\pi_1 (M,x_0)$ admits an
epimorphism onto $F_2=\Z \ast \Z=\pi_1 (Y,y_0)$.  Then, there exist
continuous maps $f\colon M\to Y$ and $g\colon Y \to M$ such that
$f\circ g$ is the identity of $Y$. We may also assume that
$f(x_0)=y_0$ and $g(y_0)=x_0$.

Let $p^Y_A\colon \widetilde{Y}\to Y$ be the covering associated to
$[F_2,F_2]$. Since $f_\ast ([G,G])\subseteq [F_2,F_2]$ and $g_\ast
([F_2,F_2])\subseteq [G,G]$, the maps $f$ and $g$ lift to continuous
maps $\widetilde f\colon \widetilde{M}\to\widetilde{Y}$, $\widetilde
g\colon \widetilde{Y}\to\widetilde{M}$ which can be chosen in such a
way that $\widetilde{f}\circ\widetilde{g}={\rm Id}_{\widetilde{Y}}$.

Let us now consider the relative homology group
$A(F_2)=H_1(\widetilde{Y};(p^Y_A)^{-1} (y_0))$.  Observe that $f_\ast$
induces an isomorphism between $K=G/[G,G]$ and $F_2/[F_2,F_2]$, which is in
turn canonically isomorphic to the group of covering automorphisms of
$\widetilde{Y}$. As a consequence, we may consider the induced action
of $K$ also on the pair $(\widetilde{Y};(p^Y_A)^{-1} (y_0))$.  We have
therefore that also $A(F_2)$ admits a natural structure of
$\Lambda$-module.

By construction, the maps $\widetilde{f}$ and $\widetilde{g}$ commute
with the action of $K$ on $\widetilde{M}$ and $\widetilde{Y}$, so
$\varphi:=\widetilde{f}_\ast\colon A(M)\to A(F_2)$ and $\psi :=
\widetilde{g}_\ast\colon A(F_2)\to A(M)$ provide morphisms of
$\Lambda$-modules such that $\psi\circ \varphi={\rm Id}_{A(F_2)}$.

$$
\xymatrix{
\widetilde{M} \ar@/^.5pc/[r]^{\widetilde{f}} \ar[d] & 
\widetilde{Y} \ar@/^.5pc/[l]^{\widetilde{g}} \ar[d]\\
M \ar@/^.5pc/[r]^f & Y \ar@/^.5pc/[l]^g ,
}\qquad 
\xymatrix{
\\
A(M)\ar@/^.5pc/[r]^\varphi & A(F_2)  \ar@/^.5pc/[l]^\psi .
}
$$

\begin{lem}\label{sommadiretta}
We have $A(M)\cong \Lambda^2 \oplus \ker\varphi$. 
In particular, $\ker\varphi$ is a finitely--presented $\Lambda$--module,
and $E_2 (M)=E_0 (\ker\varphi)$.
\end{lem}
\begin{proof}
  First observe that since $\psi\circ \varphi={\rm Id}_{A(F_2)}$ we
  have $A(M)\cong A(F_2)\oplus \ker \varphi$, so in order to prove the
  first assertion it is sufficient to observe that $A(F_2)\cong
  \Lambda^2$ (this is a consequence, for example, of Theorem~\ref{pres}). 

  Since $A(M)$ is finitely presented and $A(M)\cong \Lambda^2 \oplus
  \ker\varphi$, also $\ker\varphi$ is finitely presented, and now the
  conclusion follows from Lemma~\ref{esercizio}--(1).
\end{proof}

A result of Stallings ensures that the kernel of any surjection
 of $\pi_1 (M)$ onto $F_2$ does not depend on the chosen surjection 
(see Theorem~\ref{stallings} below). Therefore, 
$\ker\varphi\subseteq A(M)$ admits an intrinsic characterization, which
is actually independent from $\varphi$. However, we show 
in Proposition~\ref{principal2} below how 
this characterization can 
be obtained without relying on Stallings' results.
 
Let us define the {\it torsion submodule} $T(M)$ of $A(M)$ as follows:
$$
T (M)=\big\{a\in A(M)\, |\, \gamma (a)=0\ {\rm for\ some}\ \gamma\in
\Lambda\setminus\{0\}\big\}.
$$

The following Proposition implies point~(4) of
Proposition \ref{main-ob}.
\begin{prop}\label{principal2}
 We have $T(M) =\ker\varphi$. Moreover, 
\begin{enumerate}
\item
$T(M)=\ker \varphi$; 
\item 
$T(M)$ is finitely presented, and $A(M)\cong \Lambda^2\oplus T(M)$;
\item
$E_2 (M)=E_0 (T(M))$;
\item
$E_2 (M)$ is principal.
\end{enumerate}
\end{prop}
\begin{proof}
By Lemmas~\ref{sommadiretta} and~\ref{esercizio} we have that $\Delta_2 (M)=\Delta_0 (\ker\varphi)\in E_0 (\ker\varphi)\subseteq {\rm Ann} (\ker\varphi)$. Moreover, 
 since $E_2 (M)$ is unitary, the Alexander polynomial $\Delta_2
  (M)$ cannot be null, and these facts imply that $\ker\varphi$ is
  contained in $T(M)$. On the other hand, since $A(M)=\Lambda^2\oplus
  \ker\varphi$ we also have $T(M)\subseteq \ker\varphi$, whence~(1).
  Having proved (1), points (2) and (3) are simply a restatement of
  Lemma~\ref{sommadiretta}. By Lemma~\ref{def} and Theorem~\ref{pres},
  the module $A(M)$ admits a presentation of deficiency two. Since
  $A(M)\cong \Lambda^2\oplus \ker\varphi$, this readily implies that
  $\ker\varphi$ admits a square presentation matrix, and this
  implies in turn that $E_0 (\ker\varphi)$ is principal, whence (4).
\end{proof}

\subsection{The case when $M$ admits a $\partial$-connected cut system}\label{4proof}
The proof of the last point of Proposition \ref{main-ob} is a bit more
demanding and incorporates more geometric insight about $E_2(M)$.
We follow here the strategy described in~\cite{guti}, where the case of (complements of)
boundary links is treated.

Assume as above that $\cut (M)=2$, and let us fix a cut system
$\Ss=\{S_1,S_2\}$ on $M$. We can assume that the base point $x_0$ of
$M$ does not lie on the union $S_1\cup S_2$. Then, just as the maximal free covering
$\widetilde{M}_\omega$ introduced in Section~\ref{more-bpm}, also the covering
$\widetilde M$ admits a concrete description in terms of the topology
of $M\setminus (S_1\cup S_2)$.  

So, let $V$ and $S_i^\pm$ be defined as in Subsection~\ref{maximalfree:sub}, and let us
%Let us denote by $V$ the manifold with
%boundary obtained by cutting $M$ along $\calS$.  Then, the boundary of
%$V$ consists of some ``horizontal'' boundary region (given by
%$\partial M\cap \partial V$) and some ``vertical'' boundary region,
%coming from the cuts along $S_1$ and $S_2$.  More precisely, we fix an
%orientation on $S_i$, and call $S_i^+$, $S_i^-$ the vertical
%components of $\partial V$ associated to $S_i$, $i=1,2$, 
%in such a way that a positive basis of $S_i$ is completed to a positive
%basis of $M$ by adding a vector pointing towards $S_i^+$.
denote
by $k_i$ the element of $K=H_1 (M)$ satisfying $k_i\cap
S_j=\delta_{ij}$, $i,j=1,2$. Then, any loop representing $k_i$
and intersecting $S_i$ transversely in one point runs from $S_i^-$ to
$S_i^+$ in a regular neighbourhood of $S_i$, and from $S_i^+$ to
$S_i^-$ in $V$. Let us also consider the disjoint union of a countable
number $\{V_k\}_{k\in K}$ of copies of $V$, indexed by the elements of
$K=G/[G,G]$, and let us take the quotient of such a union under the
equivalence relation generated by
$$
V_k\ni x \ \sim \ y \in V_{k'} \quad \Longleftrightarrow \quad
k'=k_ik,\ x\in S_i^-\subseteq V_k,\, y\in S_i^+\subseteq V_{k'},\
{\rm and}\ x=y\ {\rm in}\ M\ .
$$

It is now easy to recognize that such a quotient is homeomorphic to
$\widetilde M$.  Also observe that the action of $K$ on $\widetilde M$
admits a very easy description: for every $k_0\in K$, the covering
translation associated to $k_0$ translates $V_k$ onto its copy $V_{k_0
  k}$, for every $k\in K$.

\smallskip

By Lemma~\ref{sommadiretta}, the ideal $E_2 (M)$ is equal to the ideal
$E_0 (\ker\varphi)$.  In what follows, we denote by $\widetilde\calS$
the set $(p_A^M)^{-1}(S_1\cup S_2)\subseteq \widetilde M$.

\begin{lem}\label{e2lemma}
 We have 
$$
\ker\varphi=T(M)\cong {\rm Im} \left(i_\ast\colon H_1 \left(\widetilde
    M\setminus \widetilde\calS\right)\to H_1 (\widetilde{M})\right).
$$
\end{lem}
\begin{proof}
  The Crowell sequences for the Alexander modules of $A(M)$ and
  of $A(F_2)$, together with the epimorphism $\varphi \colon G\to F_2$,
  gives rise to the following commutative diagram, where rows are
  exact and the last vertical arrow is an isomorphism:
$$
\xymatrix{
0 \ar[r] & H_1 (\widetilde M)=G/[G,G] \ar[r] \ar[d]^{\widetilde f_\ast} & 
A(M) \ar[r] 
\ar[d]^\varphi & \Lambda \ar[d]^\cong\\
0 \ar[r] & H_1 (\widetilde Y)=F_2/[F_2,F_2] \ar[r] & A(F_2) \ar[r] & \Lambda
}
$$
and this easily implies (by chasing the diagram) that
$T(M)=\ker\varphi$ is isomorphic to $\ker \widetilde{f}_\ast$.

Let us fix an identification of $K$ with $\mathbb{Z}^2$, set 
$$
V'=\bigsqcup_{i+j\ {\rm even}} V_{(i,j)} \subseteq \widetilde M,\quad
V''=\bigsqcup_{i+j\ {\rm odd}} V_{(i,j)}\subseteq \widetilde M,
$$
and observe that $V'\cap V''=\widetilde\calS$, $V'\cup V''=\widetilde M$.
Also observe that we have an obvious homotopy equivalence
between $\widetilde M\setminus\widetilde\calS$ and the \emph{disjoint} union
$V'\sqcup V''$. 

Also the covering $\widetilde Y$ of the figure-eight graph $Y$ admits
a similar decomposition. Putting together the Mayer--Vietoris
sequences relative to these decompositions, which are preserved by the
map $f\colon \widetilde M\to\widetilde Y$ introduced above, we get the
commutative diagram
$$
\xymatrix{
H_1 (V'\cap V'')\ar[r]^\theta & 
H_1 (\widetilde M\setminus\widetilde S)\ar[r]^{i_\ast} \ar[d] 
& H_1 (\widetilde M)=[G,G] 
\ar[r] \ar[d]^{\widetilde f_\ast} & 
H_0 (\widetilde\calS)\cong \Lambda^2 \ar[d]^\cong\\
& 0 \ar[r] & H_1(\widetilde Y)=F_2/[F_2,F_2] \ar[r] & \Lambda^2},
$$
where rows are exact, and the last vertical arrow is an isomorphism.
By chasing the diagram, it is now easy to show 
that $i_\ast (H_1(\widetilde M\setminus
\widetilde\calS))=\ker \widetilde{f}_\ast$, whence the conclusion.
\end{proof}

\begin{remark}\label{explicitpres:rem}
{\rm  The last diagram in the proof of Lemma~\ref{e2lemma} shows that
  $T(M)$ is isomorphic to the quotient of $H_1 (\widetilde
  M\setminus\widetilde S)$ by the image of $H_1 (V'\cap V'')$ via the
  map $\theta$. By definition, if $\widetilde S$ is the component of
  $\widetilde\calS$ separating $V_{(i',j')}\subseteq V'$ and
  $V_{(i'',j'')}\subseteq V''$, then $\theta$ maps every $\alpha\in
  H_1 (\widetilde S)$ into the difference $\theta'(\alpha)-\theta''
  (\alpha)$, where $\theta'$ (resp.~$\theta''$) is induced by the
  inclusion of $\widetilde S$ into $V'$ (resp.~into $V''$).  This
  remark will prove useful in Theorem~\ref{e2teo} for explicitly
  constructing a presentation of $T(M)$.}
\end{remark}

From now on, we denote by $W$ a regular neighbourhood of the set
$H\cup S_1\cup S_2$. Henceforth, we also make the

\smallskip

\noindent {\bf Standing assumption:}
The cut system $\Ss$ is $\partial$-connected.

\smallskip

We stress that this assumption plays a fundamental r\^ole in several
arguments below.

\begin{lem}\label{generators}
The inclusion
$S_1\cup S_2 \hookrightarrow W$ induces an isomorphism
$$
H_1 (S_1)\oplus H_1 (S_2) \cong H_1 (W).
$$
\end{lem}
\begin{proof}
The decomposition $W=H\cup (S_1\cup S_2)$ provides 
the Mayer--Vietoris sequence
$$
\xymatrix{
H_1 (\partial S_1\cup\partial S_2)\ar[r] & H_1 (S_1\cup S_2) 
\oplus H_1 (H)\ar[r] &
H_1 (W) \ar[r] & H_0 (\partial S_1\cup\partial S_2)}
$$
It is readily seen that our assumptions on $\partial S_i$, $i=1,2$,
imply that the last arrow is the zero map, while the first one has
exactly $H_1 (H)$ as image. This implies the conclusion.
\end{proof}

Let now $g_i$ be the genus of $S_i$, and let
$\beta_1^i,\ldots,\beta_{2g_i}^i$ be a basis of $H_1 (S_i)$. Observe
that $V=S^3\setminus W$, and let us define $(\beta_i^j)^+\in H_1 (V)$
(resp.~$(\beta_i^j)^-\in H_1 (V)$) as the element obtained by
``pushing'' $\beta_i^j$ on the positive (resp.~negative) side of $S_j$
into $V$.

We define the Seifert matrices $A^{11}$, $A^{12}$, $A^{21}$, $A^{22}$
as follows:
$$
(A^{hk})_{ij}=\lk ((\beta_i^h)^-,\beta_j^k),\qquad h,k\in\{1,2\},\
i=1,\ldots,2g_h,\ j=1,\ldots,2g_k ,
$$
where $\lk$ is the linking number in $S^3$.  Not that by the symmetry
properties of the linking number we have ${}^{t} A^{12}=A^{21}$.

Let us fix the identification $\Lambda\cong \Z[t_1^{\pm 1},t_2^{\pm 1}]$ which 
carries the element $k_i\in K$ dual to $S_i$ to the variable $t_i$, $i=1,2$.
 
\begin{lem}\label{e2teo}
The module $T(M)$ admits the presentation matrix 
$$
B=\left( \begin{array}{cc}{}^{t} A^{11} -t_1 A^{11} & (1-t_1) A^{12}\\
(1-t_2) A^{21} & {}^{t} A^{22} -t_2 A^{22}\end{array}\right) .
$$
\end{lem}
\begin{proof}
  Since $V=S^3\setminus W$, by Lemma~\ref{generators} (and Alexander
  duality) there exists a basis $\gamma_1^1,\ldots,
  \gamma_{2g_1}^1,\gamma_1^2,\ldots,\gamma_{2g_2}^2$ of $H_1 (V)$
  which is dual to the $\beta_i^j$'s, in the sense that $\lk
  (\beta_i^j,\gamma_h^k)=\delta_{ih}\delta_{jk}$, where $\lk$ is the
  linking number in $S^3$. Of course, the $\gamma_i^j$'s generate $H_1
  (\widetilde{M}\setminus \widetilde{\calS})$ as a
  $\Lambda$--module. By Remark~\ref{explicitpres:rem}, in order to
  obtain a presentation of $T(M)$ we have now to add the relations
$$
t_i \cdot (\beta_j^i)^+ =(\beta_j^i)^- ,
$$
written in terms of the $\gamma_j^i$'s. The conclusion follows.
\end{proof}

The following Corollary, together with Lemma \ref{in-coord}, achieves
the proof of the last point of Proposition \ref{main-ob}.

\begin{cor}\label{suzuki:cor}
  Given a $\partial$-connected cut system $\Ss=\{S_1,S_2\}$ of $M$,
  and fixing the identification of $\Lambda$ with $\Z[t_1^{\pm
    1},t_2^{\pm 1}]$ associated to $\Ss$, the principal ideal
    $E_2(M)$ is generated by a polynomial $\Delta_2 (M)$ satisfying
    the condition:
$$
\Delta_2 (M)(t_1^{-1},t_2^{-1})=t_1^{-m} t_2^{-n} \Delta_2 (M)(t_1,t_2),
$$ 
where $m,n$ are even.
\end{cor}
\begin{proof}
  Let $B$ be the matrix described in the statement of
  Lemma~\ref{e2teo}, and denote by $B^-$ the matrix obtained by
  replacing in $B$ every occurrence of $t_i$ with $t_i^{-1}$, $i=1,2$.
  By Lemma~\ref{e2teo} we have $\Delta_2 (M)(t_1,t_2)=\det B$,
  $\Delta_2 (M)(t_1^{-1},t_2^{-1})=\det B^-$. Moreover, $B$ can be
  obtained from $B^-$ by performing the following operations:
  multiplication of the first $2g_1$ rows by $t_1$; multiplication of
  the last $2g_2$ rows by $t_2$; transposition; multiplication of the
  whole matrix by $-1$ (this operation does not change the
  determinant, since the matrix has even order); multiplication of the
  first $2g_1$ columns by $(1-t_2)$; multiplication of the last $2g_2$
  columns by $(1-t_1)$; division of the first $2g_1$ rows by
  $(1-t_2)$; division of the last $2g_2$ rows by $(1-t_1)$. It readily
  follows that
$$
\Delta_2 (M)(t_1^{-1},t_2^{-1})=t_1^{-2g_1} t_2^{-2g_2} 
\Delta_2 (M)(t_1,t_2).
$$ 
\end{proof}

The proof of Proposition \ref{main-ob} is now complete.

\subsection{The elementary ideals associated to the graphs $\Gamma_1 (p)$}\label{lambertcomput}
For every odd prime $p$, let us set $M_1 (p)=\compl (H_1 (p))$, where
$\Gamma_1 (p)$ is the graph described in Figure~\ref{1-knotted}.
This Subsection is devote to the proof of the following:

\begin{prop}\label{via-A-obs} 
For every odd prime $p$, the  ideal $E_2 (M_1 (p))$
is principal, but not symmetric. Therefore, by Corollary~\ref{suzuki:cor},
the $M_1(p)$ does not
  admit any $(M\to W)$--boundary--preserving--map, or equivalently it does not
  admit any $\partial$-connected cut system. Hence $H_1(p)$ is
  $(3)_S$-knotted.
  \end{prop}
\begin{proof}
We begin by providing a presentation of the fundamental group
of $M_1(p)$, via the well--known Wirtinger's procedure. 
So, let us consider the diagram $\Dd$ of $\Gamma_1 (p)$ described in Figure~\ref{wirtinger:fig}, 
and let us denote by $a_i$ (resp.~$b_i,c_i,d_i$) the element
of the fundamental group represented by a loop based above the plane
of the diagram and positively encircling the arc labelled by $a_i$ (resp.~$b_i,c_i,d_i$)
(recall that an arc is an embedded curve in $\Dd$ having as
endpoints either an under-crossing or a vertex of $\Gamma_1(p)$).

\begin{center}\label{wirtinger:fig}
\begin{figure}[h!]
\input{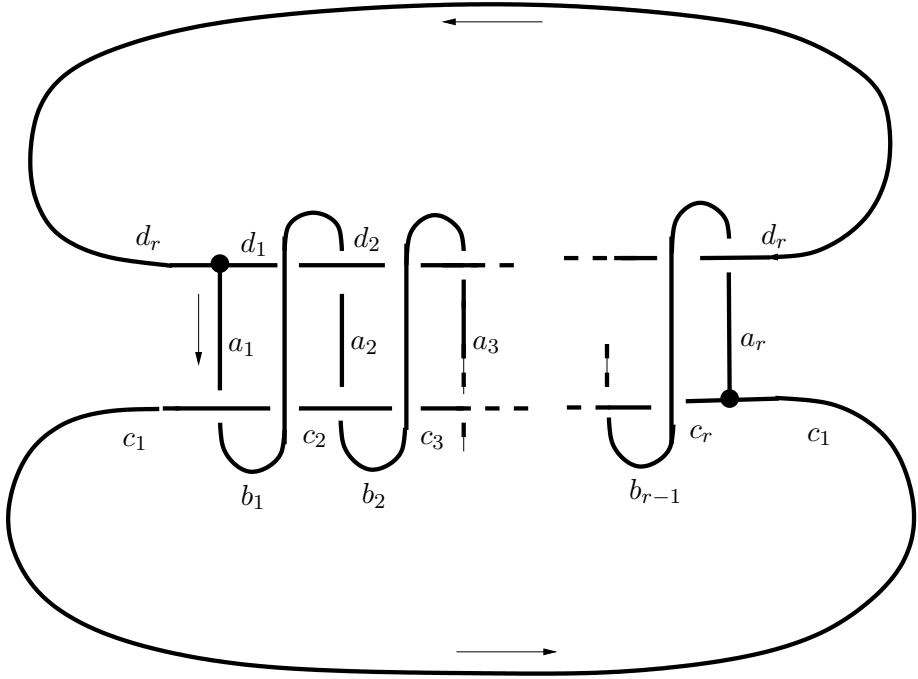}
\caption{Notation for a Wirtinger presentation of $\pi_1 (M_1(p))$.}
\end{figure}
\end{center}

Let us set $r=(p+1)/2$. 
The group $G=\pi_1 (M_1(p))$ has a presentation with generators
$a_i,b_j,c_i,d_i$, $i=1,\ldots,r$, $j=1,\ldots,r-1$, and relations
arising at every crossing and every vertex of $\Dd$. More precisely, the two vertices
determine the relations $c_1c_r^{-1}a_r$ and
$d_1d_r^{-1}a_1$, while looking at the crossings we obtain the relations
\begin{equation}\label{relations}
%\begin{array}{l}
c_ib_ic_i^{-1}a_i^{-1},\quad 
b_ic_ib_i^{-1}c_{i+1}^{-1},\quad 
b_id_ib_i^{-1}d_{i+1}^{-1},\quad
d_{i+1}b_id_{i+1}^{-1}a_{i+1}^{-1},%\ i=1,\ldots,r-1,
%c_1c_r^{-1}a_r,\ 
%d_1d_r^{-1}a_1
%\end{array}
\end{equation}
for $i=1,\ldots,r-1$.
Let us denote by $F_{4r-1}$ the free group generated by the symbols 
 $a_i,b_j,c_i,d_i$, $i=1,\ldots,r$, $j=1,\ldots,r-1$, and by
$j\colon F_{4r-1}\to G$ and $k\colon G\to G/[G,G]$
the natural projections. 
Then we have $G/[G,G]\cong H_1(M_1(p))\cong \mathbb{Z}^2$.
Moreover, it is readily seen that we may choose generators $t,s$ of
$H_1(M_1(p))$ in such a way that $k(j(c_i))=t$, $k(j(d_i))=s$
and $k(j(a_i))=k(j(b_i))=1$ for every $i$ (we denote here by $1$ the identity
of $H_1(M_1(p))$, which is considered as a multiplicative group).
By Theorem~\ref{pres}, the module $A(G)=A(M_1(p))$ admits a presentation with generators
$\overline{a}_i, \overline{b}_j,\overline{c}_i,\overline{d}_i$, $i=1,\ldots,r$,
$j=1,\ldots,r-1$, and relations which can be deduced by from the relations 
defining $G$ via Fox differential calculus.

Now, if $x,y,z$ belong to the fixed set of generators of $F_{4r-1}$, then
we have
$$
\begin{array}{cll}
k\left(j\left(\partial_x xyx^{-1}z^{-1}\right)\right)&=&1-k(j(y)),\\
k\left(j\left(\partial_y xyx^{-1}z^{-1}\right)\right)&=&k(j(x)),\\
k\left(j\left(\partial_z xyx^{-1}z^{-1}\right)\right)&=&-k(j(y))k(j(z))^{-1},%\\
%k\left(\partial_x y^{-1}z\right)&=&1,\\
%k\left(\partial_y y^{-1}z\right)&=&-k(j(x))k(j(y))^{-1},\\
%k\left(\partial_z xy^{-1}z\right)&=&k(j(x))k(j(y))^{-1},
\end{array}
$$
so the relations~\eqref{relations} induce the following relations for 
$A(G)$:
\begin{equation}\label{relations2}
\overline{b}_i=t^{-1}\overline{a}_i,\quad \overline{c}_{i+1}=(1-t)\overline{b}_i+\overline{c}_i,\quad
\overline{d}_{i+1}=(1-s)\overline{b}_i+\overline{d}_i,\quad \overline{a}_{i+1}=s\overline{b}_i ,
\end{equation}
where $i=1,\ldots,r-1$.
Moreover, the relations arising from the vertices of $\Dd$ give rise to the relations
\begin{equation}\label{relations3}
\overline{c}_r-\overline{c}_1=\overline{a}_r,\qquad
\overline{d}_r-\overline{d}_1=\overline{a}_1.
\end{equation}
It is readily seen that equations~\eqref{relations2} imply that 
\begin{equation}\label{relations4}
\begin{array}{llll}
\overline{a}_i&=&\left(st^{-1}\right)^{i-1}\overline{a}_1,\quad & i=1,\ldots,r,\\
\overline{b}_i&=&t^{-1}\left(st^{-1}\right)^{i-1}\overline{a}_1,\quad & i=1,\ldots,r-1,\\
\overline{c}_i&=&\overline{c}_1+(1-t)t^{-1}\left(\sum_{j=0}^{i-2} (st^{-1})^j\right) \overline{a}_1,\quad & i=1,\ldots,r,\\
\overline{d}_i&=&\overline{d}_1+(1-s)t^{-1}\left(\sum_{j=0}^{i-2} (st^{-1})^j\right) \overline{a}_1,\quad & i=1,\ldots,r .
\end{array}
\end{equation}
Putting together the relations~\eqref{relations4} and~\eqref{relations2}, we obtain that a presentation for $A(G)$
is given by the generators $\overline{a}_1,\overline{c}_1,\overline{d}_1$ with relations
\begin{equation*}%\label{relations5}
\left( (1-t)t^{-1} \left(\sum_{j=0}^{r-2} (st^{-1})^j\right) -(st^{-1})^{r-1}\right)\overline{a}_1=0,
\end{equation*}
\begin{equation*}%\label{relations5}
\left( (1-s)t^{-1} \left(\sum_{j=0}^{r-2} (st^{-1})^j\right) -1\right)\overline{a}_1=0.
\end{equation*}
As expected (since every fixed relation of a Wirtinger presentation is a consequence of the other
ones), these last two conditions are equivalent, and 
$$
A(G)=\Lambda^2 \oplus \Lambda / (f(s,t)),\qquad f(s,t)=(1-s)t^{-1} \left(\sum_{j=0}^{r-2} (st^{-1})^j \right)-1.
$$
By Lemma~\ref{esercizio}, we have $E_2(A(G))=E_0 (\Lambda / (f(s,t)))=(f(s,t))$, so in order to conclude
we only have to show that the ideal $(f(s,t))$ is not symmetric, \emph{i.e.}~that the ideals generated
by $f(s,t)$ and by $f(s^{-1},t^{-1})$ do not coincide. 

However, it is readily seen that every non--null principal ideal of $\Lambda=\mathbb{Z}[s^{\pm 1},t^{\pm 1}]$ 
admits a preferred generator which lies in $\mathbb{Z}[s,t]$ and 
is not divisible (in $\mathbb{Z}[s,t]$) by $s$ or $t$. Such a generator is unique up to sign.
Now, the preferred generators of $(f(s,t))$ and of $(f(s^{-1},t^{-1}))$ are respectively
$$
f_1(s,t)=t^{r-1}f(s,t)=(1-s) \left(\sum_{j=0}^{r-2} s^jt^{r-j-2} \right)-t^{r-1},
$$
$$
f_2(s,t)=s^{r-1}f(s^{-1},t^{-1})=(s-1)t\left(\sum_{j=0}^{r-2} s^jt^{r-j-2} \right)+s^{r-1}.
$$
Since $f_1(s,t)\neq \pm f_2(s,t)$, the proof of point~(1) of Proposition~\ref{via-A-obs} is complete.
\end{proof}

\subsection{The elementary ideals of the Kinoshita graph complement}\label{kinocomput}
This Subsection is devoted to the proof of the following:

\begin{prop}\label{via-A-obs2}
The elementary ideal $E_2(M_K)$ of Kinoshita's manifold $M_K$ is not principal. Therefore,
$\cut (M_K)=1$.
\end{prop}
\begin{proof}
As showed in~\cite{thu}, the manifold $M_K$ admits an ideal triangulation with two tetrahedra.
By analyzing the combinatorial structure of such a triangulation it is immediate to write down the following presentation
for the fundamental group $G$ of $M_K$ (such a presentation is much simpler than any Wirtinger
presentation associated to a diagram of $K$):  
$$
\langle
x_1,x_2,x_3\, |\, x_1x_2x_1^{-1}x_3x_1x_3^{-1}x_2x_3x_2^{-1}
\rangle.
$$
If we denote by $r$ the above relation,
then we have
$$
\begin{array}{lll}
{\partial_1 r}&=& 1-x_1x_2x_1^{-1}+x_1x_2x_1^{-1}x_3,\\
{\partial_2 r}&=& x_1+x_1x_2x_1^{-1}x_3x_1x_3^{-1}-x_1x_2x_1^{-1}x_3x_1x_3^{-1}
x_2x_3x_2^{-1},\\
{\partial_3 r}&=& x_1x_2x_1^{-1}-x_1x_2x_1^{-1}x_3x_1x_3^{-1}+x_1x_2x_1^{-1}
x_3x_1x_3^{-1}x_2 .
\end{array}
$$
Moreover, if $j\colon F(x_1,x_2,x_3)\to G$ and $k\colon G\to G/[G,G]$
are the natural projections, then we have $G/[G,G]\cong \langle
t_1,t_2\rangle $, where $t_1=k(j(x_1))$, $t_2=k(j(x_2))$ and
$k(j(x_3))=t_1^{-1}t_2^{-1}$.  Therefore
$$
k\left(j\left({\partial_1 r}\right)\right)= 1+t_1^{-1}-t_2,\quad
k\left(j\left({\partial_2 r}\right)\right)= -1+t_1+t_1t_2,\quad
k\left(j\left({\partial_3 r}\right)\right)= t_2(1-t_1+t_1t_2).
$$

Let us now compute $E_2 (M_K)$. Let $p_i=k(j(\partial_i r))$, $i=1,2,3$,
and observe that $E_2 (M_K)=(p_1,p_2,p_3)$. Since $t_1,t_2$ are
invertible in $\Lambda$, if $p_1'=t_1p_1=1+t_1-t_1t_2$,
$p_3'=t_2^{-1}p_3=1-t_1+t_1t_2$, then $E_2(M_K)=(p_1',p_2,p_3')$. 
Moreover, we have $2=p_1'+p_3'\in E_2 (M_K)$, and
$p_2=p_1'+(t_1t_2-1)\cdot 2$, so
$$
E_2 (M_K)=(2,p_1')=(2,1+t_1-t_1t_2).
$$

Since $2$ is irreducible in $\Lambda$ and does not divide
$1+t_1-t_1t_2$, we have $\Delta_2 (M_K)=1$. Therefore, if $E_2(M_K)$ were
principal, there would exist Laurent polynomials $q_1,q_2\in\Z[t_1^{\pm
  1},t_2^{\pm 1}]$ such that
$$
(1+t_1-t_1t_2)\cdot q_1 (t_1,t_2)+2q_2(t_1,t_2)=1,
$$
whence 
$$
(1+t_1+t_1t_2)\cdot \overline{q}_1 (t_1,t_2)=1\ {\rm in}\qquad \mathbb{Z}_2[t_1^{\pm 1},t_2^{\pm 1}]
$$
for a Laurent polynomial $\overline{q}_1\in\Z_2[t_1^{\pm 1},t_2^{\pm 1}]$.
Since $(1+t_1+t_1t_2)$ is not invertible in $\Z_2[t_1^{\pm 1},t_2^{\pm 1}]$, this gives a contradiction, hence
$E_2 (M_K)$ is not principal. Then it follows from Proposition 
\ref{main-ob} that the cut number of $M_K$ is not equal to 2.
\end{proof}

\subsection{Infinitely many $(4)_L$-knotted handlebodies}\label{infinitelym}
We now show how starting from a $M=\compl (H)$ such that
$\cut(M)=1$ (like
$M_K$) we can construct infinitely many other examples having cut number
equal to 1. Let
$\Gamma$ be a spine of $H$. Let $B\subseteq S^3$ be a 3--ball whose
boundary intersects transversely $\Gamma$ at two regular points, and
such that $(B,\ell)$, $\ell := B\cap \Gamma$, is an unknotted 1--1 tangle.
We call such an $\ell$ an {\it untangled arc} of $\Gamma$. We now replace $(B,\ell)$ by a 1--1 tangle $(B,T)$ such that $T$ and
$\ell$ have the same end--points, thus obtaining a new graph $\Gamma'$. If $K$ is the knot obtained as the closure of $(B,T)$ in $S^3$,  
then we denote the spatial graph $\Gamma'$ by the symbol $\Gamma\#_\ell
K$. If orientations on $K$ and on $\Gamma$ are not specified, then
there are nonequivalent ways of performing the ``sum'' $\Gamma\#_\ell
K$, but to our purposes this is not relevant.

\begin{lem}\label{newcorank}
Suppose that $\crk (\pi_1 (S^3\setminus \Gamma))=1$. Then, for
every sum $\Gamma\#_\ell
K$ as above, we have
$\crk( \pi_1 (S^3\setminus (\Gamma\#_\ell K)))=1$.
\end{lem}
\begin{proof}
  Let us set $G_\Gamma=\pi_1 (S^3\setminus \Gamma)$, $G_K=\pi_1
  (S^3\setminus K)$, $G_{\Gamma\# K}=\pi_1 (S^3\setminus
  (\Gamma\#_\ell K))$, and let $\psi\colon G_{\Gamma\# K}\to F_2$ be a
  homomorphism.  We will show that $\psi$ is not surjective.
 
 Let $(B,\ell)$ and $(B,T)$ be the tangles involved in the construction of
 $\Gamma\#_\ell K$, as above. 
 Then, the set 
$\Omega_1=S^3\setminus (\Gamma\cup B)$ is a deformation retract of $S^3\setminus \Gamma$, and intersects
a regular neighbourhood $\Omega_2$ of $B\setminus T$ in an annulus. Also observe that $\Omega_1\cup\Omega_2=S^3\setminus
(\Gamma\#_\ell K)$, so we have
the following commutative diagram 
$$
\xymatrix{
& G_\Gamma \ar[rd]_{j_1} &\\
%\ar@/^1pc/[rrd]^{\psi_1} &\\
\Z\ar[ur]_{i_1} \ar[dr]^{i_2} & & G_{\Gamma\#_\ell K}\ar[r]^\psi & F_2\\
& G_K\ar[ur]^{j_2} %\ar@/_1pc/[urr]_{\psi_2},
}
$$
where $i_2$ maps a generator of $\mathbb{Z}$ onto the class of a meridian of $S^3\setminus K$.
Since $G_K$ has corank one and every subgroup of a free group is free,
the group $\psi (j_2(G_K))$ is cyclic (possibly trivial), whence
abelian. This implies that $\psi \circ j_2$ factors through $H_1
(S^3\setminus K)\cong \Z$.  Since $i_2 (\Z)<\pi_1 (G_K)$
isomorphically maps onto $H_1 (S^3\setminus K)=G_K/[G_K,G_K]$, this
implies that
\begin{equation}\label{cont}
  \psi (j_2 (G_K))=\psi(j_2(i_2 (\Z)))=\psi (j_1(i_1(\Z))) 
  < \psi (j_1 (G_\Gamma)).
\end{equation}
By Van Kampen's Theorem, the group $G_{\Gamma\#_\ell K}$ is generated
by $j_1 (G_\Gamma)\cup j_2(G_K)$. Together with~\eqref{cont}, this
implies that $\psi (G_{\Gamma\#_\ell K})=\psi_1 (G_\Gamma)$.  But
$\psi_1 (G_\Gamma)$ cannot be the whole $F_2$ since the corank of
$\psi_1 (G_\Gamma)$ is one, so $\psi$ is not surjective.
\end{proof}

Now, let $H_1,\ldots,H_n$ be non--isotopic handlebodies whose
complements in $S^3$ have cut number equal to 1. We show how to
construct inductively a further handlebody $H_{n+1}$ 
such that $\cut(\compl (H_{n+1}))=1$ and $\compl (H_{n+1})$ is not homeomorphic to $\compl(H_i)$, $i=1,\ldots,n$
(in particular, $H_{n+1}$ is not  
isotopic to $H_i$, $i=1,\ldots,n$).
Starting for example from $M_K$, this procedure 
provides the desired infinite family of examples with cut number equal to 1, thus concluding the proof 
of Theorem~\ref{nonempty}.

So, let $\Gamma$ be a (hc)--spine of $H_n$, having
$K_1,K_2$ as constituent knots. Let $\ell$ be an untangled sub-arc of
$K_1$ which does not contain any vertex of $\Gamma$.  Let also $N$ be
the maximum of the ranks of the fundamental groups of $S^3\setminus
H_1,\ldots,S^3\setminus H_n$ (recall that the rank of a group is the
minimal number of generators in any presentation of the group), and take
a knot $K_3$ which is the sum of $N+1$ non--trivial knots.
Define $H_{n+1}$ as the regular neighbourhood of $\Gamma\#_\ell K_3$.
By Lemma~\ref{newcorank}, the complement of $H_{n+1}$ has cut number equal
to 1. Moreover, by compressing $S^3\setminus H_{n+1}$ along the
meridian disk of $H_{n+1}$ dual to $K_2$, we obtain $S^3\setminus
N(K_1\# K_3)$. This readily implies that the rank of $\pi_1
(S^3\setminus H_{n+1})$ is not smaller than the rank of
$\pi_1(S^3\setminus N(K_1\# K_3))$. But $K_1\# K_3$ is the sum of at
least $N+1$ non--trivial knots, so ${\rm rk}\, \pi_1(S^3\setminus
N(K_1\# K_3))\geq N+1$ by~\cite{weid}. In particular, $\pi_1
(S^3\setminus H_{n+1})$ is not isomorphic to $\pi_1 (S^3\setminus
H_i)$ for every $i=1,\ldots,n$.

\medskip

The proof of Theorem~\ref{nonempty}
is now complete.

\section{Perspectives}\label{where}
Let us summarize what we have achieved so far. In Section
\ref{levels} we asked some questions concerning the relations
that hold among  
extrinsic and intrinsic definitions of knotting for spatial handlebodies.
The results proved in our paper provide a quite clear picture of the situation. However,
two questions still remain open
(as usual, $M=\compl (H)$):
\smallskip
\begin{enumerate}
\item[(a)] Does the existence of a $\partial_R$-connected cut system for $M$
imply that $H$ is not $(3)_L$-knotted?

%\smallskip

\item[(b)] Does the existence of a cut system for $M$ imply that $H$ is not
$(4)_L$-knotted?
\end{enumerate}
\smallskip

Taking into account Remark~\ref{poincrem}, we expect that,
in order to answer these questions, one should probably invoke
very deep results. 
We get a potentially easier weaker version of these
questions by allowing us to work {\it up to reimbedding} of $M$.

\medskip

We observe
that the quandle obstructions we have constructed have 
proved to be sufficiently effective for
producing some ``ad hoc'' examples pertinent to our
discussion. On the other hand, they are by themselves very weak, and
the quandle invariants machinery is potentially much more
powerful. For example, the basic dihedral quandle invariants we have
used are not able to distinguish the Kinoshita's $M_K$ from the
unknotted handlebody.  In \cite{ishii2} this is done by means of more
sophisticated ``twisted'' quandle invariants that involve the
tetrahedral quandle and a suitable {\it quandle cocycle} on it. Then
the following is a natural challenge in order to test the performances
of quandle invariants theory:

\smallskip

{\it Is the whole program outlined in the present paper achievable
  only by means of sufficiently sophisticated instances of quandle
  invariants (recovering by the way also the information derived via
  the analysis handlebody patterns and the Alexander module obstructions)?}

\medskip

We conclude with a short outline about how to extend the above discussion 
to $H$ of {\it arbitrary genus} $g\geq 2$.

\begin{itemize}

\item The instances of knotting of $H$ can be straightforwardly
  generalized in terms of multi--hand--cuff spines, that are spines
  $\Gamma$ containing a maximal 3--valent open sub--tree (playing the
  r\^ole of the isthmus) and carrying a $g$--component constituent
  link $L_\Gamma$.

\item The definitions concerning the cut systems of $M$ are formally
  the same, provided that they are formed by $g$ disjoint surfaces.

\item A weaker version of Proposition \ref{1-K-char} holds, depending
on the weaker conclusions of Theorem 4 in \cite{scharlemann4}, when $g>2$.

\item  The use of the quandle obstructions can be adapted, obtaining the
same conclusions of Proposition \ref{K2}.

\item The definitions and the results of Section \ref{more-bpm} can be adapted
to the general case (but observe that the analysis of
  the longitudes is more complicated when $g>2$).

\item The discussion about the Alexander module obstructions extends,
  by realizing that the relevant elementary ideal is now $E_g(M)$.  In
  this way we get results analogous to
  Theorem~\ref{nonempty} (implying that every istance of knotting is non-empty)
and Proposition~\ref{via-A-obs} (implying that there exist handlebodies with a trivial constituent link which do not
admit $\partial$-connected cut-systems).

\item The conclusion of Theorem~\ref{ex<=>in:teo} holds in general,
and the proof essentially makes use of the full Poincar\'e conjecture.

\item So far we have extended the discussion focusing on the
  alternative for $M$ of {\it being or not of maximal cut number} equal to
  $g$. A more demanding problem (of the same type) would arise by
  filtering the $M$'s via the cut number values, pointing out a
  corresponding filtering of the instances of knotting of $H$.

\end{itemize}

\end{document}